%% file: main.tex
\documentclass[11pt,hidelinks,a4paper]{article}

\input{./aux/preamble}

\input{./aux/custom_commands}
\input{./aux/macros}

\author{Malte Londschien${}^{1,2}$ and Peter B\"uhlmann${}^{1}$\\
\vspace{0.1cm}\\
{\small${}^{1}$Seminar for Statistics, ETH Z\"urich, Switzerland}\\
{\small${}^{2}$AI Center, ETH Z\"urich, Switzerland}}
\title{Weak-instrument-robust subvector inference in instrumental variables regression: A subvector Lagrange multiplier test and properties of subvector Anderson-Rubin confidence sets}
\begin{document}

\date{March 2026}

\maketitle

\input{abstract.tex}

\normalsize

\vspace{0.2cm}
\small	
\noindent{Accepted for publication in the Journal of Econometrics.}
\normalsize

\input{sec_intro.tex}

\input{sec_main_contributions.tex}

\input{sec_numerical_analysis.tex}

\input{sec_applications.tex}
\input{sec_conclusion.tex}

\FloatBarrier

\input{acknowledgements.tex}
\bibliography{bib}

\appendix
\part*{Appendix}
\red{Appendix \ref{sec:exogenous_variables} extends the main results of the manuscript to provide an explicit treatment of included exogenous variables.
Appendix \ref{sec:proofs} contains the proofs of all results.
Finally, Appendix \ref{sec:additional_tables} contains supplementary tables for the numerical analysis and applications.}
\clearpage

\input{sec_exogenous_variables.tex}

\input{sec_auxiliary_results.tex}

\FloatBarrier
\section{Additional tables}
\label{sec:additional_tables}
\input{tables/table_guggenberger12_empirical_sizes_n50.tex}
\input{tables/table_guggenberger12_empirical_sizes_n100.tex}
\input{tables/table_card_black.tex}
\input{supplement.tex}

\clearpage

\end{document}

%% file: aux/preamble.tex
\usepackage[utf8]{inputenc}
\usepackage[english]{babel}
\usepackage[T1]{fontenc}

\usepackage{mathtools}
\usepackage{amsmath, amsthm, amssymb,amsfonts}
\usepackage{leftindex}
\usepackage{thmtools,thm-restate}
\usepackage{proof-at-the-end}
\PassOptionsToPackage{hidelinks}{hyperref}
\usepackage[hidelinks]{hyperref}
\usepackage[nameinlink,noabbrev]{cleveref}
\usepackage{placeins}
\usepackage{graphicx}
\usepackage{adjustbox}
\usepackage{float}
\usepackage{subcaption}
\usepackage[skip=2pt]{caption}
\usepackage[shortlabels]{enumitem}

\usepackage{tikz}
\usetikzlibrary{arrows.meta, positioning}

\usepackage{dsfont}
\usepackage{physics}
\usepackage{mleftright} %

\usepackage{booktabs}

\usepackage{natbib} 
\usepackage[a4paper, total={16cm, 22cm}]{geometry}

\newtheorem{model}{Model}

\newtheorem{definition}{Definition}

\usepackage{aliascnt}

\newaliascnt{appendixsection}{section}
\crefname{appendixsection}{Appendix}{Appendices}
\Crefname{appendixsection}{Appendix}{Appendices}

\newaliascnt{appendixsubsection}{subsection}
\crefname{appendixsubsection}{Appendix}{Appendices}
\Crefname{appendixsubsection}{Appendix}{Appendices}

%% file: aux/custom_commands.tex
\newcommand{\BP}{{\mathbb{P}}}

\newcommand{\BR}{{\mathbb{R}}}

\newcommand{\CN}{{\cal N}}
\newcommand{\CO}{{\cal O}}

\newcommand{\mX}{{{m_x}}}
\newcommand{\mW}{{{m_w}}}
\newcommand{\mC}{{{m_c}}}
\newcommand{\mD}{{{m_d}}}

\DeclareMathOperator{\Var}{Var} 
\DeclareMathOperator{\Cov}{Cov}

\DeclareMathOperator{\CI}{CI}

\DeclareMathOperator{\AR}{AR}
\DeclareMathOperator{\LM}{LM}
\DeclareMathOperator{\LR}{LR}

\DeclareMathOperator{\Wald}{Wald}

\DeclareMathOperator*{\argmin}{argmin}

\DeclareMathSymbol{\shortminus}{\mathbin}{AMSa}{"39}

\DeclareMathOperator*{\plim}{plim}
\DeclareMathOperator{\vecop}{vec}

\newcommand{\lambdamin}[1]{{\lambda_\mathrm{min}\mleft(#1\mright)}}

\def\liml{\mathrm{LIML}}

\DeclarePairedDelimiterX{\infdivx}[2]{(}{)}{%
  #1\;\delimsize\|\;#2%
}

\newcommand{\Id}{\mathrm{Id}}
\newcommand{\tod}{\overset{d}{\to}}
\newcommand{\toP}{\overset{\BP}{\to}}

\let\phi\varphi
\let\theta\vartheta
\let\epsilon\varepsilon

\let\leq\leqslant
\let\geq\geqslant

\def\iid{i.i.d.\ }
\def\kclass{\mathrm{k}}
\newcommand\numberthis{\addtocounter{equation}{1}\tag{\theequation}}

\newcommand\red[1]{{#1}}

\pgfkeys{/prAtEnd/malte/.style={normal}}
\pgfkeys{/prAtEnd/malte_all_end/.style={normal}}
\pgfkeys{/prAtEnd/malte_intro/.style={end, restate, one big link, one big link translated={\noindent See proof on page}, text proof={Proof}}}

%% file: aux/macros.tex
\def\figurekleibergencaption#1#2{Empirical maximal rejection frequency for $\alpha=0.05$ over 25 linearly spaced values of $\tau$ in $[0, \pi)$ for various $\lambda_1, \lambda_2$. This uses the data-generating process proposed by \citet{kleibergen2021efficient} with #1 and #2, based on 2500 simulations per $(\lambda_1, \lambda_2, \tau)$. As we are taking the maximum over $\tau$, values up to around $0.06$ are expected for size-correct tests. Blue indicates conservativeness, green indicates near-nominal level, yellow indicates slightly above nominal level.
AR (GKM) uses the critical values of \citet{guggenberger2019more}.}

\def\tableguggenbergerempiricalsizes#1{
    Empirical sizes of various tests at levels $\alpha = 0.01$ and $0.05$.
    Based on 10'000 draws from the data-generating process proposed by \citet{guggenberger2012asymptotic} with #1.
    LM (LIML) is the subvector test obtained by plugging the LIML into the Lagrange multiplier test.
    AR (GKM) uses the critical values of \citet{guggenberger2019more}.
}

%% file: abstract.tex
\begin{abstract}
    \noindent We propose a weak-instrument-robust subvector Lagrange multiplier test for instrumental variables regression.
    We show that it is asymptotically size-correct under a technical condition or as the number of instruments grows to infinity.
    This is the first weak-instrument-robust subvector test for instrumental variables regression to recover the degrees of freedom of the commonly used non-weak-instrument-robust Wald test.
    Additionally, we provide a closed-form solution for subvector confidence sets obtained by inverting the subvector Anderson-Rubin test.
    We show that they are centered around a k-class estimator.
    We show that the subvector confidence sets for single coefficients of the causal parameter are jointly bounded if and only if Anderson's likelihood-ratio \red{test} rejects the null hypothesis that the first-stage regression parameter is of reduced rank, that is, that the causal parameter is not identified.
    Finally, we show that if a confidence set obtained by inverting the Anderson-Rubin test is bounded and nonempty, it is equal to a Wald-based confidence set with a data-dependent confidence level.
    We explicitly compute this Wald-based confidence \red{set and its confidence level}.
\end{abstract}

%% file: sec_intro.tex
\section{Introduction}
Instrumental variables regression is essential for empirical economics.
Identifying the causal effect requires at least as many instruments as endogenous (confounded) variables.
Finding valid instruments is difficult in practice, so historically much research considers the causal effect of a single endogenous variable, instrumented by another.
However, the trend is towards the availability of larger and more detailed data, allowing modeling settings with more endogenous variables.

For causal effect estimation to influence policy, uncertainty quantification is essential.
Typically, tests based on the estimator's asymptotic normality are used to construct confidence sets and $p$-values.
\citet{staiger1997instrumental} note that in many applications, the signal-to-noise ratio in the first stage is low, leading to size distortion of such tests.
They prove that, for a single instrument and endogenous variable, if the expected first-stage $F$-statistic exceeds 10, asymptotically, Wald-based 95\% confidence sets have at least 85\% coverage.
This resulted in the following widely adopted heuristic: If the first-stage $F$-statistic is greater than 10, the instrument is considered strong, and Wald-based confidence sets are considered valid.

We believe such pre-testing is flawed.
Even though the heuristic is applied more generally, \citeauthor{staiger1997instrumental}'s \citeyearpar{staiger1997instrumental} results for the Wald test's size do not apply for multiple  instruments or multiple endogenous variables \citep{stock2002testing}.
In a recent paper, \citet{lee2022valid} show that that the first stage $F$-statistic rule of 10 can still lead to severe size distortion of the TSLS Wald test and that a larger cutoff is needed for size control.
On the other hand, pre-testing might be too conservative, preventing inference given weak instruments but a strong causal effect.

An alternative is using weak-instrument-robust tests.
These include the Anderson-Rubin test \citep{anderson1949estimation}, \citeauthor{kleibergen2002pivotal}'s \citeyearpar{kleibergen2002pivotal} Lagrange multiplier test, and \citeauthor{moreira2003conditional}'s \citeyearpar{moreira2003conditional} conditional likelihood-ratio test.
They have the correct size asymptotically, even if the signal-to-noise ratio in the first stage decreases such that the first-stage $F$-statistic is of constant order (weak instrument asymptotics).

For multiple endogenous variables, multivariate confidence sets for the causal parameter are difficult to interpret.
Subvector tests and confidence sets for individual coefficients of the causal parameter, similar to $t$-tests in linear regression, are more useful.
A naive approach to construct such subvector confidence sets is by projecting the multivariate confidence set onto the coefficient of interest \citep{dufour2005projection}.
However, this results in a loss of power.

Subvector variants of the Anderson-Rubin and conditional likelihood-ratio tests have been proposed by \citet{guggenberger2012asymptotic} and \citet{kleibergen2021efficient}.
However, with more instruments than endogenous variables, the limiting distributions have higher degrees of freedom than the corresponding Wald test.
Consequently, additional instruments, for example by modeling non-linear relations with interactions and non-linear basis expansions, possibly increase the size of the confidence sets.
This is undesirable, as the additional instruments should improve identification.

We propose a subvector variant of the Lagrange multiplier test.
This recovers the degrees of freedom of the Wald test, independently of the number of instruments.
While the test is not necessarily more powerful than the subvector conditional likelihood-ratio test, to our knowledge, this is the first weak-instrument-robust subvector test in instrumental variables regression to achieve this.
Our test fills a gap in the weak-instrument-robust testing literature for instrumental variables regression, as shown in \Cref{tab:overview_robust_tests}.
\input{tables/table_overview_robust_tests.tex}

\citet{dufour2005projection} note that the confidence sets obtained by inverting the Anderson-Rubin test are quadrics, as are their projections.
We provide a closed-form solution for the inverse Anderson-Rubin test confidence sets using the more powerful critical values by \citet{guggenberger2012asymptotic} and show that they are centered around a k-class estimator.
Additionally, we prove that the subvector confidence sets for a single coefficient of the causal parameter are jointly bounded if and only if \citeauthor{anderson1951estimating}'s \citeyearpar{anderson1951estimating} likelihood-ratio test rejects the hypothesis that the first-stage regression parameter is of reduced rank, that is, that the causal parameter is not identified.
Finally, we show that if the confidence set obtained by inverting the Anderson-Rubin test is bounded and nonempty, it is equal to a Wald-based confidence set around a k-class estimator.
We explicitly compute this k-class estimator's $\kappa$ parameter and the Wald confidence set's confidence level.

Implementations for k-class estimators, the (subvector) Anderson-Rubin, (conditional) like\-li\-hood-ratio, Lagrange multiplier, and Wald tests, and the confidence sets obtained by inversion can be found in our open-source software package \href{https://github.com/mlondschien/ivmodels/}{\texttt{ivmodels}} for Python.
We hope that this simplifies and thus possibly popularizes the use of weak-instrument-robust tests in empirical research.
The \texttt{ivmodels} software package is available on \href{https://pypi.org/project/ivmodels/}{\texttt{PyPI}} and \href{https://anaconda.org/conda-forge/ivmodels}{\texttt{conda-forge}}.
See the GitHub repository at \href{https://github.com/mlondschien/ivmodels/}{\texttt{github.com/mlondschien/ivmodels}} and the documentation at \href{https://ivmodels.readthedocs.io/}{\texttt{ivmodels.readthedocs.io}} for more details.
Code and instructions to reproduce the tables and figures in this paper are available at \href{https://github.com/mlondschien/ivmodels-simulations}{\texttt{github.com/mlondschien/ivmodels-simulations}}.

%% file: tables/table_overview_robust_tests.tex
\begin{table}[H]
\begin{tabular}{l | l | l | l}
    & Anderson-Rubin & Lagrange multiplier & conditional likelihood-ratio \\
    \hline
    complete & \citet{anderson1949estimation} & \citet{kleibergen2002pivotal} & \citet{moreira2003conditional} \\
    subvector & \citet{guggenberger2012asymptotic} & \bf this work & \citet{kleibergen2021efficient}
\end{tabular}
\caption{
    \label{tab:overview_robust_tests}    
    Overview of weak-instrument-robust tests for instrumental variables regression.
}
\end{table}

%% file: sec_main_contributions.tex
\section{Main results}
\label{sec:main_results}
We present our main contributions in \Cref{sec:main_results:subvector_lagrange_multiplier,sec:main_results:closed_form_subvector_confidence_sets}.
See \citep{londschien2025overview} for an overview of the literature on weak-instrument-robust inference for instrumental variables regression.

\def\consideralinearinstrumental{Consider a linear instrumental variables regression model with multiple endogenous covariates.
In a typical application, one wishes to make (subvector) inferences for each component of the causal parameter.
To model this, we split the endogenous covariates and the components of the causal parameter into two parts, formalized in \Cref{model:1}.
Our parameter of interest is $\beta_0$ and we treat $\gamma_0, \Pi_X,$ and $\Pi_W$ as nuisance parameters.}
\consideralinearinstrumental
\input{theorems/model_1.tex}

In practice, it is common that additional exogenous covariates $C$ enter the model equation.
If these are not of interest, one can reduce to \cref{model:1} by replacing $y$, $X$, $W$, and $Z$ with their residuals after regressing out $C$.
This also allows for the inclusion of an intercept by centering the data.
If the exogenous covariates $C$ are of interest, that is, one would like to make inference for their causal effect on $y$, one can include them into the model as both instruments and additional endogenous covariates.
We treat this explicitly in Appendix \ref{sec:exogenous_variables}.

In standard asymptotics, the expected first-stage $F$-statistic (or its multivariate extension, see \citeauthor{anderson1951estimating}, \citeyear{anderson1951estimating}) grows with the sample size.
\citet{staiger1997instrumental} note that in many applications, even though the sample size $n$ is large, the $F$-statistic for regressing $X$ on $Z$ is small.
In these settings, standard asymptotics no longer apply, motivating weak instrument asymptotics, where the first-stage $F$-statistic is of constant order.

We assume that a central limit theorem applies to the sums $Z^T \varepsilon$, $Z^T V_X$, and $Z^T V_W$.
\input{theorems/assumption_1.tex}
\Cref{ass:1} is standard in the literature.
It is equivalent to the assumption made by \citet{kleibergen2002pivotal} and \citet{staiger1997instrumental} and is implied by assumptions made by \citet{guggenberger2012asymptotic}, \citet{guggenberger2019more}, and \citet{kleibergen2021efficient}.
If the $(Z_i, \varepsilon_i, V_{X, i}, V_{W, i})$ are i.i.d.\ with finite second moments and if the noise terms $\varepsilon, V_{X, i}$, and $V_{W, i}$ are homoscedastic, centered, and uncorrelated with the $Z_i$, then \Cref{ass:1} holds \citep[][Lemma 1]{londschien2025overview}.
\red{Crucially, our results do not allow for general forms of heteroskedasticity.
For inference that is robust to heteroskedasticity, see \citep{guggenberger2024powerful}.}

\subsection{A subvector Lagrange multiplier test}
\label{sec:main_results:subvector_lagrange_multiplier}
We propose a subvector extension of the Lagrange multiplier test statistic from \citet{kleibergen2002pivotal}.
For any $A \in \BR^{p \times q}$, let $P_A := A (A^T A)^\dagger A^T$, where $\dagger$ denotes the Moore–Penrose inverse, be the projection matrix onto the column span of $A$.
Write $M_A := \Id_{p} - P_A$ for the projection onto its orthogonal complement.
\input{theorems/def_subvector_klm_test_statistic.tex}
\noindent This reduces to \citeauthor{kleibergen2002pivotal}'s \citeyearpar{kleibergen2002pivotal} Lagrange multiplier test %
if $\mW=0$.

\citet{guggenberger2012asymptotic} analyse the subvector Lagrange multiplier test statistic obtained by plugging in the \red{limited information maximum likelihood (LIML)} estimator \red{(see, e.g., Definition 6 of \citeauthor{londschien2025overview}, \citeyear{{londschien2025overview}})} $\hat\gamma_\liml$ using outcomes $y - X \beta$, covariates $W$, and instruments $Z$ for $\gamma$ into Equation \eqref{def:subvector_klm_test_statistic}.
They propose a data-generating process for which this subvector Lagrange multiplier test with $\chi^2(\mX)$ critical values is size distorted and thus show that it is not weak-instrument-robust.
\red{The LIML minimizes the Anderson-Rubin test statistic, that is, $\hat\gamma_\liml = \argmin_\gamma \AR(\beta, \gamma)$ (see, e.g., Corollary 9 of \citeauthor{londschien2025overview}, \citeyear{londschien2025overview}) and \citet{guggenberger2012asymptotic} show that plugging it into the Anderson-Rubin test statistic yields a size-correct test with $\chi^2(k - \mW)$ critical values.}
\red{However, the minimization over $\gamma$ in \eqref{def:subvector_klm_test_statistic}} is not equivalent to plugging in the LIML.
We show that our proposed subvector Lagrange multiplier test has the correct size under a technical condition.

\input{theorems/prop_subvector_klm_test_statistic_chi_squared_statement.tex}

The proof of \Cref{prop:subvector_klm_test_statistic_chi_squared} uses ideas from \citet{kleibergen2002pivotal} and \citet{guggenberger2012asymptotic}.
We provide some intuition below.
For the full proof see the proof of \Cref{prop:subvector_klm_test_statistic_chi_squared_exogenous} in Appendix \ref{app:proofs:subvector_lagrange_multiplier}.
This extends \Cref{prop:subvector_klm_test_statistic_chi_squared} to the case with included exogenous variables.

Under weak instrument asymptotics, \citeauthor{kleibergen2002pivotal}'s \citeyearpar{kleibergen2002pivotal} Lagrange multiplier test statistic can be shown to be equal to
\begin{equation*}
    \LM(\beta) = \frac{1}{\hat\sigma_\varepsilon^2(\beta)} \| P_{\Psi_{\tilde X(\beta)}}  (\Psi_\varepsilon + \Psi_{X} (\beta - \beta_0) ) \|^2,
\end{equation*}
where $\Psi_\varepsilon \sim \CN(0, \sigma^2_\varepsilon \cdot \Id_k), \ \Psi_X \sim \CN(Q^{1/2} \Pi, \Omega_{V} \otimes \Id_k),$ and $\Psi_{\tilde X(\beta)} \sim \CN(\mu(\beta), \tilde\Omega(\beta) \otimes \Id_k)
$ \red{for some $\tilde \Omega(\beta)$,}
are asymptotically jointly Gaussian.
Importantly, $\tilde X(\beta)$, \red{equal to $\tilde S(\beta)$ from \Cref{def:subvector_klm_test_statistic} when $\mW=0$,} is constructed such that $\Psi_\varepsilon + \Psi_{X} (\beta_0 - \beta)$ and \red{$\Psi_{\tilde X(\beta)} := (Z^T Z)^{-1/2} Z^T \tilde X(\beta)$} are asymptotically independent.
Thus, under the null $\beta = \beta_0$, conditionally on $\Psi_{\tilde X(\beta_0)}$, we have that $\| P_{\Psi_{\tilde X(\beta_0)}} \Psi_\varepsilon \|^2 \to_\BP \sigma_\varepsilon^2 \chi^2(m)$, where $m = \rank(P_{\Psi_{\tilde X(\beta_0)}})$.
Integrating over $\Psi_{\tilde X(\beta_0)}$, this holds also unconditionally and $\LM(\beta_0) \to_d \chi^2(m)$ follows from $\hat\sigma_\varepsilon^2(\beta_0) \to_\BP \sigma_\varepsilon^2$.

In contrast, \citet{guggenberger2012asymptotic} construct $\gamma^\dagger$ in such a manner that
\begin{align*}
    (k - \mW) \AR(\beta_0, \gamma^\dagger) &= \frac{1}{\hat\sigma_\varepsilon^2(\beta_0, \gamma^\dagger)} \| \Psi_\varepsilon + \Psi_W (\gamma_0 - \gamma^\dagger) \|^2 
    = \frac{\mathrm{const(\gamma^\dagger)}}{\hat\sigma_\varepsilon^2(\beta_0, \gamma^\dagger)} \| M_{\Psi_{\tilde W(\beta_0, \gamma_0)}} (\Psi_\varepsilon  ) \|^2,
\end{align*}
where $\Psi_\varepsilon \sim \CN(0, \sigma^2_\varepsilon \cdot \Id_k)$ \red{ and $\Psi_{\tilde W(\beta_0, \gamma_0)} := (Z^T Z)^{-1/2} Z^T \tilde W(\beta_0, \gamma_0) \sim \CN(\mu, \tilde\Omega \otimes \Id_k)$} are asymptotically jointly Gaussian and independently distributed, \red{for $(\tilde X(\beta, \gamma), \tilde W(\beta, \gamma) ) := \tilde S(\beta, \gamma)$ and some invertible $\tilde \Omega$}.
Thus, under the null $\beta = \beta_0$, conditionally on $\Psi_{\tilde W(\beta_0, \gamma_0)}$ we have that $\| M_{\Psi_{\tilde W(\beta_0, \gamma_0)}} \Psi_\varepsilon \|^2 \to_\BP \sigma_\varepsilon^2 \chi^2(k - \mW)$, where $k - \mW = \rank(M_{\Psi_{\tilde W(\beta_0, \gamma_0)}}) = \rank(\Id_k - P_{\Psi_{\tilde W(\beta_0, \gamma_0)}})$.
This then also holds unconditionally.
The subvector Anderson-Rubin test statistic is then bounded from above by a $\chi^2(k - \mW) / (k - \mW)$ random variable, as $\plim \hat\sigma_\varepsilon^2(\beta, \gamma^\dagger) / \mathrm{const(\gamma^\dagger)} \geq \sigma_\varepsilon^2$ and $\min_\gamma \AR(\beta_0, \gamma) \leq \AR(\beta_0, \gamma^\dagger)$.

Note that the $\gamma^\dagger$ is constructed such that $\Psi_\varepsilon + \Psi_W (\gamma_0 - \gamma^\dagger) = M_{\Psi_{\tilde W(\beta_0, \gamma_0)}} \Psi_\varepsilon \neq M_{\Psi_{\tilde W(\beta_0, \gamma^\dagger)}} \Psi_\varepsilon$.
The random variable $\Psi_{\tilde W(\beta_0, \gamma^\dagger)}$ is asymptotically independent of $\Psi_\varepsilon + \Psi_{W} (\gamma_0 - \gamma^\dagger)$, but not of $\Psi_\varepsilon$, and so, \red{for the subvector Lagrange multiplier test}, the argument based on conditioning on $\Psi_{\tilde W(\beta_0, \gamma^\dagger)}$ does not work.

The \red{subvector} Lagrange multiplier test statistic is equal to $\min_\gamma \LM(\beta, \gamma)$, where $$\LM(\beta, \gamma) = \frac{1}{\hat\sigma_\varepsilon^2(\beta, \gamma)} \| P_{[\Psi_{\tilde X(\beta, \gamma)}, \Psi_{\tilde W(\beta, \gamma)}]} (\Psi_\varepsilon + \Psi_{X} (\beta - \beta_0) + \Psi_W (\gamma_0 - \gamma)) \|^2.$$
If we plug in $\gamma^\dagger$ as above, we obtain $\LM(\beta_0, \gamma^\dagger) = \frac{1}{\hat\sigma_\varepsilon^2(\beta_0, \gamma^\dagger)} \| P_{[\Psi_{\tilde X(\beta_0, \gamma^\dagger)}, \Psi_{\tilde W(\beta_0, \gamma^\dagger)}]} M_{\Psi_{\tilde W(\beta_0, \gamma_0)}} \Psi_\varepsilon \|^2$.
The first projection is onto the column span of $[\Psi_{\tilde X(\beta_0, \gamma^\dagger)}, \Psi_{\tilde X(\beta_0, \gamma^\dagger)}]$, which is not independent of $\Psi_\varepsilon$.
\Cref{tc:subvector_klm} ensures the existence of a $\gamma^\star$ such that 
$$\LM(\beta_0, \gamma^\star) \toP \frac{1}{\hat\sigma_\varepsilon^2(\beta_0, \gamma^\star)} \| P_{M_{Q^{1/2} \Pi_W} [\Psi_{\tilde X(\beta_0, \gamma^\star)}, \Psi_{\tilde W(\beta_0, \gamma^\star)}]} (\Psi_\varepsilon + \Psi_{V_W} (\gamma_0 - \gamma^\star)) \|^2,$$
where $\Psi_\varepsilon^\star := \Psi_\varepsilon + \Psi_{V_W} (\gamma_0 - \gamma^\star) \sim \CN(0, \sigma_{\varepsilon^\star}^2 \cdot \Id_k)$ is independent of $[\Psi_{\tilde X(\beta_0, \gamma^\star)}, \Psi_{\tilde W(\beta_0, \gamma^\star)}]$.
We can thus apply the conditioning argument as above to show that
$$\LM(\beta_0, \gamma^\star) \tod \chi^2(\mX), \text{ where }\mX = m - \mW = \rank(P_{M_{Q^{1/2} \Pi_W} [\Psi_{\tilde X(\beta_0, \gamma^\star)}, \Psi_{\tilde W(\beta_0, \gamma^\star)}]}).$$
The result then follows from $\hat\sigma_\varepsilon^2(\beta_0, \gamma^\star) \to_\BP \sigma_{\varepsilon^\star}^2$ and $\min_\gamma \LM(\beta_0, \gamma) \leq \LM(\beta_0, \gamma^\star)$.

We were unable to prove that the \Cref{tc:subvector_klm} holds in general with high probability and leave this as future work.
However, we empirically observe that the statement of \Cref{prop:subvector_klm_test_statistic_chi_squared} holds in great generality for finite sample sizes and fixed number of instruments $k$, see also the numerical analyses in \Cref{sec:numerical_analysis}.
We thus conjecture the following:
\input{theorems/conj_subvector_klm_test_statistic.tex}

Recall the issue with using $\gamma^\dagger$ as in \citep{guggenberger2012asymptotic}: $\Psi_{\tilde W(\beta_0, \gamma^\dagger)}$ is independent of $\Psi_\varepsilon + \Psi_W (\gamma_0 - \gamma^\dagger)$, but not of $\Psi_\varepsilon$.
However, we can show that as $k \to \infty$, we have that $\Psi_{\tilde W(\beta_0, \gamma^\dagger)} \to_\BP \Psi_{\tilde W(\beta_0, \gamma_0)}$.
This yields the following result.

\input{theorems/prop_subvector_klm_test_statistic_chi_squared_many_instruments.tex}

\noindent For the proof of \Cref{prop:subvector_klm_test_statistic_chi_squared_many_instruments}, see the proof of \Cref{prop:subvector_klm_test_statistic_chi_squared_many_instruments_exogenous} in Appendix \ref{app:proofs:subvector_lagrange_multiplier}.
\Cref{prop:subvector_klm_test_statistic_chi_squared_many_instruments_exogenous} extends \Cref{prop:subvector_klm_test_statistic_chi_squared_many_instruments} to the case with included exogenous variables.

Note that \Cref{ass:1} as $k = k(n) \to \infty$ is a stronger assumption than \Cref{ass:1}.
Lemma 1 of \citet{londschien2025overview} can be extended to allow for a $k = k(n)$ that slowly grows with $n$, for example $k(n) = \log(n)$, by invoking a central limit theorem in slowly growing dimension, see also \citet{bentkus2003dependence}.
We note that an arbitrarily slow growth of $k = k(n)$ in \Cref{ass:1} is sufficient for \Cref{prop:subvector_klm_test_statistic_chi_squared_many_instruments}.

\Cref{prop:subvector_klm_test_statistic_chi_squared_many_instruments,prop:subvector_klm_test_statistic_chi_squared} directly imply the following corollary:
\input{theorems/cor_subvector_lagrange_multiplier_test.tex}

Thus, under \Cref{tc:subvector_klm} or if $k = k(n) \to \infty$, the weak-instrument-robust subvector Lagrange multiplier test achieves the same degrees of freedom as the commonly used Wald test that is not robust to weak instruments.
To our knowledge, this is the ﬁrst weak-instrument-robust subvector test for the causal parameter in instrumental variables regression to achieve this.
The subvector Lagrange multiplier test fills a gap in the literature on weak-instrument-robust inference for instrumental variables regression, as shown in \Cref{tab:overview_robust_tests}.

We numerically analyze the size (and power) of various tests under (a variant of) \citeauthor{guggenberger2012asymptotic}'s (\citeyear{guggenberger2012asymptotic}) data-generating process in \Cref{sec:numerical_analysis:guggenberger12_size,sec:numerical_analysis:guggenberger12_power}.
Our proposed subvector Lagrange multiplier test is size-correct and substantially more powerful than their subvector Anderson-Rubin test, which uses $\chi^2(k - \mW)$ critical values.
We compare properties of various tests in instrumental variables regression in \Cref{tab:comparison_tests}.
\input{tables/table_comparison_tests.tex}

In practice, finding valid instruments is challenging and applications where the number of instruments is much larger than the number of endogenous covariates are rare.
However, we encourage improving identification by modeling non-linear relations between instruments and endogenous covariates with interactions and non-linear basis expansions.
If using the non-weak-instrument-robust Wald test, the addition of such instruments possibly decreases the concentration parameter $\mu^2 / k := n \cdot \lambdamin{\Omega_V^{-1} \Pi^T Q \Pi} / k$ \citep{staiger1997instrumental,stock2002testing},
degrading the the test's approximation quality and possibly leading to size distortion.
If using the weak-instrument-robust Anderson-Rubin test, the critical values increase with the number of instruments and additional instruments possibly increase the size of confidence intervals.
These issues do not affect the subvector Lagrange multiplier test, which is robust to weak instruments and whose critical values do not depend on the number of instruments, as described by \Cref{cor:subvector_lagrange_multiplier_test}.

In \Cref{sec:applications}, we apply various tests to data from \citet{card1993using} and \citet{tanaka2010risk}.
The addition of instrument interactions generally leads to smaller $p$-values for the Lagrange multiplier test, which is not the case for the Anderson-Rubin and conditional likelihood-ratio tests.
\red{As a score-type test, the Lagrange multiplier test statistic is zero at all stationary points of the likelihood or equivalently of the Anderson-Rubin test statistic.
As pointed out by a reviewer and established by \citet{andrews2006optimal} and \citet{kleibergen2007generalizing}, there are $m + 1$ such stationary points.
In practice, this results in reduced power around such stationary points, see also \Cref{sec:numerical_analysis:guggenberger12_power}.
Confidence sets obtained by inversion of the subvector Lagrange multiplier test are possibly disjoint, with the number of components equal to the number of small eigenvalues of the concentration matrix (the number of weakly identified endogenous variables) plus one.
We observe this phenomenon in \Cref{sec:applications}.
}

Even if there is only one endogenous covariate, subvector inference is needed to construct weak-instrument-robust confidence sets for included exogenous covariates.
In a typical instrumental variables specification, inference is needed for the causal effect of endogenous covariates after accounting for included exogenous covariates.
However, there are also use cases where the focus is on the causal effect of an exogenous covariate after accounting for an endogenous covariate, and instruments are required to identify the nuisance parameter \citep{thams2022identifying}.
For example, \citet{tanaka2010risk} discuss the effect of the effect of (exogenous) gender on risk preferences, after accounting for (endogenous) income.

\citet{kleibergen2021efficient} suggest to test a hypothesis on the causal effect of an included exogenous covariate with their subvector conditional likelihood-ratio test by including it both as an endogenous covariate and an instrument.
However, both their proof and that of the full vector conditional likelihood-ratio test \citep{moreira2003conditional} depend on invertibility of the covariance estimate $\hat\Omega = \frac{1}{n - k} X^T M_Z X$.
This is not the case if an included exogenous covariate is included into both $X$ and $Z$.
\citet[][Proposition 54]{londschien2025overview} explicitly computes the limiting distribution of the full vector conditional likelihood-ratio test statistic with included exogenous covariates.
Their limiting distribution is strictly smaller than the one that would be obtained by simply including the included exogenous covariates as both endogenous covariates and instruments.
That is, testing hypotheses on the causal effect of included exogenous covariates by simply including them into both instruments and endogenous variables and applying \citeauthor{moreira2003conditional}'s (\citeyear{moreira2003conditional}) conditional likelihood-ratio test is conservative and possibly leads to unnecessarily large confidence sets.
It is not clear whether \citeauthor{londschien2025overview}'s (\citeyear{londschien2025overview}) result on the distribution of the conditional likelihood-ratio test statistic with included exogenous covariates generalizes to \citeauthor{kleibergen2021efficient}'s \citeyearpar{kleibergen2021efficient} subvector test.

In contrast, we treat the setting with included exogenous covariates explicitly and prove versions of \Cref{prop:subvector_klm_test_statistic_chi_squared} and \Cref{prop:subvector_klm_test_statistic_chi_squared_many_instruments} that allow for included exogenous variables under test in \Cref{prop:subvector_klm_test_statistic_chi_squared_exogenous} and \Cref{prop:subvector_klm_test_statistic_chi_squared_many_instruments_exogenous} in Appendix \ref{sec:exogenous_variables}.

The optimization problem in \Cref{def:subvector_klm_test_statistic} is non-convex and thus difficult to solve.
In our \texttt{ivmodels} software package, which was used for the simulations in \Cref{sec:numerical_analysis}, we use the quasi-Newton method of Broyden–Fletcher–Goldfarb–Shannon (BFGS).
This is the default option of \texttt{scipy.optimize.minimize} in Python.
While this is not guaranteed to find the global optimum, in practice we observe that it suffices to find a local optimum next to the LIML for the test to be size-correct.
See Appendix \ref{sec:optimization} for details.

\subsection{Closed-form subvector confidence sets obtained by inverting the Anderson-Rubin test}
\label{sec:main_results:closed_form_subvector_confidence_sets}
The Anderson-Rubin test is an important test for instrumental variables regression.
\citet{anderson1949estimation} showed that the statistic is asymptotically $\chi^2(k) / k$ distributed, independently of instrument strength.
\citet{guggenberger2012asymptotic} proposed the following subvector variant.
\input{theorems/def_subvector_anderson_rubin_test_statistic.tex}
\citet{guggenberger2012asymptotic} showed that under \Cref{model:1} and \Cref{ass:1}, the subvector Anderson-Rubin test statistic is asymptotically bounded from above by a ${\chi^2(k - \mW) / (k - \mW)}$ distributed random variable, also independently of instrument strength.
This yields weak-instrument-robust confidence sets for subvectors of the causal parameter via test inversion:
$$
\CI_{\AR}(1 - \alpha) := \{
    \beta \in \BR^{\mX} \mid (k - \mW) \cdot \AR(\beta) \leq F^{-1}_{\chi^2(k - \mW)}(1 - \alpha)
\}.
$$
We give an explicit form for these confidence sets.
\input{theorems/prop_cis_closed_forms_intro.tex}
\noindent For the proof of \Cref{prop:cis_closed_forms_intro}, see the proof of \Cref{prop:cis_closed_forms} in Appendix \ref{app:proofs:closed_form_subvector_confidence_sets}.
In \Cref{prop:cis_closed_forms}, we also compute the closed-form solutions of the subvector inverse Wald and likelihood-ratio tests.
To derive conditions for the subvector inverse Anderson-Rubin test confidence sets to be bounded, we introduce the following technical condition.
\input{theorems/tc_2.tex}
\noindent If the noise in \Cref{model:1} is absolutely continuous with respect to the Lebesgue measure, then \Cref{tc:2} holds with probability one.
\input{theorems/prop_cis_bounded_intro.tex}
\noindent For the proof of \Cref{prop:cis_bounded_intro}, see the proof of \Cref{prop:cis_bounded} in Appendix \ref{app:proofs:closed_form_subvector_confidence_sets}.
\Cref{prop:cis_bounded} additionally provides a condition for the (subvector) inverse likelihood-ratio test confidence set to be bounded.
If the \Cref{tc:2} \ref{tc:2:a} does not hold, then the subvector inverse Anderson-Rubin test confidence set for $1 - \alpha = F^{-1}_{\chi^2(k - \mW)}(\lambda)$ may be bounded.
See also \Cref{fig:figure_tc3_counterexample} in \red{Appendix}.

One is typically interested in subvector inference for a single component of the causal parameter, that is, $\mX = 1$, as this is easier to interpret.
We characterize when the subvector inverse Anderson-Rubin test confidence sets for individual components are (jointly) bounded and nonempty (thus, forming confidence intervals).
\input{theorems/cor_ar_bounded_iff_rank_test_rejects_intro.tex}
\noindent For the proof of \Cref{cor:ar_bounded_iff_rank_test_rejects_intro}, see the proof of \Cref{cor:ar_bounded_iff_rank_test_rejects} in Appendix \ref{app:proofs:closed_form_subvector_confidence_sets}.
\Cref{cor:ar_bounded_iff_rank_test_rejects} additionally gives a condition for the (subvector) inverse likelihood-ratio test confidence sets to be bounded.

\red{
\Cref{cor:ar_bounded_iff_rank_test_rejects_intro} summarizes properties of the confidence sets obtained by inverting the Anderson-Rubin test that follow directly from \Cref{prop:cis_closed_forms_intro}.
We hope that this characterization assists in the interpretation of subvector inference with the Anderson-Rubin test.
As pointed out by a reviewer, these properties are implicit in the geometric construction of subvector statistics in \citet{guggenberger2012asymptotic} and \citet{kleibergen2021efficient}.
If $m = \mX = 1$, \citeauthor{anderson1951estimating}'s \citeyearpar{anderson1951estimating} likelihood-ratio test of reduced rank reduces to the $F$-test of the first-stage regression of $X$ on $Z$ using $\chi^2(k)$ critical values.
For this special case, \citet{kleibergen2007generalizing} notes that the inverse Anderson-Rubin test confidence sets are unbounded if and only if the first-stage $F$-test rejects the null hypothesis of underidentification.
For the setting with multiple endogenous covariates ($\mW > 0$), \citet{kleibergen2021efficient} shows in their Theorem 12a that for values of $\beta$ far from $\beta_0$, the subset Anderson-Rubin statistic equals \citeauthor{anderson1951estimating}'s \citeyearpar{anderson1951estimating} likelihood-ratio test for reduced rank, implying \Cref{cor:ar_bounded_iff_rank_test_rejects_intro} (b).
Property (d) is a direct consequence of the definition of the subvector Anderson-Rubin test and the critical values first derived by \citet{guggenberger2012asymptotic}.}

\Cref{cor:ar_bounded_iff_rank_test_rejects_intro} \ref{cor:ar_bounded_iff_rank_test_rejects_intro:b} only holds if $\mX = 1$, that is, if we are testing a single component of the causal parameter.
Otherwise, the condition for boundedness of the inverse Anderson-Rubin test and \citeauthor{anderson1951estimating}'s \citeyearpar{anderson1951estimating} likelihood-ratio test for reduced rank use different critical values.
See also \Cref{fig:inverse_ar_different_alpha} (middle panel), where the 2-dimensional 90\% confidence set for the causal parameter is unbounded, even though the $p$-value of the \citeauthor{anderson1951estimating}'s \citeyearpar{anderson1951estimating} likelihood-ratio test for reduced rank is 0.07.
As predicted by \Cref{cor:ar_bounded_iff_rank_test_rejects_intro}, the 1-dimensional 90\% subvector confidence sets are bounded.

\input{figures/figure_inverse_ar_different_alpha.tex}

On the other hand, the critical values $F^{-1}_{\chi^2(k-m + 1)}(1 - \alpha)$ in \Cref{cor:ar_bounded_iff_rank_test_rejects_intro} \ref{cor:ar_bounded_iff_rank_test_rejects_intro:d} do not match those of the LIML variant of the J-statistic  $F^{-1}_{\chi^2(k-m)}(1 - \alpha)$.
In particular, rejection of the LIML J-test does not imply that the inverse Anderson-Rubin confidence sets are empty.

\citet{dufour2005projection} note that if $\mW = 0$, confidence sets obtained by inverting the Anderson-Rubin test are quadrics.
They propose constructing subvector confidence sets by projection but do not give an explicit closed form for the confidence sets.
Nor do they show that the confidence sets are centered around a k-class estimator.
Additionally, they do not use the more powerful $\chi^2(k - \mW)$ critical values by \citet{guggenberger2012asymptotic}, resulting in a substantial loss of power, see \Cref{fig:inverse_ar_different_alpha}.

See \Cref{def:wald_test_statistic,def:lr_test_statistic} in Appendix \ref{app:proofs:closed_form_subvector_confidence_sets} for the definitions of the Wald and likelihood-ratio test statistics.
Let $\CI_{\Wald, \hat\beta_k(\kappa)}(1 - \alpha)$ be the $1-\alpha$ confidence set obtained by inverting the subvector Wald test centered around the k-class estimator $\hat\beta_k(\kappa)$ with parameter $\kappa$ and let $\CI_{\LR}(1 - \alpha)$ be the $1 - \alpha$ confidence set obtained by inverting the subvector likelihood-ratio test,
As for the subvector Anderson-Rubin test, these are quadrics (see \Cref{prop:cis_closed_forms} in Appendix \ref{app:proofs:closed_form_subvector_confidence_sets}).
Notably, they differ only in the critical value and the $\kappa$ parameter of the k-class estimator they are centered around.
We derive the following result:
\input{theorems/prop_inverse_ar_equal_to_wald.tex}
Let $\hat\kappa_\liml$ be the $\kappa$ parameter of the LIML's k-class parametrization ${\hat\beta_\liml = \hat\beta_\kclass(\hat\kappa_\liml)}$ \citep[see, e.g.,][Proposition 10]{londschien2025overview}.
As ${J_\liml = (n - k) (\hat\kappa_\liml - 1)}$, if ${J_\liml \leq F^{-1}_{\chi^2(k-\mW)}(1 - \alpha)}$, then ${\kappa_{\AR}(\alpha) \geq \hat\kappa_\liml}$.

%% file: theorems/model_1.tex
\begin{theoremEnd}[malte,restate command=modelone]{model}
    \label{model:1}
    Let $y_i = X_i^T \beta_0 + W_i^T \gamma_0 + \varepsilon_i \in \BR$ with $X_i = Z_i^T \Pi_X + V_{X, i} \in \BR^\mX$ and $W_i = Z_i^T \Pi_W + V_{W, i} \in \BR^\mW$ for random vectors $Z_i \in \BR^k, V_{X, i}\in \BR^\mX, V_{W, i}\in\BR^\mW$, and $\varepsilon_i \in \BR$ for $i=1\ldots, n$ and parameters $\Pi_X \in \BR^{k \times \mX}$, $\Pi_W \in \BR^{k \times \mW}$, $\beta_0 \in \BR^\mX$, and $\gamma_0 \in \BR^\mW$.
    We call the $Z_i$ \emph{instruments}, the $X_i$ \emph{endogenous covariates of interest}, the $W_i$ \emph{endogenous covariates not of interest}, and the $y_i$ \emph{outcomes}.
    The $V_{X, i}, V_{W, i}$, and $\varepsilon_i$ are \emph{errors}.
    These need not be independent across observations.
    Let $Z \in \BR^{n \times k}, X \in \BR^{n \times \mX}, W \in \BR^{n \times \mW}$, and $y \in \BR^n$ be the matrices of stacked observations.

    In \emph{strong instrument asymptotics}, we assume that $\Pi := \begin{pmatrix} \Pi_X & \Pi_W \end{pmatrix}$ is fixed and of full column rank $m := \mX + \mW$.
    In \emph{weak instrument asymptotics} \citep{staiger1997instrumental}, we assume that $\sqrt{n} \, \Pi = \sqrt{n}  \begin{pmatrix} \Pi_X & \Pi_W \end{pmatrix}$ is fixed and of full column rank $m$. Thus, $\Pi = \CO(\frac{1}{\sqrt{n}})$.
    Both asymptotics imply that $k \geq m$.
\end{theoremEnd}

%% file: theorems/assumption_1.tex
\begin{theoremEnd}[malte,restate command=assumptionone]{assumption}
    \label{ass:1}
    Let
    $$
    \Psi := \begin{pmatrix} \Psi_{\varepsilon} & \Psi_{V_X} & \Psi_{V_W} \end{pmatrix} := (Z^T Z)^{-1/2} Z^T \begin{pmatrix} \varepsilon & V_X & V_W \end{pmatrix} \in \BR^{k \times (1 + m)}.
    $$
    Assume there exist $\Omega \in \BR^{(1 + m) \times (1 + m)}$ and $Q \in \BR^{k \times k}$ positive definite such that, as $n \to \infty$,
    \begin{align*}
        &\mathrm{(a)} \ \ \frac{1}{n} \begin{pmatrix}\varepsilon & V_X & V_W\end{pmatrix}^T \begin{pmatrix}\varepsilon & V_X & V_W \end{pmatrix} \toP \Omega, \\
        &\mathrm{(b)} \ \ \vecop(\Psi) \tod \CN(0, \Omega \otimes \Id_k), \text{ and }\\
        &\mathrm{(c)} \ \ \frac{1}{n} Z^T Z \toP Q,
    \end{align*}
    where $\Cov(\vecop(\Psi)) = \Omega \otimes \Id_k$ means $\Cov(\Psi_{i, j}, \Psi_{i', j'}) = 1_{i = i'} \cdot \Omega_{j, j'}$.
\end{theoremEnd}

%% file: theorems/def_subvector_klm_test_statistic.tex
\begin{theoremEnd}[malte,restate command=defsubvectorklmteststatistic]{definition}
    \label{def:subvector_klm_test_statistic}
    Let 
    $$
    \tilde S(\beta, \gamma) := \begin{pmatrix} X & W \end{pmatrix} - (y - X \beta - W \gamma) \frac{(y - X \beta - W \gamma)^T M_Z \begin{pmatrix} X & W \end{pmatrix}}{(y - X \beta - W \gamma)^T M_Z (y - X \beta - W \gamma)}.
    $$
    The subvector Lagrange multiplier test statistic is
    \begin{equation}
        \label{eq:subvector_klm_test_statistic}
        \LM(\beta) := (n - k) \min_{\gamma \in \BR^\mW} \frac{(y - X \beta - W \gamma)^T P_{P_Z \tilde S(\beta, \gamma)} (y - X \beta - W \gamma)}{(y - X \beta - W \gamma)^T M_Z (y - X \beta - W \gamma)}.
    \end{equation}
\end{theoremEnd}

%% file: theorems/prop_subvector_klm_test_statistic_chi_squared_statement.tex
\begin{theoremEnd}{technical_condition}
    \label{tc:subvector_klm}
    Assume there exists a $\gamma^\star \in \BR^\mW$ such that 
    $$
        \gamma^\star = \gamma_0 + (\Pi_W^T Z^T P_{P_Z \tilde S(\beta_0, \gamma^\star)} W)^{-1} \Pi_W^T Z^T P_{P_Z \tilde S(\beta_0, \gamma^\star)} \varepsilon,
    $$
    or, equivalently,
    $$
    \Pi_W^T Z^T P_{P_Z \tilde S(\beta_0, \gamma^\star)} (\varepsilon + W(\gamma_0 - \gamma^\star)) = 0.
    $$    
\end{theoremEnd}
\begin{theoremEnd}{theorem}%
    \label{prop:subvector_klm_test_statistic_chi_squared}
    Consider \Cref{model:1} and assume that \Cref{ass:1} and \Cref{tc:subvector_klm} hold.
    Under the null $\beta = \beta_0$, under both strong and weak instrument asymptotics, the subvector Lagrange multiplier test statistic is bounded from above by a random variable that is asymptotically $\chi^2(\mX)$ distributed.
\end{theoremEnd}

%% file: theorems/conj_subvector_klm_test_statistic.tex
\begin{theoremEnd}[malte,restate command=conjsubvectorklmteststatistic]{conjecture}
    \label{conj:subvector_klm_test_statistic_chi_squared}
    \Cref{tc:subvector_klm} holds with probability tending to 1 as $n \to \infty$.
\end{theoremEnd}

%% file: theorems/prop_subvector_klm_test_statistic_chi_squared_many_instruments.tex
\begin{theoremEnd}{theorem}%
    \label{prop:subvector_klm_test_statistic_chi_squared_many_instruments}
    Consider \Cref{model:1} and assume that \Cref{ass:1} holds as $k = k(n) \to \infty$ as $n\to\infty$.
    Here, one needs to reformulate \cref{ass:1} (b) with convergence of the difference between $\Psi$ and a Gaussian random variable to zero for all probabilities of Borel-measurable sets (see also \citeauthor{bentkus2003dependence}, \citeyear{bentkus2003dependence} and \citeauthor{chernozhukov2017central}, \citeyear{chernozhukov2017central}).
    Under the null $\beta = \beta_0$, under both strong and weak instrument asymptotics, the subvector Lagrange multiplier test statistic is bounded from above by a random variable that is asymptotically $\chi^2(\mX)$ distributed.
\end{theoremEnd}

%% file: theorems/cor_subvector_lagrange_multiplier_test.tex
\begin{theoremEnd}[category=index,malte_intro,restate command=corsubvectorlagrangemultipliertest]{corollary}%
    \label{cor:subvector_lagrange_multiplier_test}
    Consider \cref{model:1} and assume the conditions of \cref{prop:subvector_klm_test_statistic_chi_squared} or \cref{prop:subvector_klm_test_statistic_chi_squared_many_instruments}.
    Let $F_{\chi^2(\mX)}$ be the cumulative distribution function of a chi-squared random variable with $\mX$ degrees of freedom.
    Under both strong and weak instrument asymptotics, a test that rejects the null $H_0: \beta = \beta_0$ whenever $\LM(\beta) > F^{-1}_{\chi^2(\mX)}(1 - \alpha)$ has asymptotic size at most $\alpha$.
\end{theoremEnd}

%% file: tables/table_comparison_tests.tex
\begin{table}[htpb]
    \begin{adjustbox}{max width=\textwidth}
    \begin{tabular}{l|c|c|c|c|c|c}
        & degrees of & subvector degrees & weak instru- & closed-form & possibly & possibly \\
        & freedom & of freedom & ment robust? & confidence set? & empty? & unbounded? \\
        \hline
        Wald & $m$ & $\mX$ & No & Yes & No & No \\
        AR & $k$ & $k - \mW$ & Yes & Yes & Yes & Yes \\
        LM & $m$ & $\mX$ & Yes & No & No & Yes \\
        CLR & in-between & in-between & Yes & No & No & Yes \\
    \end{tabular}
\end{adjustbox}
\caption{
    \label{tab:comparison_tests}
    Properties of various tests for the causal parameter in instrumental variables regression.
    }
\end{table}

%% file: theorems/def_subvector_anderson_rubin_test_statistic.tex
\begin{theoremEnd}[malte, restate command=subvectorandersonrubinteststatistic]{definition}[\citeauthor{guggenberger2012asymptotic}, \citeyear{guggenberger2012asymptotic}]
    \label{def:subvector_anderson_rubin_test_statistic}
    The subvector Anderson-Rubin test statistic is
    \begin{align*}
        \AR(\beta) &:=  \min_{\gamma \in \BR^\mW} \frac{n - k}{k - \mW} \frac{(y - X \beta - W \gamma)^T P_Z (y - X \beta - W \gamma)}{(y - X \beta - W \gamma)^T M_Z (y - X \beta - W \gamma)} \\
        &= \frac{n - k}{k - \mW} \frac{(y - X \beta - W \hat\gamma_\liml)^T P_Z (y - X \beta - W \hat\gamma_\liml)}{(y - X \beta - W \hat\gamma_\liml)^T M_Z (y - X \beta - W \hat\gamma_\liml)}, %
    \end{align*}
    where $\hat\gamma_\liml$ is the LIML estimator using outcomes $y - X \beta$, covariates $W$, and instruments $Z$.
\end{theoremEnd}

%% file: theorems/prop_cis_closed_forms_intro.tex
\begin{theoremEnd}[normal]{proposition}%
    \label{prop:cis_closed_forms_intro}
    The inverse Anderson-Rubin test confidence sets describe quadrics in $\BR^{\mX}$.
    Let $S := \begin{pmatrix} X & W \end{pmatrix}$ and $B \in \BR^{\mX \times m}$ have ones on the diagonal and zeros elsewhere such that $BS = X$.
    Let $F_{\chi^2(k - \mW)}$ denote the cumulative distribution function of a $\chi^2(k - \mW)$ random variable and let
    $$
        \hat\beta_\kclass(\kappa) := \left( S^T (\Id_n - \kappa M_Z) S \right)^{\dagger} S^T (\Id_n - \kappa M_Z) y,
    $$
    where $\dagger$ denotes the Moore-Penrose pseudoinverse, be the k-class estimator for $(\beta_0^T, \gamma_0^T)^T$ using covariates $S$, instruments $Z$, and outcome $y$.
    Define
    \begin{align*}
    \kappa_{\AR}(\alpha) &= 1 + F^{-1}_{\chi^2(k - \mW)}(1 - \alpha) / (n - k) , \ \ \ \hat\sigma^2(\kappa) := \frac{1}{n-k} \| M_Z (y - S \hat\beta_\kclass(\kappa)) \|^2, \ \text{ and} \\
    A(\kappa) &:= \left( B \left(S^T (\Id_n - \kappa M_Z) S \right)^{-1} B^T \right)^{-1} =
    X^T (\Id_n - \kappa M_Z) X \\
    & \hspace{1.5cm} - X^T (\Id_n - \kappa M_Z) W \left(W^T (\Id_n - \kappa M_Z) W \right)^{-1} W^T (\Id_n - \kappa M_Z) X.
    \end{align*}
    If $\mW > 0$, let $\kappa_\mathrm{max} := \lambdamin{(W^T M_Z W)^{-1} W^T P_Z W} + 1$, where $\lambda_\mathrm{min}$ denotes the minimal eigenvalue of a matrix.
    Else $\kappa_{\max} := \infty$.
    \vspace*{0.2cm}
    
    \noindent If $\kappa_{\AR}(\alpha) \leq \kappa_\mathrm{max}$, then
    \begin{multline*}
        \CI_{\AR}(1 - \alpha) = \Big\{ \beta \in \BR^{\mX} \mid \left(\beta - B \hat\beta_\kclass(\kappa_{\AR}(\alpha))\right)^T A(\kappa_{\AR}(\alpha)) \left(\beta - B\hat\beta_\kclass(\kappa_{\AR}(\alpha)) \right) \\
        \leq \hat\sigma^2(\kappa_{\AR}(\alpha)) \cdot \big(F^{-1}_{\chi^2(k - \mW)}(1 - \alpha) - k \AR(\hat\beta_\kclass(\kappa_{\AR}(\alpha)))\big)\Big\}.
    \end{multline*}
    Else, $\CI_{\AR}(1 - \alpha) = \BR^{\mX}$.
\end{theoremEnd}

%% file: theorems/tc_2.tex
\begin{theoremEnd}[malte, restate command=tctwo]{technical_condition}
    \label{tc:2}
    Let $S = \begin{pmatrix} X & W \end{pmatrix}$.
    \begin{enumerate}[(a)]
        \item \label{tc:2:a} If $\mW > 0$, then
            $$
                \lambda := \lambdamin{ (S^T M_Z S)^{-1} S^T P_Z S} < \lambdamin{ (X^T M_Z X)^{-1} X^T P_Z X},
            $$
            or, equivalently,
            $$
                \lambdamin{ X^T (P_Z + (1 - \lambda) M_Z) X } > 0.
            $$
        \item Let $S^{(-i)}$ be $S$ with the $i$-th column removed. Then, for all $i = 1, \ldots, m$,
            $$
                \lambda < \lambdamin{ (S^{(-i)T} M_Z S^{(-i)})^{-1} S^{(-i)T} P_Z S^{(-i)} }
            $$
            or, equivalently,
            $$
                \lambdamin{ S^{(-i)T} (P_Z + (1 - \lambda) M_Z) S^{(-i)} } > 0.
            $$ \label{tc:2:b}
    \end{enumerate}
\end{theoremEnd}

%% file: theorems/prop_cis_bounded_intro.tex
\begin{theoremEnd}[malte_intro]{proposition}%
    \label{prop:cis_bounded_intro}
    Assume \Cref{tc:2} \ref{tc:2:a} \red{holds}. Let $\alpha > 0$.
    Let
    \begin{align*}
        J_\liml &:= (k - \mW) \cdot \min_b \AR(b) \\
        &= (n - k) \ \lambdamin{ \left(\begin{pmatrix} X & W & y \end{pmatrix}^T M_Z \begin{pmatrix} X & W & y \end{pmatrix} \right)^{-1}\begin{pmatrix} X & W & y \end{pmatrix}^T P_Z \begin{pmatrix} X & W & y \end{pmatrix}}
    \end{align*}
    be the LIML variant of the J-statistic \citep{guggenberger2012asymptotic} and let
    $$
        \lambda = (n - k) \ \lambdamin{ \left(\begin{pmatrix} X & W \end{pmatrix}^T M_Z \begin{pmatrix} X & W \end{pmatrix} \right)^{-1}\begin{pmatrix} X & W \end{pmatrix}^T P_Z \begin{pmatrix} X & W \end{pmatrix}}.
    $$
    be \citeauthor{anderson1951estimating}'s \citeyearpar{anderson1951estimating} likelihood-ratio test statistic of reduced rank.
    Then the (subvector) inverse Anderson-Rubin test confidence set $\CI_{\AR}(1 - \alpha)$ is bounded and nonempty if and only if
    \begin{align*}
        J_\liml \leq F^{-1}_{\chi^2(k - \mW)}(1 - \alpha) < \lambda \ \Leftrightarrow \
        F_{\chi^2(k - \mW)} \left(J_\liml \right) \leq 1 - \alpha < F_{\chi^2(k - \mW)} \left(  \lambda \right).
    \end{align*}
\end{theoremEnd}

%% file: theorems/cor_ar_bounded_iff_rank_test_rejects_intro.tex
\begin{theoremEnd}[malte_intro]{corollary}%
    \label{cor:ar_bounded_iff_rank_test_rejects_intro}
    Consider structural equations $y_i = S_i^T \beta_0 + \varepsilon_i$ and $S_i = Z_i^T \Pi + V_i$.
    We are interested in inference for each component of the causal parameter $\beta_0$.
    Thus, for each covariate $i = 1, \ldots, m$, we separate $X \leftarrow S^{(i)}$ and $W \leftarrow S^{(-i)}$ and construct confidence intervals $\CI_{\AR}^{(i)}(1 - \alpha)$ for the $i$-th component of $\beta_0$.
    Assume \Cref{tc:2} \ref{tc:2:b} \red{holds.}
    Then
    \begin{enumerate}[(a)]
        \item The subvector inverse Anderson-Rubin test confidence sets are jointly (un)bounded. That is, if $\CI_{\AR}^{(i)}(1 - \alpha)$ is (un)bounded for any $i$, then $\CI_{\AR}^{(i)}(1 - \alpha)$ is (un)bounded for all $i$.
        \item The subvector inverse Anderson-Rubin test confidence sets at level $\alpha$ are bounded (and thus confidence intervals) if and only if \citeauthor{anderson1951estimating}'s \citeyearpar{anderson1951estimating} likelihood-ratio test rejects the null hypothesis that $\Pi$ is of reduced rank at level $\alpha$. \label{cor:ar_bounded_iff_rank_test_rejects_intro:b}
        \item The subvector inverse Anderson-Rubin test confidence sets are jointly (non)empty. That is, if $\CI_{\AR}^{(i)}(1 - \alpha)$ is (non)empty for any $i$, then $\CI_{\AR}^{(i)}(1 - \alpha)$ is (non)empty for all $i$.
        \item The subvector inverse Anderson-Rubin test confidence sets at level $\alpha$ are jointly empty if and only if the LIML variant of the J-statistic \citep{guggenberger2012asymptotic} $J_\liml > F^{-1}_{\chi^2(k - m + 1)}(1 - \alpha)$. \label{cor:ar_bounded_iff_rank_test_rejects_intro:d}
    \end{enumerate}
\end{theoremEnd}

%% file: figures/figure_inverse_ar_different_alpha.tex
\begin{figure}[!ht]
    \centering
    \includegraphics[width=\textwidth]{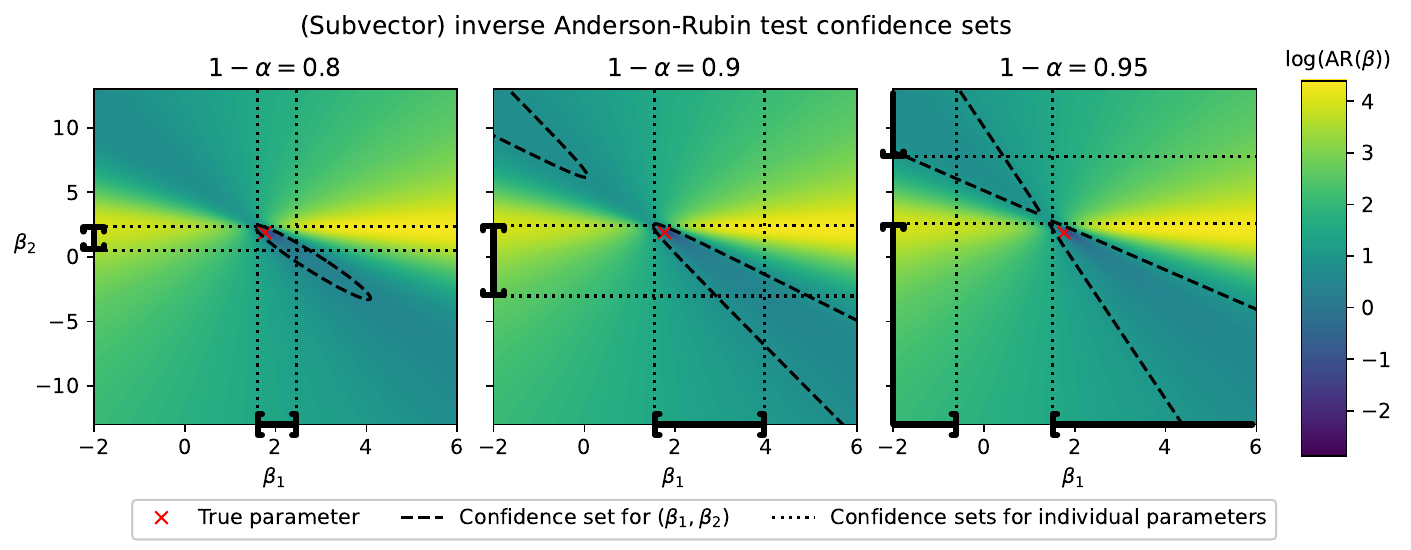}
    \caption{
        \label{fig:inverse_ar_different_alpha}
        The data was sampled from a Gaussian causal model with $n=100$, $k=3$, and $\mX = \mW = 1$.
        The $p$-value of \citeauthor{anderson1951estimating}'s \citeyearpar{anderson1951estimating} test of reduced rank is $0.07$.
        The multivariate confidence sets for $\beta = (\beta_1, \beta_2)$ are bounded for $1 - \alpha = 0.8$ and unbounded for $1 - \alpha = 0.9, 0.95$.
        The subvector confidence sets for the individual coefficients are jointly bounded for $1 - \alpha = 0.8, 0.9$ and unbounded for $1 - \alpha = 0.95$.
        Notably, for $1 - \alpha = 0.9$, the confidence set for $\beta = (\beta_1, \beta_2)$ is unbounded, while the subvector confidence sets for the individual coefficients are bounded.
        In all cases, the confidence sets for the individual coefficients are substantially smaller than projection-based confidence sets \citep{dufour2005projection}.
    }
\end{figure}

%% file: theorems/prop_inverse_ar_equal_to_wald.tex
\begin{theoremEnd}[malte_intro,restate,restate command=propinversearequaltowald,category=propinversearequaltowald]{proposition}
    \label{prop:inverse_ar_equal_to_wald}
    Let $\alpha > 0$.
    Assume that $J_\liml \leq F^{-1}_{\chi^2(k-\mW)}(1 - \alpha) < \lambda$ such that $\CI_{\AR}$ is bounded and nonempty by \Cref{prop:cis_bounded_intro}.
    Let
    \begin{multline*}        
    s(\kappa) := 
    y^T (\Id_n - \kappa M_Z ) y %
    - y^T (\Id_n - \kappa M_Z ) X (X^T (\Id_n - \kappa M_Z ) X)^{-1} X^T (\Id_n - \kappa M_Z ) y.
    \end{multline*}
    Recall that $\kappa_{\AR}(\alpha) = 1 + F^{-1}_{\chi^2(k - \mW)}(1 - \alpha) / (n-k)$ and define
    \begin{align*}
    \alpha_{\Wald \mid \AR}(\alpha) &:= 1 - F_{\chi^2(\mX)}(-s(\kappa_{\AR}(\alpha)) / \hat\sigma^2_{\Wald}(\kappa_{\AR}(\alpha))) \text{ and }
    \\\alpha_{\LR \mid \AR}(\alpha) &:= 1 - F^{-1}_{\chi^2(\mX)}( F^{-1}_{\chi^2(k - \mW)}(1 - \alpha) - J_\liml),
    \end{align*}
    where $\hat\sigma^2_{\Wald}(\kappa) = \frac{1}{n-m}\| y - X \hat\beta_\kclass(\kappa) \|^2$. Then
    $$
    \CI_{\AR}(1 - \alpha) = \CI_{\LR}(1 - \alpha_{\LR \mid \AR}(\alpha)) = \CI_{\Wald_{\hat\beta_\kclass(\kappa_{\AR}(\alpha))}}(1 - \alpha_{\Wald \mid \AR}(\alpha)).
    $$
        
\end{theoremEnd}
\begin{proofEnd}
    We apply \Cref{prop:cis_closed_forms}. Write $s := s(\kappa_{\AR}(\alpha))$.
    \paragraph*{Step 1: $J_\liml \leq F^{-1}_{\chi^2(k - \mW)}(1 - \alpha) < \lambda$ implies $s \leq 0$.}

    From 
    \begin{multline*}
    J_\liml \leq F^{-1}_{\chi^2(k - \mW)}(1 - \alpha) < \lambda \\
    \Leftrightarrow \lambdamin{ (S^T M_Z S)^{-1} S^T P_Z S } \leq \kappa_{\AR}(\alpha) - 1 < \lambdamin{ (\begin{pmatrix} S & y \end{pmatrix}^T M_Z \begin{pmatrix} S & y \end{pmatrix})^{-1} \begin{pmatrix} S & y \end{pmatrix}^T P_Z \begin{pmatrix} S & y \end{pmatrix} } \\
    \overset{\Cref{lem:kappa_pos_definite}}{\Leftrightarrow}
    \lambdamin{ S^T (P_Z - \kappa_{\AR}(\alpha) M_Z) S } \geq 0 > \lambdamin{ \begin{pmatrix} S & y \end{pmatrix}^T (P_Z - \kappa_{\AR}(\alpha) M_Z) \begin{pmatrix} S & y \end{pmatrix} }.
    \end{multline*}
    The eigenvalues of $\begin{pmatrix} S & y \end{pmatrix}^T (P_Z - \kappa_{\AR}(\alpha) M_Z) \begin{pmatrix} S & y \end{pmatrix}$ and $ S^T (P_Z - \kappa_{\AR}(\alpha) M_Z) S $ interleave and thus
    \begin{multline*}    
    \lambda_2 \left(\begin{pmatrix} S & y \end{pmatrix}^T (P_Z - \kappa_{\AR}(\alpha) M_Z) \begin{pmatrix} S & y \end{pmatrix} \right) \geq \lambdamin{ S^T (P_Z - \kappa_{\AR}(\alpha) M_Z) S } \geq 0 \\
    \Rightarrow \det( \begin{pmatrix} S & y \end{pmatrix}^T (P_Z - \kappa_{\AR}(\alpha) M_Z) \begin{pmatrix} S & y \end{pmatrix} ) \leq 0.
    \end{multline*}
    Applying the formula for the determinant of a block matrix, we get
    \begin{multline*}
        \det( S^T (P_Z - \kappa_{\AR}(\alpha) M_Z) S ) \cdot s = \det( \begin{pmatrix} S & y \end{pmatrix}^T (P_Z - \kappa_{\AR}(\alpha) M_Z) \begin{pmatrix} S & y \end{pmatrix} ) \Rightarrow s \leq 0.
    \end{multline*}

    \paragraph*{Step 2: $\CI_{\AR}(1 - \alpha) = \CI_{\Wald_{\hat\beta_\kclass(\kappa_{\AR}(\alpha))}}(1 - \alpha_{\Wald \mid \AR}(\alpha))$}

    We calculate
    \begin{multline*}
    \sigma_{\AR}^2(\kappa_{\AR}(\alpha)) \cdot ( F^{-1}_{\chi^2(k - \mW)}(1 - \alpha) - k \AR(\hat\beta_\kclass(\kappa_{\AR}(\alpha)))) \\
    = (\kappa_{\AR}(\alpha) - 1) \cdot \| M_Z (y - S \hat\beta_\kclass(\kappa_{\AR}(\alpha)) ) \|^2 - \| P_Z (y - S \hat\beta_\kclass(\kappa_{\AR}(\alpha)) ) \|^2 \\
    = \left( y^T - y^T \left(P_Z - \kappa_{\AR}(\alpha) M_Z \right) S \left( S^T \left( P_Z - \kappa_{\AR}(\alpha) M_Z \right) S \right)^{-1} S^T \right) \left(P_Z - \kappa_{\AR}(\alpha) M_Z \right) \\
    \left( y - S \left( S^T \left( P_Z - \kappa_{\AR}(\alpha) M_Z \right) S \right)^{-1} S^T \left( P_Z - \kappa_{\AR}(\alpha) M_Z \right) y \right)
    = - s.
    \end{multline*}
    
    Clearly $\hat\sigma^2_{\Wald}(\kappa_{\AR}(\alpha)) \geq 0$ and thus 
    \begin{multline*}        
    \hat\sigma^2_{\Wald}(\kappa_{\AR}(\alpha)) \cdot F^{-1}_{\chi^2(\mX)}(1 - \alpha_{\Wald \mid \AR}(\alpha)) = - s \\
    = \sigma_{\AR}^2(\kappa_{\AR}(\alpha)) \cdot ( F^{-1}_{\chi^2(k - \mW)}(1 - \alpha) - k \AR(\hat\beta_\kclass(\kappa_{\AR}(\alpha)))).
    \end{multline*}

    \paragraph*{Step 3: $\CI_{\AR}(1 - \alpha) = \CI_{\LR}(1 - \alpha_{\LR \mid \AR}(\alpha))$}

    From \Cref{prop:cis_closed_forms}, the two confidence sets are equal if
    $$
    \kappa_{\AR}(\alpha) = \kappa_{\LR}(\alpha_{\LR \mid \AR}(\alpha)) \Leftrightarrow F^{-1}_{\chi^2(k - \mW)}(1 - \alpha) = J_\liml + F^{-1}_{\chi^2(\mX)}(1 - \alpha_{\LR \mid \AR}(\alpha)).
    $$
    By assumption $J_\liml \leq F^{-1}_{\chi^2(k - \mW)}(1 - \alpha)$, so we can solve for $\alpha_{\LR \mid \AR}(\alpha)$.
\end{proofEnd} %

%% file: sec_numerical_analysis.tex
\section{Numerical analysis of the subvector Lagrange multiplier test}
\label{sec:numerical_analysis}
In this section, we compare the empirical size and power of the proposed subvector Lagrange multiplier test with other subvector tests.

We start by comparing the empirical sizes of the different tests for the data-generating process proposed by \citet{guggenberger2012asymptotic} in \Cref{sec:numerical_analysis:guggenberger12_size}.
As noted by \citet{guggenberger2012asymptotic}, the subvector variant of the Lagrange multiplier test obtained by plugging in the LIML estimator is size-inflated.
The Lagrange multiplier test proposed by us is size-correct, \red{meaning that it does not reject the null hypothesis more often than the nominal level $\alpha$.}

Next, we compare the power of the different tests using a variation of \citeauthor{guggenberger2012asymptotic}'s (\citeyear{guggenberger2012asymptotic}) data-generating process in \Cref{sec:numerical_analysis:guggenberger12_power}.
The subvector Lagrange multiplier test proposed by us has power comparable to the subvector conditional likelihood-ratio test and among the highest power out of the tests that are size-correct according to \Cref{sec:numerical_analysis:guggenberger12_size}.

Finally, we compare the empirical sizes of the different tests for the data-generating process proposed by \citet{kleibergen2021efficient} in \Cref{sec:numerical_analysis:kleibergen19}.
Again, the subvector Lagrange multiplier test we propose is size-correct.

\subsection{Empirical sizes of subvector tests for \citeauthor{guggenberger2012asymptotic}'s (\citeyear{guggenberger2012asymptotic}) data-generating process}
\label{sec:numerical_analysis:guggenberger12_size}
\citet{guggenberger2012asymptotic} suggest a data-generating process under which the subvector variant of the Lagrange multiplier test obtained by plugging in the LIML instead of direct minimization as in \Cref{def:subvector_klm_test_statistic} is size inflated.
In this data-generating process 
\begin{equation*}
\begin{split}
    &\mX = \mW = 1 \\
    &k \in \{5, 10, 15, 20, 30\} \\
    &Q = \Id_k
\end{split} \quad
\begin{split}
    &\sqrt{n} \| \Pi_X \| = 100 \\
    &\sqrt{n} \| \Pi_W \| = 1 \\
    & \frac{ | \langle \Pi_W,  \Pi_X \rangle | } {\| \Pi_W \| \| \Pi_X \|} = 0.95
\end{split} \quad
\begin{pmatrix} \varepsilon_i \\ V_{X, i} \\ V_{W, i} \end{pmatrix} \overset{\text{i.i.d.}}{\sim} \mathcal{N}\left(0, \begin{pmatrix} 1 & 0 & 0.95 \\ 0 & 1 & 0.3 \\ 0.95 & 0.3 & 1 \end{pmatrix} \right).
\end{equation*}
We simulate according to this model (see Appendix \ref{sec:simulation_details} for details) 10'000 times, each time drawing $n=1000$ observations, and report the empirical sizes of the different subvector tests in \Cref{tab:guggenberger12_empirical_sizes}.
\input{tables/table_guggenberger12_empirical_sizes.tex}

We replicate the size-distortion found by \citet{guggenberger2012asymptotic} for the subvector Lagrange multiplier test obtained by \red{plugging} in the LIML estimator.
However, our proposed subvector Lagrange multiplier test is size-correct.
\Cref{fig:guggenberger12_qqplots} in Appendix \ref{sec:additional_figures} shows quantile-quantile plots comparing the $p$-values of the different tests to the uniform distribution.
Empirical sizes of the different subvector tests for $n=50$ and $100$ are in Appendix \ref{sec:additional_tables}.
The subvector Anderson-Rubin, conditional likelihood-ratio, and Lagrange multiplier tests are size-correct as long as $n/k > 2$.

\subsection{Power of subvector tests for a variation of \citeauthor{guggenberger2012asymptotic}'s (\citeyear{guggenberger2012asymptotic}) data-generating process}
\label{sec:numerical_analysis:guggenberger12_power}
Next, we compare the power of different tests using a variation of \citeauthor{guggenberger2012asymptotic}'s (\citeyear{guggenberger2012asymptotic}) data-generating process.
In the original process, identification is so weak that all size-correct tests have no power above the significance level.
To address this, we increase identification by scaling $\Pi_W$ such that $\sqrt{n} \| \Pi_W \| = 10$.
We then plot the rejection frequencies of weak-instrument-robust tests in \Cref{fig:guggenberger12_power}.
\input{figures/figure_guggenberger12_power.tex}

The subvector Anderson-Rubin test with \citeauthor{guggenberger2019more}'s \citeyearpar{guggenberger2019more} critical values is indistinguishable from the test using $\chi^2(k-\mW)$ critical values.
Both are size-correct but have substantially less power than the (also size-correct) subvector conditional likelihood-ratio and subvector Lagrange multiplier tests.
The power curves of the subvector Lagrange multiplier and the subvector conditional likelihood-ratio tests overlap, except for a small region around $\beta = 0.92$.
Here, the subvector Lagrange multiplier incurs a loss in power due to the test being a \red{quadratic} function of the \red{Anderson-Rubin test statistic's derivative}.

\red{In \Cref{fig:guggenberger12_power_rho=0.8,fig:guggenberger12_power_rho=0.99} in Appendix, we show power curves for the same data-generating process but with correlations $\rho := \langle \Pi_W, \Pi_X \rangle / ( \| \Pi_W \| \| \Pi_X \| ) = 0.8$ and $0.99$.
As $\rho$ grows towards 1, identification becomes weaker and the power of all tests decreases.
Additionally, the subvector Lagrange multiplier test's power loss around $\beta = 0.92$ becomes more pronounced.
This power loss occurs because the derivative of the Anderson-Rubin test statistic is often close to zero around $\beta = 0.92$ when $\rho$ is close to 1, see also the discussion in \Cref{sec:main_results:subvector_lagrange_multiplier}.
}

In \Cref{fig:guggenberger12_d_power}, we compare the power of subvector tests for the causal effect of an included exogenous regressor.
We use the same data-generating process as in \Cref{fig:guggenberger12_power}, but test for the causal effect $\delta$ for the first components of the $k=10$ instruments $Z$.
That is, we set $\mX = 0, W \leftarrow \begin{pmatrix} X & D \end{pmatrix}$, $D = Z_1$, and $Z = \begin{pmatrix} Z_2 & Z_2 & \ldots & Z_{10} \end{pmatrix}$ (see also Appendix \ref{sec:exogenous_variables}).
In our variation of \citeauthor{guggenberger2012asymptotic}'s (\citeyear{guggenberger2012asymptotic}) data-generating process, the instruments are valid and thus the true causal effect $\delta_0 = 0$.
\red{We implement the subvector CLR as suggested by \citet{kleibergen2021efficient} by including $D$ as both an endogenous regressor and an instrument.}
See also the discussion at the end of \Cref{sec:main_results:subvector_lagrange_multiplier}.
\input{figures/figure_guggenberger12_d_power.tex}

The power of the subvector Anderson-Rubin test with \citeauthor{guggenberger2019more}'s \citeyearpar{guggenberger2019more} critical values is indistinguishable from the test using $\chi^2(8)$ critical values.
The power of the subvector conditional likelihood-ratio and subvector Lagrange multiplier tests are also very similar, with the subvector Lagrange multiplier test being minimally more powerful for $|\delta| \leq 0.15$.

\subsection{Empirical sizes of subvector tests for \citeauthor{kleibergen2021efficient}'s (\citeyear{kleibergen2021efficient}) data-generating process}
\label{sec:numerical_analysis:kleibergen19}
Let $\Omega_{VV \cdot \varepsilon} := \Cov( \begin{pmatrix} V_{X, i} & V_{W, i} \end{pmatrix} ) - \Cov( \begin{pmatrix} V_{X, i} & V_{W, i} \end{pmatrix}, \varepsilon_i ) \Var(\varepsilon_i)^{-1}  \Cov( \varepsilon_i, \begin{pmatrix} V_{X, i} & V_{W, i} \end{pmatrix} )$ be the covariance of $\begin{pmatrix} V_{X, i} & V_{W, i} \end{pmatrix}$ conditional on $\varepsilon_i$.
\citet{kleibergen2021efficient} prove that under \Cref{ass:1}, the asymptotic distribution of the subvector (conditional) likelihood-ratio test statistic is fully captured by
$$
R \Lambda^T \Lambda R^T := n \Omega_{VV \cdot \varepsilon}^{-1/2} \Pi^T Q \Pi \Omega_{VV \cdot \varepsilon}^{-1/2},
$$
where $R$ is orthogonal and $\Lambda$ is diagonal.
They parametrize
$$
\Lambda = \begin{pmatrix} \sqrt{\lambda_1} & 0 \\ 0 & \sqrt{\lambda_2} \end{pmatrix} \text{ for } 0 \leq \lambda_1, \lambda_2 \leq 100 \text{ and } R = \begin{pmatrix} \cos(\tau) & -\sin(\tau) \\ \sin(\tau) & \cos(\tau) \end{pmatrix} \text{ for } 0 \leq \tau < \pi.
$$
We use the same parametrization (see Appendix \ref{sec:simulation_details} for details) and present the empirical sizes of the different tests for $k=100$ and $\Omega = \Cov( \begin{pmatrix} \varepsilon & V_X & V_W \end{pmatrix}) = \Id_3$ in \Cref{fig:kleibergen19_identity_k100}.
Unlike the subvector Anderson-Rubin and (conditional) likelihood-ratio statistics, the distribution of the subvector Lagrange multiplier statistic is not fully captured by $R \Lambda^T \Lambda R^T$.
We include additional plots for $\Omega=\Id_3$ and $k=5, 20$ and $\Omega$ as in \citet{guggenberger2012asymptotic} and $k=5, 20, 100$ in \Cref{fig:kleibergen19_guggenberger12_k5_20,fig:kleibergen19_identity_k5_20,fig:kleibergen19_guggenberger12_k100} in Appendix \ref{sec:additional_figures}.
\citeauthor{guggenberger2012asymptotic}'s (\citeyear{guggenberger2012asymptotic}) data-generating process corresponds to $\lambda_1 \approx 1, \lambda_2 \approx 120'000$, and $\tau \approx 0$.
In all settings, the subvector Lagrange multiplier test we propose is size-correct.
\input{figures/figure_kleibergen19_identity_k100.tex}

%% file: tables/table_guggenberger12_empirical_sizes.tex
\begin{table}[h]
    \center
    \begin{tabular}{l | rrrrr | rrrrr}
         & \multicolumn{5}{r}{$\alpha = 0.01$} & \multicolumn{5}{|r}{$\alpha = 0.05$} \\
        \multicolumn{1}{ r |}{$k=$ } & 5 & 10 & 15 & 20 &  30 & 5 & 10 & 15 & 20 & 30 \\
        \hline
        AR & 0.4\% & 0.2\% & 0.1\% & 0.1\% & 0.1\% & 2.7\% & 2.0\% & 1.5\% & 1.1\% & 0.9\% \\
        AR (GKM) & 1.2\% & 0.9\% & 0.8\% & 0.8\% & 0.7\% & 4.9\% & 5.1\% & 4.9\% & 4.5\% & 4.0\% \\
        CLR & 0.3\% & 0.1\% & 0.2\% & 0.1\% & 0.1\% & 2.7\% & 1.9\% & 1.5\% & 1.2\% & 1.2\% \\
        LM (ours) & 0.3\% & 0.2\% & 0.2\% & 0.3\% & 0.2\% & 1.9\% & 1.7\% & 1.7\% & 2.1\% & 1.8\% \\
        LM (LIML) & 1.2\% & 2.6\% & 4.1\% & 5.3\% & 8.3\% & 6.2\% & 9.6\% & 13\% & 16\% & 19\% \\
    \end{tabular}
        
    \caption{
        \label{tab:guggenberger12_empirical_sizes}
        \tableguggenbergerempiricalsizes{$n=1000$}
    }
\end{table}

%% file: figures/figure_guggenberger12_power.tex
\begin{figure}[ht]
    \centering
    \includegraphics[width=0.8\textwidth]{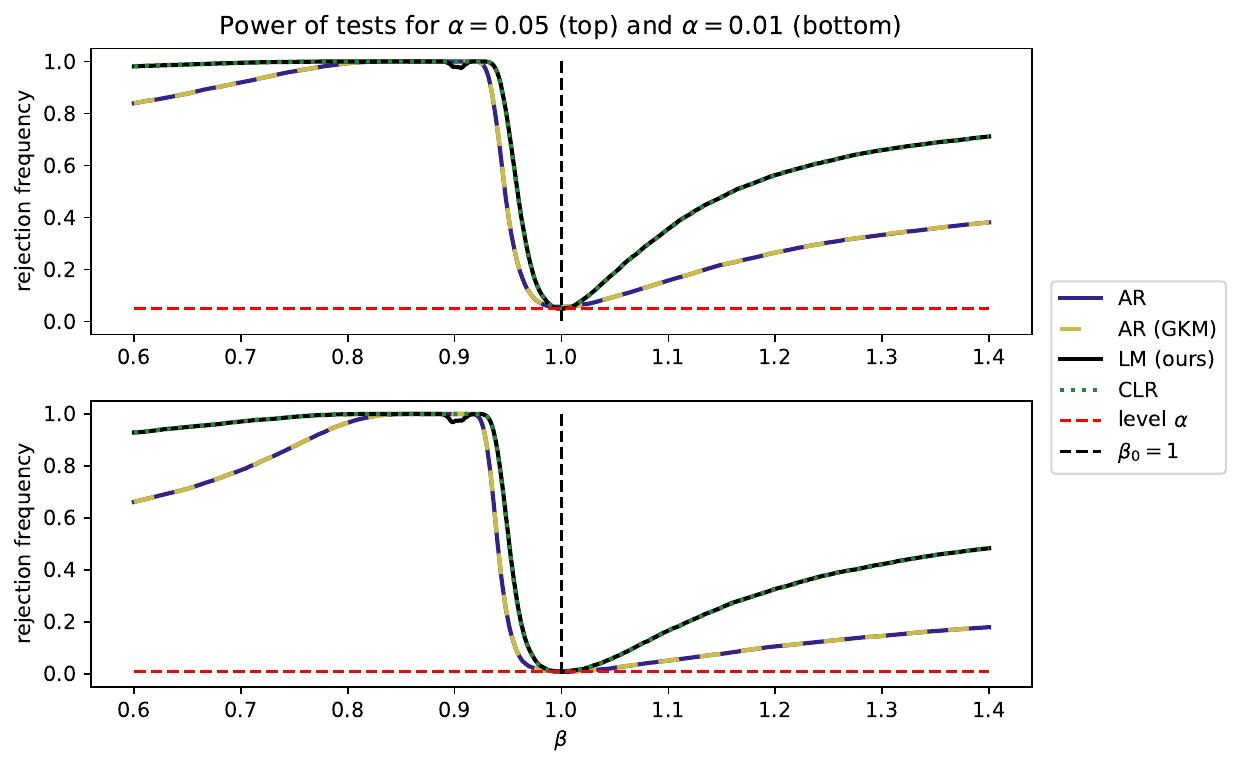}
    \caption{
        \label{fig:guggenberger12_power}
        Power curves of various weak-instrument-robust subvector tests, based on 10'000 simulations from the data-generating process proposed by \citet{guggenberger2012asymptotic} with $n=1000, k=10$, but with $\Pi_W$ scaled such that $\sqrt{n} \| \Pi_W \| = 10$.
        AR (GKM) uses the critical values of \citet{guggenberger2019more}.
    }
\end{figure}

%% file: figures/figure_guggenberger12_d_power.tex
\begin{figure}[ht]
    \centering
    \includegraphics[width=0.8\textwidth]{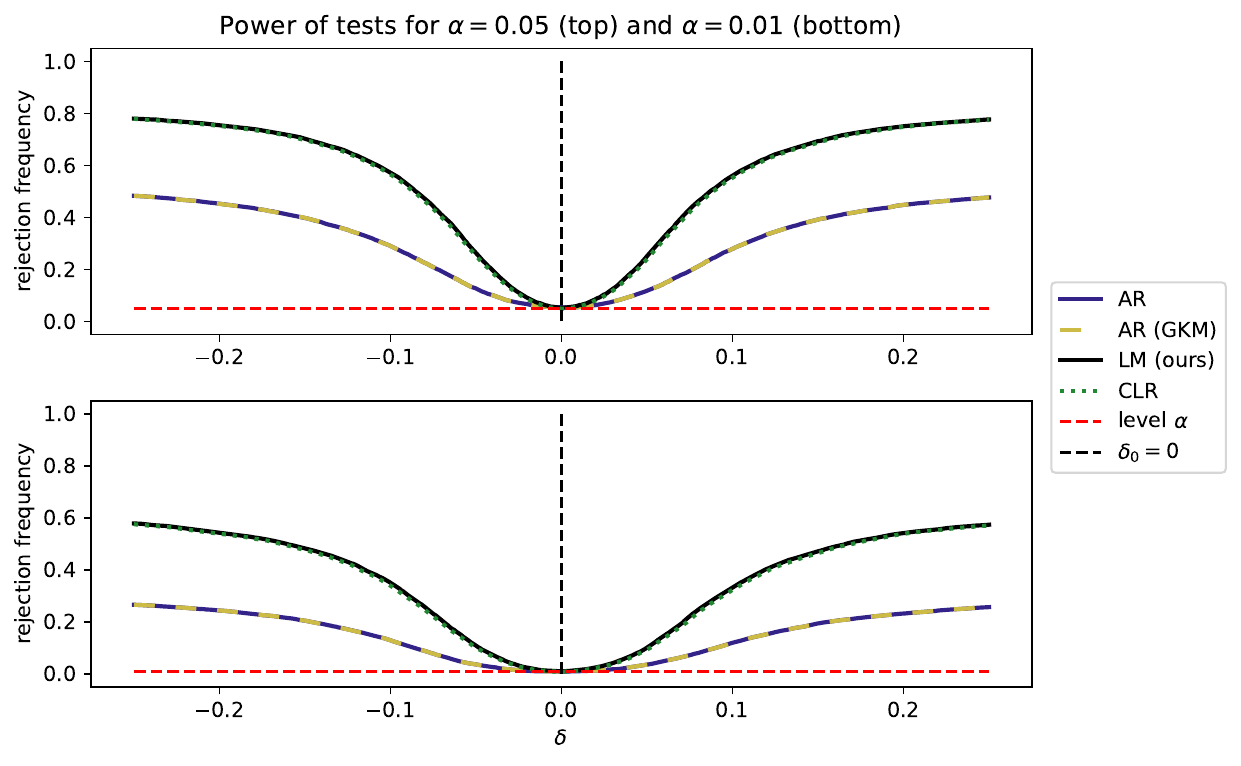}
    \caption{
        \label{fig:guggenberger12_d_power}
        Power curves for the causal effect of an included exogenous regressor based on 10'000 simulations from the data-generating process proposed by \citet{guggenberger2012asymptotic} with $n=1000, k=10$, but with $\Pi_W$ scaled such that $\sqrt{n} \| \Pi_W \| = 10$.
        AR (GKM) uses the critical values of \citet{guggenberger2019more}.
        We apply the CLR test by including the exogenous regressor $D$ into both the set of endogenous covariates of interest $X$ and instruments $Z$.
    }
\end{figure}

%% file: figures/figure_kleibergen19_identity_k100.tex
\begin{figure}[htpb]
    \centering
    \includegraphics[width=0.9\textwidth]{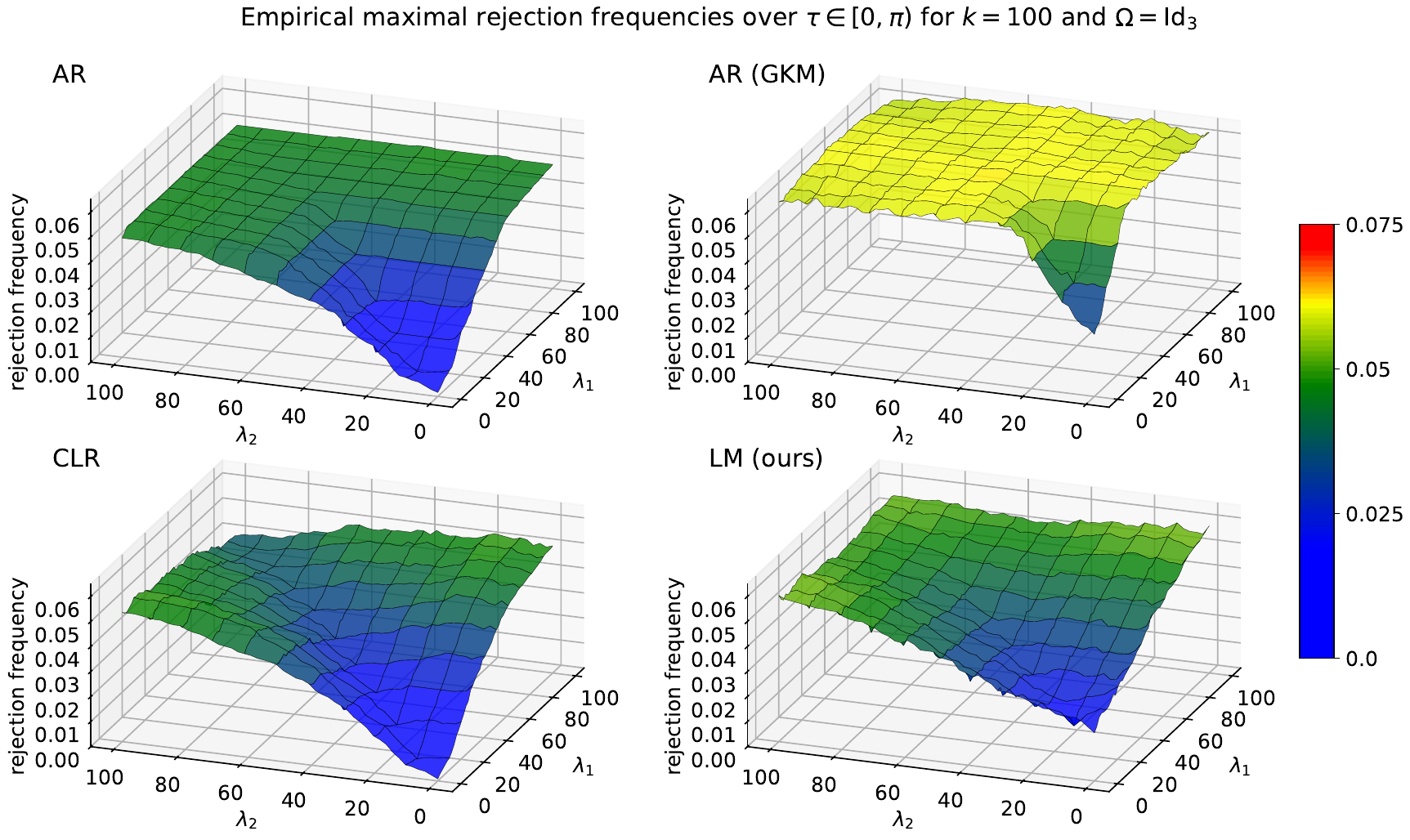}
    \caption{
        \label{fig:kleibergen19_identity_k100}
        \figurekleibergencaption{$k=100$}{$\Omega = \Id_3$}
    }
\end{figure}

%% file: sec_applications.tex
\section{Applications}
\label{sec:applications}
We apply the Wald, subvector Anderson-Rubin, subvector conditional likelihood-ratio, and our proposed subvector Lagrange multiplier test to two real-world datasets.
All estimates, test statistics, and $p$-values were generated using our \texttt{ivmodels} software package for Python.
We provide code and instructions to reproduce results at \href{https://github.com/mlondschien/ivmodels-simulations}{\texttt{github.com/mlondschien/ivmodels-simulations}}.

\subsection{\citeauthor{card1993using} \citeyearpar{card1993using}: Using geographic variation in college proximity to estimate the return to schooling}
\label{sec:applications:card}
\citeauthor{card1993using} \citeyearpar{card1993using} estimates the causal effect of length of education on hourly wages.
Their dataset is based on a 1976 subsample of the National Longitudinal Survey of Young Men and contains 3010 observations.
It includes years of education obtained by 1976, hourly wages in 1976, age in 1976, and indicators of whether the individual lived close to a public four-year college, a private four-year college, or a two-year college in 1966.
The dataset also contains variables or indicators on race, metropolitan area, region, and family background.
\citet{card1993using} defines (potential) experience as \texttt{age - education - 6}.

Like \citet{card1993using} and \citet{kleibergen2021efficient}, we set log-wages as the outcome, include education, experience, and experience squared as endogenous regressors, and race, metropolitan area, region, and family background indicators or variables as included exogenous regressors.

\subsubsection*{Inference for the causal effect of education of log-wages}
We compare three specifications, varying by the number of instruments: (i) Using proximity to a four-year college (public or private), age, and age squared as instruments ($k=3$).
This is similar to column 6 B of Table 3 of \citep{card1993using}.
(ii) Using the three college proximity indicators, age, and age squared as instruments ($k=5$).
This is equal to \citeauthor{kleibergen2021efficient}'s \citeyearpar{kleibergen2021efficient} specification.
(iii) Using age, age squared, and the interaction of the college proximity indicators with an indicator of whether both parents obtained less than 12 years of education as instruments ($k=8$).
\citet{card1993using} argue that this is a valid set of instruments, conditionally on family background variables, and report results using these instruments in column 3 of their Table 5.
We present estimates, test statistics, $p$-values, and confidence sets for the causal effect of education on log-wages in \Cref{tab:applications:card} (top).

\input{tables/table_card.tex}

For all specifications identification is weak, but not too weak to prohibit \red{informative} inference.
In particular, for all specifications and tests, the causal effect of education on log-wages is significant at level $\alpha = 0.01$.
Additional instruments improve identification, as measured by \citeauthor{anderson1951estimating}' \citeyearpar{anderson1951estimating} likelihood-ratio test statistic of reduced rank, which increases from specification (i) to (ii) to (iii).
However, as the critical value of the Anderson-Rubin and conditional likelihood-ratio tests increase with the number of instruments, the smallest $p$-values and confidence sets are achieved at specification (ii).
This is in contrast to the Lagrange multiplier test, whose $p$-value is smallest for specification (iii).
Still, the $p$-values of the conditional likelihood-ratio tests are smaller than those of the Lagrange multiplier test.
Also, the confidence sets obtained by inversion of the Lagrange multiplier test are the only ones that are not intervals and contain negative numbers.
This is due to the effect of the score being zero both at the minimum and maximum of the likelihood.
\red{The existence of two disjoint components in the confidence sets, given education as the single weakly identified endogenous variable, is consistent with the discussion in \Cref{sec:main_results:subvector_lagrange_multiplier}.}
See \Cref{fig:card_0} in Appendix \ref{sec:additional_figures} for the $p$-values of the tests for varying values of $\beta$.

\subsubsection*{Inference of the causal effect of education of log-wages for non-blacks}
A natural extension of the above specifications investigates how the causal effect of education on log-wages differs by race.
For this, we keep specifications (i), (ii), and (iii), but add the interaction between education and a race indicator to the set of endogenous variables and interact the instruments other than age and age squared with race.
We present estimates, test statistics, $p$-values, and confidence sets for the causal effect of education on log-wages for non-blacks in \Cref{tab:applications:card} (bottom).

Again, identification is weak, but not too weak to prohibit \red{informative} inference.
\citeauthor{anderson1951estimating}'s \citeyearpar{anderson1951estimating} likelihood-ratio test statistic of reduced rank increases from specification (i) to (ii) to (iii).
Again, the smallest $p$-values and confidence sets for the Anderson-Rubin and conditional likelihood-ratio test are achieved at specification (ii), whereas the Lagrange multiplier test achieves the smallest $p$-value at specification (iii) with $k=14$.
The $p$-value of the Lagrange multiplier test for specification (iii) is the smallest among all specifications and weak-instrument-robust tests.
Still, the confidence sets obtained by inversion of the Lagrange multiplier test are disconnected and include negative numbers.
If one discards the negative part of the confidence set, for example as such values are unlikely according to economic theory, the remaining confidence interval is smallest among those obtained by inverting weak-instrument-robust tests. 

The causal effect of education on wages for blacks was not significant at level $\alpha=0.2$ for any of the tests of interest.
We present estimates, test statistics, $p$-values, and confidence sets for the causal effect of education on log-wages for blacks in \Cref{tab:applications:card_3} in Appendix \ref{sec:additional_tables}.

\subsection{\citeauthor{tanaka2010risk} \citeyearpar{tanaka2010risk}: Risk and time preferences: Linking experimental and household survey data from vietnam}
\label{sec:applications:tanaka}
\citet{tanaka2010risk} study causes of risk preferences in Vietnam. Individuals from 25 households were interviewed for each of 289 villages.
\citet{tanaka2010risk} work with a subsample of in total 181 households from 9 villages.
From the interviews, they estimate measures of risk preferences, including the curvature of the utility function.
This is the dependent variable in the specifications of our interest.
Also measured is household income, which \citet{tanaka2010risk} split into the village means and the household's income relative to the village mean.
These regressors are assumed to be endogenous.
Included exogenous regressors are gender, age, education, distance to market, and whether the individual is Chinese or not.
Finally, \citet{tanaka2010risk} use rainfall in the village and an indicator of whether the head of the household can work as instruments for the household income variables.

We consider two specifications: (i) \citeauthor{tanaka2010risk}'s \citeyearpar{tanaka2010risk} original specification as described above (Table 4 column 2 and Table 5 column 2 bottom in their paper) and (ii) a specification where we include the interaction of the head of the household being able to work with rainfall as an additional instrument.
\citet{tanaka2010risk} discuss the estimated conditional causal effect for village mean income, income relative to the village mean, and gender.
We report estimates, confidence sets, and $p$-values for the causal effect of these variables in \Cref{tab:applications:tanaka}.
To make inference on gender, assumed to be exogenous, with the Anderson-Rubin and conditional likelihood-ratio test, we include gender as both an endogenous variable and an instrument.

\input{tables/table_tanaka.tex}
In specification (i), where the number of instruments is equal to the number of endogenous regressors, the LIML estimate is equal to the TSLS estimate and the Lagrange multiplier and conditional likelihood-ratio test are equal to the Anderson-Rubin test.
Adding the interaction of rainfall and the head of the household being able to work as an instrument does not substantially improve identification:
\citeauthor{anderson1951estimating}'s \citeyearpar{anderson1951estimating} test statistic of reduced rank increases from 6.051 to 6.071, but the $p$-value and the concentration parameter (which depend on the number of instruments) decrease.

\citet{tanaka2010risk} note that using the Wald test, the conditional causal effect of mean village income is significant at the 10\% level and that the conditional causal effects of income relative to the village mean and gender are not significant.
These findings still hold in specification (i) if using weak-instrument-robust tests.
After adding the interaction of rainfall and the head of the household being able to work as an additional instrument (specification ii), the $p$-value from the Lagrange multiplier test for the conditional causal effect of village mean income decreases from 0.061 to 0.060, whereas the $p$-values of the Anderson-Rubin and conditional likelihood-ratio tests increases from 0.061 to 0.150 and 0.105.
That is, the Lagrange multiplier test is the only weak-instrument-robust test that rejects the null of no conditional causal effect of village mean income on risk preferences at level 10\% for specification (ii).

We also report 90\% confidence sets.
The subvector inverse Lagrange multiplier test confidence set for the causal effect of the village mean income is the only weak-instrument-robust confidence set that does not include zero for specification (ii).
\red{It is composed of three disjoint components, consistent with two weakly identified endogenous variables and the discussion in \Cref{sec:main_results:subvector_lagrange_multiplier}.}
See \Cref{fig:tanaka} in Appendix \ref{sec:additional_figures} for the $p$-values of the tests for varying values of $\beta$.

%% file: tables/table_card.tex
\begin{table}[h]
    \begin{center}
      
    \begin{tabular}{lrrrrrr}
        \multicolumn{7}{l}{\emph{\citeauthor{card1993using} \citeyearpar{card1993using}: causal effect of education on log-wages}} \\
        \rule{0pt}{3ex} & \multicolumn{2}{c}{\ (i) $k=3$} & \multicolumn{2}{c}{\ \ (ii) $k=5$} & \multicolumn{2}{c}{\ \ (iii) $k=8$} \\
        \midrule
        estimate (TSLS) &  \multicolumn{2}{r}{0.132 (0.049)} & \multicolumn{2}{r}{0.145 (0.045)} & \multicolumn{2}{r}{0.130 (0.036)}  \\
        estimate (LIML) &  \multicolumn{2}{r}{0.132 (0.049)} & \multicolumn{2}{r}{0.172 (0.056)} & \multicolumn{2}{r}{0.148 (0.043)}  \\
        
        \rule{0pt}{4ex}test & \multicolumn{1}{l}{\ \ statistic} & \multicolumn{1}{l}{\ $p$-value} & \multicolumn{1}{l}{\ \ \ \ statistic} & \multicolumn{1}{l}{\ $p$-value} & \multicolumn{1}{l}{\ \ \ \ statistic} & \multicolumn{1}{l}{\ $p$-value}  \\
        \midrule
        Wald (TSLS) & 7.20 & 0.00728 & 10.53 & 0.00118 & 13.02 & 3.08e-4  \\
         \multicolumn{1}{r}{\ \ 95\% confidence set} & \multicolumn{2}{r}{[0.04, 0.23]} & \multicolumn{2}{r}{[0.06, 0.23]} & \multicolumn{2}{r}{[0.06, 0.20]} \\
        \rule{0pt}{3ex}Wald (LIML) & 7.20 & 0.00728 & 9.37 & 0.0022 & 11.92 & 5.55e-4  \\
          \multicolumn{1}{r}{\ \ 95\% confidence set} & \multicolumn{2}{r}{[0.04, 0.23]} & \multicolumn{2}{r}{[0.06, 0.28]} & \multicolumn{2}{r}{[0.06, 0.23]} \\
        \rule{0pt}{3ex}AR &6.83 & 0.00895 & 5.03 & 0.00175 & 2.88 & 0.00835 \\
          \multicolumn{1}{r}{\ \ 95\% confidence set} & \multicolumn{2}{r}{[0.04, 0.28]} & \multicolumn{2}{r}{[0.08, 0.36]} & \multicolumn{2}{r}{[0.04, 0.37]}  \\
        \rule{0pt}{3ex}CLR & 6.83 & 0.00895 & 10.84 & 0.00252 & 12.10 & 0.00304 \\
        \multicolumn{1}{r}{\ \ 95\% confidence set} & \multicolumn{2}{r}{[0.04, 0.28]} & \multicolumn{2}{r}{[0.07, 0.40]} & \multicolumn{2}{r}{[0.06, 0.29]} \\
        \rule{0pt}{3ex}LM (ours) & 6.83 & 0.00895 & 5.74 & 0.0166 & 7.63 & 0.00575 \\
        \multicolumn{1}{r}{\ \ 95\% confidence set} & \multicolumn{2}{r}{[0.04, 0.28]} & \multicolumn{2}{r}{$[-0.60, -0.06] \ \cup \ $} & \multicolumn{2}{r}{$[-0.75, -0.12] \ \cup \ $} \\
         & \multicolumn{2}{r}{\ } & \multicolumn{2}{r}{$[0.06, 0.47]$} & \multicolumn{2}{r}{$[0.06, 0.30]$}  \\
        \midrule
        rank  & 12.03 & 5.25e-4 & 15.47 & 0.00145 & 22.29 & 0.00107\\
        $J_\liml$ & & & 4.25 & 0.12 & 5.17 & 0.396  \\
      \end{tabular}
      \vspace{0.1cm}
      \vspace{0.3cm}
      \begin{tabular}{lrrrrrr}
        \multicolumn{7}{l}{\emph{\citeauthor{card1993using} \citeyearpar{card1993using}: causal effect of education on log-wages for non-blacks}} \\
        \rule{0pt}{3ex} & \multicolumn{2}{c}{\ (i) $k=4$} & \multicolumn{2}{c}{\ \ (ii) $k=8$} & \multicolumn{2}{c}{\ \ (iii) $k=14$} \\
        \midrule
        estimate (TSLS) &  \multicolumn{2}{r}{0.128 (0.050)} & \multicolumn{2}{r}{0.132 (0.032)} & \multicolumn{2}{r}{0.116 (0.026)}  \\
        estimate (LIML) &  \multicolumn{2}{r}{0.128 (0.050)} & \multicolumn{2}{r}{0.145 (0.036)} & \multicolumn{2}{r}{0.124 (0.028)}  \\
        
        \rule{0pt}{4ex}test & \multicolumn{1}{l}{\ \ statistic} & \multicolumn{1}{l}{\ $p$-value} & \multicolumn{1}{l}{\ \ \ \ statistic} & \multicolumn{1}{l}{\ $p$-value} & \multicolumn{1}{l}{\ \ \ \ statistic} & \multicolumn{1}{l}{\ $p$-value}  \\
        \midrule
        Wald (TSLS) & 6.51 & 0.0107 & 16.83 & 4.09e-05 & 20.21 & 6.95e-06 \\
         \multicolumn{1}{r}{\ \ 95\% confidence set} & \multicolumn{2}{r}{[0.03, 0.23]} & \multicolumn{2}{r}{[0.07, 0.19] } & \multicolumn{2}{r}{[0.07, 0.17]} \\
        \rule{0pt}{3ex}Wald (LIML)  & 6.51 & 0.0107 & 15.87 & 6.78e-05 & 19.13 & 1.22e-05 \\
          \multicolumn{1}{r}{\ \ 95\% confidence set} & \multicolumn{2}{r}{[0.03, 0.23]} & \multicolumn{2}{r}{[0.07, 0.22]} & \multicolumn{2}{r}{[0.07, 0.18]} \\
        \rule{0pt}{3ex}AR & 6.33 & 0.0119 & 4.35 & 5.88e-4 & 2.29 & 0.00866 \\
          \multicolumn{1}{r}{\ \ 95\% confidence set} & \multicolumn{2}{r}{[0.03, 0.28] } & \multicolumn{2}{r}{ [0.07, 0.27]} & \multicolumn{2}{r}{[0.03, 0.27] }  \\
        \rule{0pt}{3ex}CLR  & 6.33 & 0.0119 & 16.29 & 2.37e-4 & 17.95 & 2.57e-4 \\
        \multicolumn{1}{r}{\ \ 95\% confidence set} & \multicolumn{2}{r}{[0.03, 0.28] } & \multicolumn{2}{r}{[0.08, 0.25]} & \multicolumn{2}{r}{[0.06, 0.20]} \\
        \rule{0pt}{3ex}LM (ours)  & 6.33 & 0.0119 & 11.30 & 7.73e-4 & 13.93 & 1.89e-4 \\
        \multicolumn{1}{r}{\ \ 95\% confidence set} & \multicolumn{2}{r}{[0.03, 0.28] } & \multicolumn{2}{r}{$ [-0.54, -0.13] \ \cup \ $} & \multicolumn{2}{r}{$[-0.88, -0.28]\ \cup \ $} \\
         & \multicolumn{2}{r}{\ } & \multicolumn{2}{r}{$[0.07, 0.25]$} & \multicolumn{2}{r}{$[0.07, 0.20]$}  \\
        \midrule
        rank  & 12.14 & 4.93e-4 & 28.81 & 2.53e-05 & 45.91 & 3.35e-06 \\
        $J_\liml$ & & & 5.45 & 0.244 & 7.21 & 0.706 \\
      \end{tabular}
      \vspace{0.1cm}
      \caption{
        \label{tab:applications:card}
        Estimates, test statistics, $p$-values, and confidence sets for the causal effect of education on wages.
        See \Cref{sec:applications:card} for details.
      }
    \end{center}

\end{table}

%% file: tables/table_tanaka.tex
\begin{table}[htbp]
  \centering
  \begin{tabular}{lrrrrrr}
    \multicolumn{7}{l}{\emph{\citeauthor{tanaka2010risk} \citeyearpar{tanaka2010risk}: Causal effect of income and gender on risk preferences (i)}} \\
    \rule{0pt}{3ex}& \multicolumn{2}{c}{\ mean income} & \multicolumn{2}{c}{\ relative income} & \multicolumn{2}{c}{\ gender (male)} \\
    \midrule
    estimate (TSLS = LIML)  & 0.010 & (0.006) & 0.049 & (0.148) & -0.006 & (0.059)  \\
    
    \rule{0pt}{3ex}test & \multicolumn{1}{l}{\ statistic} & \multicolumn{1}{l}{\ $p$-value} & \multicolumn{1}{l}{\ \ statistic} & \multicolumn{1}{l}{\ $p$-value} & \multicolumn{1}{l}{\ \ statistic} & \multicolumn{1}{l}{\ $p$-value} \\
    \midrule
    \rule{0pt}{3ex}Wald (TSLS = LIML) & 3.411 & 0.065 & 0.111 & 0.739 & 0.011 & 0.917 \\
    \multicolumn{1}{r}{\ 90\% confidence set} & \multicolumn{2}{r}{[0.001, 0.019]} & \multicolumn{2}{r}{[-0.194, 0.292]} & \multicolumn{2}{r}{[-0.104, 0.091]} \\
    \rule{0pt}{3ex}AR = CLR = LM & 3.507 & 0.061 & 0.113 & 0.736 & 0.011 & 0.917  \\
    \multicolumn{1}{r}{\ 90\% confidence set} & \multicolumn{2}{r}{[0.001, 0.020]} & \multicolumn{2}{r}{[-0.234, 0.426]} & \multicolumn{2}{r}{[-0.103, 0.096]} \\
    \midrule
    &&&&&& \\
    \multicolumn{7}{l}{\emph{\citeauthor{tanaka2010risk} \citeyearpar{tanaka2010risk}: Causal effect of income and gender on risk preferences (ii)}}\\
    
    \rule{0pt}{3ex}& \multicolumn{2}{c}{\ mean income} & \multicolumn{2}{c}{\ relative income} & \multicolumn{2}{c}{\ gender (male)}  \\
    \midrule
    estimate (TSLS)  & 0.010 & (0.006) & 0.047 & (0.147) & -0.006 & (0.059)  \\
    estimate (LIML)  & 0.010 & (0.006) & 0.049 & (0.150) & -0.006 & (0.059)  \\
    
    \rule{0pt}{3ex}test & \multicolumn{1}{l}{\ statistic} & \multicolumn{1}{l}{\ $p$-value} & \multicolumn{1}{l}{\ \ statistic} & \multicolumn{1}{l}{\ $p$-value} & \multicolumn{1}{l}{\ \ \ \ statistic} & \multicolumn{1}{l}{\ $p$-value} \\
    \midrule
    \rule{0pt}{3ex}Wald (TSLS) & 3.493 & 0.062 & 0.101 & 0.751 & 0.010 & 0.921 \\
    \multicolumn{1}{r}{\ 90\% confidence set} & \multicolumn{2}{r}{[0.001, 0.019]} & \multicolumn{2}{r}{[-0.196, 0.289]} & \multicolumn{2}{r}{[-0.104, 0.092]}\\
    \rule{0pt}{3ex}Wald (LIML) & 3.501 & 0.061 & 0.106 & 0.745 & 0.010 & 0.921 \\
    \multicolumn{1}{r}{\ 90\% confidence set} & \multicolumn{2}{r}{[0.001, 0.020]} & \multicolumn{2}{r}{[-0.198, 0.296]} & \multicolumn{2}{r}{[-0.104, 0.092]}\\
    \rule{0pt}{3ex}AR & 1.894 & 0.150 & 0.164 & 0.849 & 0.114 & 0.892 \\
    \multicolumn{1}{r}{\ 90\% confidence set} & \multicolumn{2}{r}{[-0.001, 0.023]} & \multicolumn{2}{r}{[-0.435, 0.887]} & \multicolumn{2}{r}{[-0.129, 0.133]} \\
    \rule{0pt}{3ex}CLR & 3.570 & 0.105 & 0.108 & 0.767 & 0.010 & 0.929 \\
    \multicolumn{1}{r}{\ 90\% confidence set} & \multicolumn{2}{r}{[-0.000, 0.022]} & \multicolumn{2}{r}{[-0.341, 0.641]} & \multicolumn{2}{r}{[-0.129, 0.133]}\\
    \rule{0pt}{3ex}LM (ours) & 3.542 & 0.060 & 0.104 & 0.747 & 0.010 & 0.921 \\
    \multicolumn{1}{r}{\ 90\% confidence set} & \multicolumn{2}{r}{$[-\infty, -0.002] \,\cup$} & \multicolumn{2}{r}{$[-\infty, -1.016] \,\cup$} & \multicolumn{2}{r}{$[-\infty, -0.141] \,\cup$} \\
    & \multicolumn{2}{r}{$[0.001, 0.020] \,\cup$} & \multicolumn{2}{r}{$[-0.261, 0.473] \,\cup$} &  \multicolumn{2}{r}{$[-0.103, 0.097] \,\cup$}  \\
    & \multicolumn{2}{r}{$[0.024, \infty]$} & \multicolumn{2}{r}{$[14.017, \infty]$} & \multicolumn{2}{r}{$[0.139, \infty]$} \\
    \midrule

   \end{tabular}
   \caption{
      \label{tab:applications:tanaka}
      Estimated coefficients, test statistics, $p$-values and confidence sets for the causal effect of income and gender on risk preferences.
      See \Cref{sec:applications:tanaka} for details.
      For specification (i) \citeauthor{anderson1951estimating}'s \citeyearpar{anderson1951estimating} test statistic of reduced rank is 6.071 with a $p$-value of 0.014.
      The J-statistic is undefined.
      For specification (ii) \citeauthor{anderson1951estimating}'s \citeyearpar{anderson1951estimating} test statistic of reduced rank is 6.041 with a $p$-value of 0.049.
      The LIML variant of the J-statistic is 0.219 with a $p$-value of 0.64.
    }
\end{table}

%% file: sec_conclusion.tex
\section{Conclusion}
We introduced a novel weak-instrument-robust subvector Lagrange multiplier test for instrumental variables regression.
This test recovers the degrees of freedom of the standard Wald test, a property not achieved by previous weak-instrument-robust subvector tests.
We show that this test is asymptotically size-correct, either under a technical condition or as the number of instruments grows to infinity.
Numerical simulations confirm that the subvector Lagrange multiplier test maintains correct size and exhibits substantial power, comparable to the subvector conditional likelihood-ratio test.
Its ability to handle a large number of instruments without an increase in critical values makes it particularly suitable for modern applications where non-linear relationships are modeled through basis expansions.

In addition to the new test, we provide an analysis of confidence sets obtained by inverting the subvector Anderson-Rubin test.
We offer a closed-form solution for these confidence sets, revealing that they are centered around a k-class estimator.
A significant direct consequence is that for single coefficients, these subvector confidence sets are jointly bounded if and only if \citeauthor{anderson1951estimating}'s \citeyearpar{anderson1951estimating} likelihood-ratio test rejects the null hypothesis of underidentification.
We also establish a direct relationship between bounded inverse Anderson-Rubin test confidence sets and Wald-based confidence sets.

We apply various subvector tests to data from \citet{card1993using} and \citet{tanaka2010risk}.
There, the subvector Lagrange multiplier test can yield smaller weak-instrument-robust $p$-values when additional instruments through interactions are included.

Future work could focus on proving that the technical condition required for the general size-correctness of the subvector Lagrange multiplier test holds with high probability, as our empirical results suggest.
Overall, the subvector Lagrange multiplier test and the properties of subvector Anderson-Rubin confidence sets presented in this paper fill a gap in the literature and provide practitioners with additional powerful and reliable methods for subvector inference in the presence of weak instruments.

%% file: acknowledgements.tex
\section*{Acknowledgements}
Malte Londschien is supported by the ETH Foundations of Data Science.
We would like to thank Christoph Schultheiss, Cyrill Scheidegger, Fabio Sigrist, Felix Kuchelmeister, Frank Kleibergen, Gianna Wolfisberg, Jonas Peters, Juan Gamella, Leonard Henckel, Markus Ulmer, Maybritt Schillinger, Michael Law, Yuansi Chen, and Zijian Guo for helpful discussions and comments.
We also thank the co-editor, associate editor, and two anonymous reviewers for their constructive comments.

%% file: sec_exogenous_variables.tex
\section{An explicit treatment of included exogenous covariates}
\label{sec:exogenous_variables}
In practice, one often has additional exogenous variables (controls) $C$ entering the model.
For ease of exposition, many texts on instrumental variables regression reduces to the setting without such exogenous variables by considering the residuals $M_C Z$, $M_C X$, and $M_C y$ of $Z$, $X$, and $y$ after regressing out $C$.
However, one might also be interested in making inference on the parameters corresponding to exogenous variables.

We extend results in \Cref{sec:main_results:subvector_lagrange_multiplier} to the more general setting with included exogenous variables and prove that indeed one can reduce to the setting without exogenous variables by considering the residuals after regressing out the exogenous variables.
If one would like to make inference on the parameters corresponding to the exogenous variables, we prove that one can simply include them into both the instruments and endogenous variables.
For the full conditional likelihood-ratio test, this leads to a conservative result for the distribution of the test statistic \citep[see][Proposition 54]{londschien2025overview}.

\label{sec:exogenous_variables:assumptions}
\input{theorems/model_3.tex}
\input{figures/figure_iv_graph_exogenous.tex}
\input{theorems/assumption_2.tex}

Note that \Cref{ass:2} is equivalent to \Cref{ass:1} after including the exogenous included variables $C$, $D$ into the instruments $Z \leftarrow \begin{pmatrix} Z & D & C \end{pmatrix}$.
If we furthermore include them into the endogenous variables $X \leftarrow \begin{pmatrix} X & D \end{pmatrix}, W \leftarrow \begin{pmatrix} W & C \end{pmatrix}$, then \Cref{ass:1} no longer holds, as the last $\mD$ components of $V_X$ and the last $\mC$ components of $V_W$ are zero, and the corresponding columns of $\Psi$ are not Gaussian, but zero.

We reformulate and prove the statements in \Cref{sec:main_results:subvector_lagrange_multiplier} with included exogenous regressors.

\input{theorems/def_subvector_klm_test_statistic_exogenous.tex}
We first reduce to a model without the included exogenous regressors not of interest $C$.
\input{theorems/lem_lm_without_c.tex}
That is, up to the extra $\mC$ degrees of freedom in $n - k - \mC - \mD$, the Lagrange multiplier test statistic from \Cref{def:subvector_klm_test_statistic_exogenous} applied to $Z, X, W, C, D$, and $y$ is equal to the Lagrange multiplier test statistic applied to the residuals $M_C Z, M_C X, M_C W, M_C D$, and $M_C y$ after regressing out $C$.
To prove that it suffices to show that the Lagrange multiplier test statistic has the correct distribution for $\mC = 0$, we also need the following result.
\input{theorems/lem_ass_1_after_reduction.tex} 
\input{theorems/prop_subvector_klm_test_statistic_chi_squared_statement_exogenous.tex}
\input{theorems/prop_subvector_klm_test_statistic_chi_squared_many_instruments_exogenous.tex}

%% file: theorems/model_3.tex
\begin{model}
    \label{model:3}
    Let $y_i = X_i^T \beta_0 + W_i^T \gamma_0 + C_i^T \alpha_0 + D_i^T \delta + \varepsilon_i \in \BR$ with $X_i = Z_i^T \Pi_{ZX} + C_i^T \Pi_{CX} + D_i^T \Pi_{DX} + V_{X, i} \in \BR^\mX$ and $W_i = Z_i^T \Pi_{ZW} + C_i^T \Pi_{CW} + D_i^T \Pi_{DW} + V_{W, i} \in \BR^\mW$ for random vectors $Z_i \in \BR^k, C_i\in\BR^\mC, D_i \in \BR^\mD, V_{X, i}\in \BR^\mX,  V_{W, i}\in \BR^\mW$, and $\varepsilon_i \in \BR$ for $i=1\ldots, n$ and parameters $\Pi_{ZX} \in \BR^{k \times \mX}$, $\Pi_{CX} \in \BR^{\mC \times \mX}$, $\Pi_{DX} \in \BR^{\mD \times \mX}$, $\Pi_{ZW} \in \BR^{k \times \mW}$, $\Pi_{CW} \in \BR^{\mC \times \mW}$, $\Pi_{DW} \in \BR^{\mD \times \mW}$, $\beta_0 \in \BR^\mX$, $\gamma_0 \in \BR^\mW$, $\alpha_0 \in \BR^\mC$, and $\delta_0 \in \BR^\mD$.
    We call the $Z_i$ \emph{instruments}, the $X_i$ \emph{endogenous covariates of interest}, the $W_i$ \emph{endogenous covariates not of interest}, the $C_i$ \emph{exogenous (included) covariates not of interest}, the $C_i$ \emph{exogenous (included) covariates of interest}, and the $y_i$ \emph{outcomes}.
    The $V_{X, i}$, $V_{W, i}$ and $\varepsilon_i$ are \emph{errors}.
    These need not be independent across observations.
    Let $Z, X, W, C, D$, and $y$ be the matrices of stacked observations

    In \emph{strong instrument asymptotics}, we assume that $\begin{pmatrix} \Pi_{ZX} & \Pi_{ZW} \end{pmatrix}$ is fixed and of full rank column rank $\mX + \mW$.
    In \emph{weak instrument asymptotics} \citep{staiger1997instrumental}, we assume $\sqrt{n} \, \begin{pmatrix} \Pi_{ZX} & \Pi_{ZW} \end{pmatrix}$ is fixed and of full column rank $\mX + \mW$.
    Thus, $\begin{pmatrix} \Pi_{ZX} & \Pi_{ZW} \end{pmatrix} = \CO(\frac{1}{\sqrt{n}})$.
    Both asymptotics imply that $k \geq \mX + \mW$.
\end{model}

%% file: figures/figure_iv_graph_exogenous.tex
\begin{figure}[!h]
    \centering
    \begin{subfigure}[b]{0.45\textwidth}
        \begin{tikzpicture}[
            node distance=2cm and 2cm,
            >=Stealth,
            every node/.style={draw, circle, minimum size=1cm, inner sep=0pt},
            dashednode/.style={draw, circle, minimum size=1cm, inner sep=0pt, dashed},
            dashedarrow/.style={->, dashed},
            dashedtwinarrow/.style={<->, dashed}
        ]

        \node (Z) at (0,2) {Z};
        \node (X) at (2,0) {X};
        \node (Y) at (5.5,0) {Y};
        \node (C) at (0, -0.5) {C};
        \node[dashednode] (U) at (3.5,2) {U};

        \draw[->] (Z) -- (X) node[pos=0.3, below, draw=none] {$\Pi_{ZX}$};
        \draw[->] (X) -- (Y) node[pos=0.4, below, draw=none, yshift=5] {$\beta_0$};
        \draw[dashedarrow] (U) -- (X);
        \draw[dashedarrow] (U) -- (Y);
        \draw[dashedtwinarrow] (Z) -- (C);
        \draw[->] (C) -- (X) node[pos=0.7, yshift=3, below, draw=none] {$\Pi_{CX}$};
        \draw[->, bend right=25] (C) edge node[below, pos=0.4, draw=none, yshift=5] {$\alpha_0$} (Y);
        \end{tikzpicture}
    \end{subfigure}
    \begin{subfigure}[b]{0.45\textwidth}
        \begin{tikzpicture}[
            node distance=2cm and 2cm,
            >=Stealth,
            every node/.style={draw, circle, minimum size=1cm, inner sep=0pt},
            dashednode/.style={draw, circle, minimum size=1cm, inner sep=0pt, dashed},
            dashedarrow/.style={->, dashed},
            dashedtwinarrow/.style={<->, dashed}
        ]

        \node (Z) at (0,2) {Z};
        \node (W) at (2,1.5) {W};
        \node (X) at (2,0) {X};
        \node (Y) at (5.5,0) {Y};
        \node[dashednode] (U) at (3.5,2) {U};
        \node (C) at (0, -0.5) {C};
        \node (D) at (1.25, -1.75) {D};

        \draw[->] (Z) -- (X); %
        \draw[->] (Z) -- (W); %
        \draw[->] (X) -- (Y) node[pos=0.5, below, draw=none, yshift=5] {$\beta_0$};
        \draw[->] (W) -- (Y) node[pos=0.4, below, draw=none, yshift=5] {$\gamma_0$};
        \draw[dashedtwinarrow] (W) -- (X);
        \draw[dashedarrow] (U) -- (W);
        \draw[dashedarrow] (U) -- (X);
        \draw[dashedarrow] (U) -- (Y);

        \draw[dashedtwinarrow] (Z) -- (C);
        \draw[dashedtwinarrow] (D) -- (C);
        \draw[dashedtwinarrow] (Z) -- (D);
    
        \draw[->] (C) -- (X);
        \draw[->] (C) -- (W);
        \draw[->] (D) -- (X);
        \draw[->, bend left=20] (D) to (W);

        \draw[->, bend right=25] (C) edge node[below, pos=0.4, draw=none, yshift=5] {$\alpha_0$} (Y);
        \draw[->, bend right=20] (D) edge node[below, midway, draw=none, yshift=5] {$\delta_0$} (Y);

        \end{tikzpicture}
    \end{subfigure}
    \caption{
        \label{fig:iv_graph_exogenous}
        Causal graphs visualizing \cref{model:1}.
        On the left, $\mW=\mD=0$.
        On the right, we split the endogenous variables into endogenous variables of interest $X$ and endogenous variables not of interest $W$, and split the exogenous variables into exogenous variables of interest $D$ and exogenous variables not of interest $C$.
        First-stage parameters are not labeled on the right to avoid clutter.
    }
\end{figure}

%% file: theorems/assumption_2.tex
\begin{theoremEnd}[malte,restate command=assumptiontwo]{assumption}
    \label{ass:2}
    Let
    \begin{align*}
    \Psi &:= \begin{pmatrix} \Psi_{\varepsilon} & \Psi_{V_X} & \Psi_{V_W} \end{pmatrix}\\
    &:= \left(\begin{pmatrix} Z & D & C \end{pmatrix}^T \begin{pmatrix} Z & D & C \end{pmatrix} \right)^{-1/2} \begin{pmatrix} Z & D & C \end{pmatrix}^T \begin{pmatrix} \varepsilon & V_X & V_W \end{pmatrix} \in \BR^{(k + \mC + \mD) \times (1 + \mX + \mW)}.
    \end{align*}
    Assume there exist $\Omega \in \BR^{(1 + \mX + \mW) \times (1 + \mX + \mW)}$ and $Q \in \BR^{(k + \mC + \mD) \times (k + \mC + \mD)}$ positive definite such that, as $n \to \infty$,
    \begin{align*}
        &\mathrm{(a)} \ \ \frac{1}{n} \begin{pmatrix}\varepsilon & V_X & V_W\end{pmatrix}^T \begin{pmatrix}\varepsilon & V_X & V_W \end{pmatrix} \toP \Omega, \\
        &\mathrm{(b)} \ \ \vecop(\Psi) \tod \CN(0, \Omega \otimes \Id_{k + \mC + \mD}), \text{ and }\\
        &\mathrm{(c)} \ \ \frac{1}{n} \begin{pmatrix} Z & D & C \end{pmatrix}^T  \begin{pmatrix} Z & D & C \end{pmatrix} \toP Q.
    \end{align*}
\end{theoremEnd}

%% file: theorems/def_subvector_klm_test_statistic_exogenous.tex
\begin{theoremEnd}[malte,restate command=defsubvectorklmteststatisticexogenous]{definition}
    \label{def:subvector_klm_test_statistic_exogenous}
    Let 
    \begin{multline*}
        \tilde S(\beta, \gamma, \alpha, \delta) := \begin{pmatrix} X & W & D & C\end{pmatrix} - (y - X \beta - W \gamma - C \alpha - D \delta) \\
        \frac{(y - X \beta - W \gamma - C \alpha - D \delta)^T M_{[Z, C, D]} \begin{pmatrix} X & W & D & C \end{pmatrix}}{(y - X \beta - W \gamma - C \alpha - D \delta )^T M_{[Z, C, D]} (y - X \beta - W \gamma - C \alpha - D \delta )}.
    \end{multline*}
    The subvector Lagrange multiplier test statistic is
    \begin{multline}
        \label{eq:subvector_klm_test_statistic_exogenous}
        \LM(\beta, \delta) := (n - k - \mC - \mD) \\
        \min_{\gamma \in \BR^\mW, \ \alpha \in \BR^\mC} \frac{(y - X \beta - W \gamma - C \alpha - D \delta)^T P_{P_{[Z, C, D]} \tilde S(\beta, \gamma, \alpha, \delta)} (y - X \beta - W \gamma - C \alpha - D \delta)}{(y - X \beta - W \gamma - C \alpha - D \delta)^T M_{[Z, C, D]} (y - X \beta - W \gamma - C \alpha - D \delta)}.
    \end{multline}
\end{theoremEnd}

%% file: theorems/lem_lm_without_c.tex
\begin{theoremEnd}[malte_intro,category=lm_exogenous]{lemma}
    \label{lem:lm_without_c}
    The subvector Lagrange multiplier test statistic \eqref{eq:subvector_klm_test_statistic_exogenous} is equal to
    \begin{align*}
        &\LM(\beta, \delta)
        &= (n - k - \mC - \mD)
        \min_{\gamma \in \BR^\mW} \frac{ \| P_{P_{[M_C Z, M_C D]} M_C \tilde S(\beta, \gamma, \delta)} M_C (y - X \beta - W \gamma - D \delta) \|^2 }{ \| M_{[M_C Z, M_C D]} M_C (y - X \beta - W \gamma) \|^2},
    \end{align*}
    where
    \begin{multline*}
        M_C \tilde S(\beta, \gamma, \delta) = M_C \begin{pmatrix} X & W & D \end{pmatrix} - M_C (y - X \beta - W \gamma - D \delta) \\
        \frac{(y - X \beta - W \gamma)^T M_C M_{[M_C Z, M_C D]} M_C \begin{pmatrix} X & W & 0 \end{pmatrix}}{(y -  X \beta - W \gamma )^T M_C M_{[M_C Z, M_C D]} M_C (y - X \beta - W \gamma )}.
    \end{multline*}
\end{theoremEnd}
\begin{proofEnd}
We apply \Cref{prop:lm_derivative} with $Z \leftarrow \begin{pmatrix} Z & C & D \end{pmatrix}$, $X \leftarrow \begin{pmatrix} X & W & C & D \end{pmatrix}$, and $\beta = (\beta^T, \gamma^T, \alpha^T, \delta^T)$.
Together with $\tilde C(\beta, \gamma, \alpha, \delta) = C$ and $M_{[Z, C, D]} C = 0$, this implies that the derivative of the expression in \eqref{eq:subvector_klm_test_statistic_exogenous} with respect to $\alpha$ is
$$
-2 (n - k - \mC - \mD) \frac{ C^T (y - X \beta - W \gamma - C \alpha - D \delta)}{(y - X \beta - W \gamma)^T M_{[Z, C, D]} (y - X \beta - W \gamma)}.
$$
This is minimized at $\hat\alpha := (C^T C)^{-1} C^T (y - X \beta - W \gamma - D \delta)$ and
\begin{align*}
        &(n - k - \mC - \mD) \LM(\beta, \delta)\\
        &= \min_{\gamma \in \BR^\mW} \frac{\| P_{P_{[Z, C, D]} \tilde S(\beta, \gamma, \delta, \hat\alpha)} (y - X \beta - W \gamma - C \hat\alpha - D \delta) \|^2 }{\| M_{[Z, C, D]} (y - X \beta - W \gamma) \|^2} \\
        &= \min_{\gamma \in \BR^\mW} \frac{\| P_{P_{[Z, C, D]} \tilde S(\beta, \gamma, \delta, \alpha)} M_C (y - X \beta - W \gamma - D \delta) \|^2 }{\| M_{[Z, C, D]} (y - X \beta - W \gamma) \|^2 } \\
        &\overset{\Cref{lem:chain_projections}}{=} 
        \min_{\gamma \in \BR^\mW} \frac{ \| P_{P_{[M_C Z, M_C D]} M_C \tilde S(\beta, \gamma, \delta)} (M_C y - M_C X \beta - M_C W \gamma - M_C D \delta) \|^2 }{ \| M_{[M_C Z, M_C D]} (M_C y - M_C X \beta - M_C W \gamma) \|^2},
\end{align*}
where
\begin{multline*}
    M_C \tilde S(\beta, \gamma, \delta) = \begin{pmatrix} M_C X & M_C W & M_C D \end{pmatrix} - (M_C y - M_C X \beta - M_C W \gamma - M_C D \delta) \\
    \frac{(M_C y - M_C X \beta - M_C W \gamma)^T M_{[M_C Z, M_C D]} \begin{pmatrix} M_C X & M_C W & 0 \end{pmatrix}}{(M_C y - M_C  X \beta - M_C W \gamma )^T M_{[M_C Z, M_C D]} (M_C y - M_C X \beta - M_C W \gamma )}.
\end{multline*}
Here we applied \Cref{lem:chain_projections} to obtain $M_{[Z, C, D]} = M_{[M_C Z, M_C D]} M_C$ and $P_{P_{[Z, C, D]} \tilde S(\beta, \gamma, \delta, \hat\alpha)} M_C = P_{P_{[M_C Z, M_C D]} M_C \tilde S(\beta, \gamma, \delta)}$.
\end{proofEnd}

%% file: theorems/lem_ass_1_after_reduction.tex
\begin{theoremEnd}[malte_intro,category=lm_exogenous]{lemma}%
    \label{lem:ass_2_reduction}
    If \Cref{ass:2} applies to $Z, X, W, C, D$, and $y$, then it also applies to $M_C Z, M_C X, M_C D,$ and $M_C y$, the residuals after regressing out $C$.
\end{theoremEnd}
\begin{proofEnd}
Let
$$
\tilde\Psi_{[Z, D]} := \begin{pmatrix} Z & D \end{pmatrix}^T \begin{pmatrix} \varepsilon & V_X & V_W \end{pmatrix},\ \tilde\Psi_C := C^T \begin{pmatrix} \varepsilon & V_X & V_W \end{pmatrix} \text{, and } Q = \begin{pmatrix} Q_{[Z, D]} & Q_{[Z, D], C} \\ Q_{C, [Z, D]} & Q_C \end{pmatrix}.
$$
Then,
$$
\frac{1}{\sqrt{n}} \tilde \Psi := \frac{1}{\sqrt{n}} \begin{pmatrix} \tilde\Psi_{[Z, D]} \\ \tilde \Psi_C \end{pmatrix} = \frac{1}{\sqrt{n}}\left(\begin{pmatrix} Z & D & C \end{pmatrix}^T \begin{pmatrix} Z & D & C \end{pmatrix} \right)^{1/2} \Psi  \tod \CN(0, \Omega \otimes Q)
$$
by \Cref{ass:2} (b) and (c).

\paragraph{(c)} Calculate
\begin{align*}
    \begin{pmatrix} M_C Z & M_C D \end{pmatrix}^T  \begin{pmatrix} M_C Z & M_C D \end{pmatrix} &= \begin{pmatrix} Z & D \end{pmatrix}^T  \begin{pmatrix} Z & D \end{pmatrix} - \begin{pmatrix} Z & D \end{pmatrix}^T C (C^T C)^{-1} C^T \begin{pmatrix} Z & D \end{pmatrix} \\
    &\overset{\Cref{ass:2} (c)} \toP Q_{[Z, D]} - Q_{[Z, D], C} Q_{C}^{-1} Q_{C, [Z, D]} =: \tilde Q.
\end{align*}
$\tilde Q$ is the Schur complement of $Q_{C}$ in $Q$ and thus positive definite.

\paragraph{(a)} Calculate
\begin{multline*}
    \frac{1}{n} \begin{pmatrix} M_C \varepsilon & M_C V_X & M_C V_W \end{pmatrix}^T \begin{pmatrix} M_C \varepsilon & M_C V_X & M_C V_W \end{pmatrix} \\
    = \frac{1}{n} \begin{pmatrix} \varepsilon & V_X & V_W \end{pmatrix}^T \begin{pmatrix} \varepsilon & V_X & V_W \end{pmatrix} -  \frac{1}{n} \begin{pmatrix} \varepsilon & V_X & V_W \end{pmatrix}^T P_C  \begin{pmatrix} \varepsilon & V_X & V_W \end{pmatrix},
\end{multline*}
where $ \begin{pmatrix} \varepsilon & V_X & V_W \end{pmatrix}^T P_C \begin{pmatrix} \varepsilon & V_X & V_W \end{pmatrix} = (\tilde \Psi_C)^T (Z^T Z)^{-1} (\tilde \Psi_C)\toP \frac{1}{n}(\tilde \Psi_C)^T Q_C^{-1} (\tilde \Psi_C) = O(1)$ by \Cref{ass:2} (b).

\paragraph*{(b)}
We have
\begin{align*}
\begin{pmatrix} M_C Z & M_C D \end{pmatrix}^T & \begin{pmatrix} M_C \varepsilon & M_C V_X & M_C V_W \end{pmatrix} \\
&= \begin{pmatrix} Z & D \end{pmatrix}^T \begin{pmatrix} \varepsilon & V_X & V_W \end{pmatrix}
 - \begin{pmatrix} Z & D \end{pmatrix}^T C (C^T C)^{-1} C^T \begin{pmatrix} \varepsilon & V_X & V_W \end{pmatrix} \\
 &= \tilde \Psi_{[Z, D]}
 - \begin{pmatrix} Z & D \end{pmatrix}^T C (C^T C)^{-1} \tilde \Psi_{C} \\
 &\toP  \tilde \Psi_{[Z, D]} - Q_{[Z, D], C} Q_{C}^{-1} {}_{C} \tilde \Psi_C
 = \begin{pmatrix} \Id & - Q_{[Z, D]} Q_C^{-1} \end{pmatrix} \tilde \Psi
\end{align*}
As $\frac{1}{\sqrt{n}} \Cov(\tilde \Psi) = \Omega \otimes Q$, we have 
$$
\frac{1}{\sqrt{n}} \Cov( \begin{pmatrix} \Id & - Q_{[Z, D]} Q_C^{-1} \end{pmatrix}  \tilde \Psi) = \Omega \otimes \begin{pmatrix} \Id & - Q_{[Z, D]} Q_C^{-1} \end{pmatrix} Q \begin{pmatrix} \Id & - Q_{[Z, D]} Q_C^{-1} \end{pmatrix} = \Omega \otimes \tilde Q.
$$
Together with (c) from above, this yields
$$
\left( \begin{pmatrix} M_C Z & M_C D \end{pmatrix}^T \begin{pmatrix} M_C Z & M_C D \end{pmatrix} \right)^{1/2} \begin{pmatrix} M_C Z & M_C D \end{pmatrix}^T \begin{pmatrix} M_C \varepsilon & M_C V_X & M_C V_W \end{pmatrix} \tod \CN(0, \Omega \otimes \Id).
$$
\end{proofEnd}

%% file: theorems/prop_subvector_klm_test_statistic_chi_squared_statement_exogenous.tex
\begin{theoremEnd}[malte_intro,category=lm_exogenous]{technical_condition}
    \label{tc:subvector_klm_exogenous}
    Assume there exists a $\gamma^\star \in \BR^\mW$ such that 
    $$
    \gamma^\star = \gamma_0 + \left((Z \Pi_{ZW} + D \Pi_{DW})^T \ P_{P_{M_C [Z, D]} \tilde S(\beta_0, \gamma^\star, \delta_0)} W \right)^{-1} (Z \Pi_{ZW} + D \Pi_{DW})^T P_{P_{M_C [Z, D]} \tilde S(\beta_0, \gamma^\star, \delta_0)} \varepsilon,
    $$
    or, equivalently,
    $$
    (Z\Pi_{ZW} + D \Pi_{DW})^T P_{P_{M_C [Z, D]} \tilde S(\beta_0, \gamma^\star, \delta_0)} (\varepsilon + W(\gamma_0 - \gamma^\star)) = 0.
    $$    
\end{theoremEnd}
\begin{theoremEnd}[malte_intro,category=lm_exogenous]{theorem}%
    \label{prop:subvector_klm_test_statistic_chi_squared_exogenous}
    Consider \Cref{model:3} and assume that \Cref{ass:2} and \Cref{tc:subvector_klm_exogenous} hold.
    Under the null $\beta = \beta_0$ and $\delta = \delta_0$, under both strong and weak instrument asymptotics, the subvector Lagrange multiplier test statistic is bounded from above by a random variable that is asymptotically $\chi^2(\mX + \mD)$ distributed.
\end{theoremEnd}
\begin{proofEnd}
    \Cref{tc:subvector_klm_exogenous} is in terms of $M_C Z, M_C D, M_C X, M_C W$, and $M_C \varepsilon$.
    As $Z, C, D, X, W$, and $\varepsilon$ satisfy \Cref{ass:2}, so do $M_C Z, M_C D, M_C X, M_C W$, and $M_C \varepsilon$ by \Cref{lem:ass_2_reduction}.
    Using \Cref{lem:lm_without_c} we thus reduce to $\mC = 0$ and drop the $M_C$ in the following.

    We proceed in 4 steps:
    \begin{itemize}
        \item[1.] We show that $\gamma^\star$ from \Cref{tc:subvector_klm_exogenous} satisfies
        \begin{multline*}
        (y - X\beta_0 - W \gamma^\star - D \delta_0)^T P_{P_{[Z, D]} \tilde S(\beta_0, \gamma^\star, \delta_0)} (y - X\beta_0 - W \gamma^\star - D \delta_0) \\
        = \Psi_{\varepsilon^\star}^T (P_{\Psi_{\tilde S}} - P_{P_{\Psi_{\tilde S}} Q^{1/2} \left( \begin{smallmatrix} \Pi_{ZW} \\ \Pi_{DW} \end{smallmatrix} \right) }) \Psi_{\varepsilon^\star},
        \end{multline*}
        where $\Psi_{\varepsilon^\star} \in \BR^k$ and $\Psi_{\tilde S} \in \BR^{k \times m}$ are random variables.
        \item[2.] We show that under \Cref{ass:1}, the random variables' $\Psi_{\varepsilon^\star}$ and $\Psi_{\tilde V} := \Psi_{\tilde S} - \left( \left( \begin{smallmatrix} Z & D \end{smallmatrix} \right)^T \left( \begin{smallmatrix} Z & D \end{smallmatrix} \right) \right)^{1/2} \Pi $ rows are asymptotically centered Gaussian, asymptotically uncorrelated, and thus asymptotically independent.
        \item[3.] We argue, both for weak and strong instrument asymptotics, that
        $$\Psi_{\varepsilon^\star}^T (P_{\Psi_{\tilde S}} - P_{P_{\Psi_{\tilde S}} Q^{1/2} \left(\begin{smallmatrix} \Pi_{ZW} \\ \Pi_{DW} \end{smallmatrix} \right) }) \Psi_{\varepsilon^\star} \tod \sigma_{\varepsilon^\star}^2 \chi^2(\mX + \mD).$$
        \item[4.] We conclude as
        $$\LM(\beta_0) \leq (n - k - \mD) \frac{ \| P_{P_{[Z, D]} \tilde S(\beta_0, \gamma^\star, \delta_0)} (y - X \beta_0 - W \gamma^\star - D \delta_0) \|^2 }{ \| M_{[Z, D]} (y - X \beta_0 - W \gamma^\star - D \delta_0) \|^2}
        $$
        and $\frac{1}{n-k - \mD} \| M_{[Z, D]} (y - X \beta_0 - W \gamma^\star - D \delta_0) \|^2 \toP \sigma_{\varepsilon^\star}^2$.
    \end{itemize}
    In the following, write $\Pi_{W} = \left( \begin{smallmatrix} \Pi_{ZW} \\ \Pi_{DW} \end{smallmatrix} \right)$.
    \paragraph*{Step 1:}
    By \Cref{tc:subvector_klm}, there exists some $\gamma^\star \in \BR^{\mW}$ such that 
    $$
    ((Z \ D) \Pi_W)^T P_{P_{[Z, D]} \tilde S(\beta_0, \gamma^\star, \delta_0)} (\varepsilon + W(\gamma_0 - \gamma^\star)) = 0.
    $$
    Then
    \begin{align*}
        &\gamma^\star - \gamma_0 = ( ((Z \ D) \Pi_W)^T P_{P_{[Z, D]} \tilde S(\beta_0, \gamma^\star, \delta_0)} (Z \ D) \Pi_W)^{-1} \Pi_W^T (Z \ D)^T P_{P_{[Z, D]} \tilde S(\beta_0, \gamma^\star, \delta_0)} (\varepsilon + V_W (\gamma_0 - \gamma^\star)) \\
        &\ \Rightarrow P_{P_{[Z, D]} \tilde S(\beta_0, \gamma^\star, \delta_0)} (Z \ D) \Pi_W (\gamma_0 - \gamma^\star) = - P_{P_{P_{[Z, D]} \tilde S(\beta_0, \gamma^\star, \delta_0)} (Z \ D) \Pi_W} (\varepsilon + V_W (\gamma_0 - \gamma^\star))
    \end{align*}
    and thus
    \begin{multline*}
    P_{P_{[Z, D]} \tilde S(\beta_0, \gamma^\star, \delta_0)} (y - X\beta_0 - W \gamma^\star - D \delta_0) = P_{P_{[Z, D]} \tilde S(\beta_0, \gamma^\star, \delta_0)} (\varepsilon + V_W (\gamma_0 - \gamma^\star) + (Z \ D) \Pi_W (\gamma_0 - \gamma^\star)) \\
    = P_{P_{[Z, D]} \tilde S(\beta_0, \gamma^\star, \delta_0)} M_{P_{P_{[Z, D]} \tilde S(\beta_0, \gamma^\star, \delta_0)} (Z \ D) \Pi_W} (\varepsilon + V_W (\gamma_0 - \gamma^\star))\\
    =   M_{P_{P_{[Z, D]} \tilde S(\beta_0, \gamma^\star, \delta_0)} (Z \ D) \Pi_W} P_{P_{[Z, D]} \tilde S(\beta_0, \gamma^\star, \delta_0)} (\varepsilon + V_W (\gamma_0 - \gamma^\star))  \label{eq:1}\numberthis
    \end{multline*}
    as $P_{P_{[Z, D]} \tilde S(\beta_0, \gamma^\star, \delta_0)}$ and $M_{P_{P_{[Z, D]} \tilde S(\beta_0, \gamma^\star, \delta_0)} (Z \ D) \Pi_W} $ commute.
    Define $\varepsilon^\star := \varepsilon + V_W (\gamma_0 - \gamma^\star)$.
    Write
    \begin{align*}
        &\Psi_{\varepsilon^\star} := \Psi_{\varepsilon} + \Psi_{V_W} (\gamma_0 - \gamma^\star) = \left( ( Z \ D )^T ( Z \ D )\right)^{-1/2} ( Z \ D )^T \varepsilon^\star \\
        &\Psi_{\tilde V} := \begin{pmatrix} \Psi_{V_X} & \Psi_{V_W} \end{pmatrix} - \Psi_{\varepsilon^\star}\frac{{\varepsilon^\star}^T M_{[Z, D]} \begin{pmatrix}X & W \end{pmatrix}}{{\varepsilon^\star}^T M_{[Z, D]} \varepsilon^\star}, \text{ and} \\
        &\Psi_{\tilde S} := \begin{pmatrix} \Psi_{\tilde V} & 0 \end{pmatrix} + \left( ( Z \ D )^T ( Z \ D )\right)^{1/2} \begin{pmatrix} \Pi_{ZX} & \Pi_{DW} & 0 \\ \Pi_{DX} & \Pi_{DW} & \Id_\mD \end{pmatrix} \\
        &\hspace{2cm} = \left( ( Z \ D )^T ( Z \ D )\right)^{-1/2} ( Z \ D )^T \tilde S(\beta_0, \gamma^\star, \delta_0)
    \end{align*}
    Expand
    \begin{align*}
    P_{P_{[Z, D]} \tilde S(\beta_0, \gamma^\star, \delta_0)} &=
    ( Z \ D ) \left( ( Z \ D )^T ( Z \ D )\right)^{-1/2} \Psi_{\tilde S} (\Psi_{\tilde S}^T \Psi_{\tilde S})^{-1} \Psi_{\tilde S}^T \left( ( Z \ D )^T ( Z \ D )\right)^{-1/2} ( Z \ D )^T \\
    &= ( Z \ D ) \left( ( Z \ D )^T ( Z \ D )\right)^{-1/2} P_{\Psi_{\tilde S}} \left( ( Z \ D )^T ( Z \ D )\right)^{-1/2} ( Z \ D )^T
    \end{align*}
    such that
    $$
    P_{P_{[Z, D]} \tilde S(\beta_0, \gamma^\star, \delta_0)} (Z \ D) \Pi_W = Z \left( ( Z \ D )^T ( Z \ D )\right)^{-1/2} P_{\Psi_{\tilde S}} \left( ( Z \ D )^T ( Z \ D )\right)^{1/2} \Pi_W
    $$
    and
    \begin{multline*}
    M_{P_{P_{[Z, D]} \tilde S(\beta_0, \gamma^\star, \delta_0)} (Z \ D) \Pi_W} = (\Id_n - P_{P_{P_{[Z, D]} \tilde S(\beta_0, \gamma^\star, \delta_0)} (Z \ D) \Pi_W
    }) \\
    = \Id_n - (Z \ D) \left( ( Z \ D )^T ( Z \ D )\right)^{-1/2} P_{\Psi_{\tilde S}} \left( ( Z \ D )^T ( Z \ D )\right)^{1/2} \Pi_W \\
   \left(\Pi_W^T \left( ( Z \ D )^T ( Z \ D )\right)^{1/2} P_{\Psi_{\tilde S}}^T  P_{\Psi_{\tilde S}} \left( ( Z \ D )^T ( Z \ D )\right)^{1/2} \Pi_W\right)^{-1} \\
     \Pi_W^T \left( ( Z \ D )^T ( Z \ D )\right)^{1/2} P_{\Psi_{\tilde S}}^T \left( ( Z \ D )^T ( Z \ D )\right)^{-1/2} ( Z \ D )^T \\
    = \Id_n - Z \left( ( Z \ D )^T ( Z \ D )\right)^{-1/2} P_{P_{\Psi_{\tilde S}} \left( ( Z \ D )^T ( Z \ D )\right)^{1/2} \Pi_W} \left( ( Z \ D )^T ( Z \ D )\right)^{-1/2} ( Z \ D )^T.
    \end{multline*}
    We combine this to get
    \begin{align*}
        (y - &X \beta_0 - W \gamma^\star - D \delta_0)^T P_{P_{[Z, D]} \tilde S(\beta_0, \gamma^\star, \delta_0)} (y - X \beta_0 - W \gamma^\star - D \delta_0) \\       
        &\overset{(\ref{eq:1})}{=} {\varepsilon^\star}^T P_{P_{[Z, D]} \tilde S(\beta_0, \gamma^\star, \delta_0)} M_{P_{P_{[Z, D]} \tilde S(\beta_0, \gamma^\star, \delta_0)}(Z \ D) \Pi_W} P_{P_{[Z, D]} \tilde S(\beta_0, \gamma^\star, \delta_0)}  \varepsilon^\star \\
        &= \Psi_{\varepsilon^\star}^T 
        (P_{\Psi_{\tilde S}} - P_{P_{\Psi_{\tilde S}} \left( ( Z \ D )^T ( Z \ D )\right)^{1/2} \Pi_W}) \Psi_{\varepsilon^\star}.
        \end{align*}
    \paragraph*{Step 2:}
    By \Cref{ass:1}
    $$
    \vecop\begin{pmatrix} \Psi_{\varepsilon} & \Psi_{V_X} & \Psi_{V_W} \end{pmatrix} \tod \CN(0, \Omega \otimes \Id_k)
    $$
    and thus
    $$
    \vecop\begin{pmatrix} \Psi_{\varepsilon^\star} & \Psi_{V_X} & \Psi_{V_W} \end{pmatrix} \tod \CN(0, \underbrace{\begin{pmatrix} 1 & 0 & 0 \\ 0 & \Id_\mX & 0 \\ \gamma_0 - \gamma & 0 & \Id_\mW \end{pmatrix}^T \Omega \begin{pmatrix} 1 & 0 & 0 \\ 0 & \Id_\mX & 0 \\ \gamma_0 - \gamma & 0 & \Id_\mW \end{pmatrix}}_{=: \Omega^\star} \otimes \Id_k).
    $$
    By \Cref{ass:1} (a), $\frac{{\varepsilon^\star}^T M_{[Z, D]} (X \ W)}{{\varepsilon^\star}^T M_{[Z, D]} \varepsilon^\star} \toP \frac{\Omega^\star_{\varepsilon^\star, (V_X \ V_W)}}{\Omega^\star_{\varepsilon^\star}}$, implying that
    \begin{align*}
    \Cov(\Psi_{\varepsilon^\star}, \Psi_{\tilde V}) &= \Cov( \Psi_{\varepsilon^\star}, \begin{pmatrix} \Psi_{V_X} & \Psi_{V_W} \end{pmatrix} - \Psi_{\varepsilon^\star} \frac{{\varepsilon^\star}^T M_{[Z, D]} \begin{pmatrix} X & W \end{pmatrix}}{{\varepsilon^\star}^T M_{[Z, D]} \varepsilon^\star}) \\
    &\overset{\BP}{\to}
    \Cov(\Psi_{\varepsilon^\star}, \begin{pmatrix} \Psi_{V_X} & \Psi_{V_W} \end{pmatrix} - \Psi_{\varepsilon^\star} \frac{\Omega^\star_{\varepsilon^\star, (V_X \ V_W)}}{\Omega^\star_{\varepsilon^\star}}) \\
    &\toP \Omega^\star_{\varepsilon^\star, (V_X \ V_W)} - \Omega^\star_{\varepsilon^\star} \frac{\Omega^\star_{\varepsilon^\star, (V_X \ V_W)}}{\Omega^\star_{\varepsilon^\star}} = 0.
    \end{align*}
    Thus $\Psi_{\varepsilon^\star}$ and $\Psi_{\tilde V}$ are asymptotically uncorrelated and by their asymptotic normality thus asymptotically independent.

    \paragraph*{Step 3 for strong instrument asymptotics:} Here, $\Pi_W$ is constant and of full rank as $\Pi_{ZW}$ is of full rank.
    We calculate 
    \begin{align*}
    \plim \frac{1}{\sqrt{n}} \Psi_{\tilde S} &= \plim \frac{1}{\sqrt{n}} \left( ( Z \ D )^T ( Z \ D )\right)^{1/2} \begin{pmatrix} \Pi_{ZX} & \Pi_{ZW} & 0 \\ \Pi_{DX} & \Pi_{DW} & \Id_\mD \end{pmatrix} + \begin{pmatrix} \frac{1}{\sqrt{n}} \Psi_{\tilde V} & 0 \end{pmatrix} \\
    &= Q^{1/2} \begin{pmatrix} \Pi_{ZX} & \Pi_{ZW} & 0 \\ \Pi_{DX} & \Pi_{DW} & \Id_\mD \end{pmatrix}
    \end{align*}
    as $\Psi_{\tilde V}$ is asymptotically Gaussian and thus $\frac{1}{\sqrt{n}} \Psi_{\tilde V} \toP 0$.
    Then,
    \begin{align*}
    \Psi_{\varepsilon^\star}^T (P_{\Psi_{\tilde S}} &- P_{P_{\Psi_{\tilde S}} \left( ( Z \ D )^T ( Z \ D )\right)^{1/2} \Pi_W}) \Psi_{\varepsilon^\star} \\
    &= \Psi_{\varepsilon^\star}^T (P_{\frac{1}{\sqrt{n}}\Psi_{\tilde S}} - P_{P_{\frac{1}{\sqrt{n}}\Psi_{\tilde S}} \frac{1}{\sqrt{n}}\left( ( Z \ D )^T ( Z \ D )\right)^{1/2} \Pi_W}) \Psi_{\varepsilon^\star} \\
    &\toP \Psi_{\varepsilon^\star}^T (P_{Q^{1/2} \left( \begin{smallmatrix} \Pi_{ZX} & \Pi_{ZW} & 0 \\ \Pi_{DX} & \Pi_{DW} & \Id_\mD \end{smallmatrix} \right) } - P_{Q^{1/2}  \left( \begin{smallmatrix} \Pi_{ZW} \\ \Pi_{DW} \end{smallmatrix} \right) } ) \Psi_{\varepsilon^\star} \\
    &=
    \Psi_{\varepsilon^\star}^T P_{M_{Q^{1/2} \Pi_W} Q^{1/2} \left( \begin{smallmatrix} \Pi_{ZX} & 0 \\ \Pi_{DX} & \Id_\mD \end{smallmatrix} \right) } \Psi_{\varepsilon^\star} \tod \Var(\Psi_{\varepsilon^\star}) \cdot \chi^2(\mX + \mD),
    \end{align*}
    as $\Psi_{\varepsilon^\star} \tod \CN(0, \Omega^\star_{\varepsilon^\star} \cdot \Id_k)$ and $\rank\left( M_{Q^{1/2} \Pi_W}  Q^{1/2} \left( \begin{smallmatrix} \Pi_{ZX} & 0 \\ \Pi_{DX} & \Id_\mD \end{smallmatrix} \right) \right) = \mX + \mD$.

    \paragraph*{Step 3 for weak instrument asymptotics:} Here, $\sqrt{n} \begin{pmatrix} \Pi_{ZX} & \Pi_{ZW} \\ \Pi_{DX} & \Pi_{DW} \end{pmatrix} = \begin{pmatrix} \Pi_{ZX}^0 & \Pi_{ZW}^0 \\ \Pi_{DX}^0 & \Pi_{DW}^0 \end{pmatrix}$ is fixed.
    We calculate 
    \begin{align*}
    \plim \Psi_{\tilde S} &= \plim \left( ( Z \ D )^T ( Z \ D )\right)^{1/2} \begin{pmatrix} \Pi_{ZX} & \Pi_{ZW} & 0 \\ \Pi_{DX} & \Pi_{DW} & \Id_\mD \end{pmatrix} + \begin{pmatrix} \Psi_{\tilde V} & 0 \end{pmatrix} \\
    &= Q^{1/2} \begin{pmatrix} \Pi^0_{ZX} & \Pi^0_{ZW} & 0 \\ \Pi^0_{DX} & \Pi^0_{DW} & \Id_\mD \end{pmatrix} + \begin{pmatrix} \Psi_{\tilde V} & 0 \end{pmatrix}
    \end{align*}
    and
    \begin{align*}
        \Psi_{\varepsilon^\star}^T (P_{\Psi_{\tilde S}} -& P_{P_{\Psi_{\tilde S}} \left( ( Z \ D )^T ( Z \ D )\right)^{1/2} \Pi_W}) \Psi_{\varepsilon^\star}  \\
        &\toP \Psi_{\varepsilon^\star}^T (P_{Q^{1/2} \left(\begin{smallmatrix} \Pi^0_{ZX} & \Pi^0_{ZW} & 0 \\ \Pi^0_{DX} & \Pi^0_{DW} & \Id_\mD \end{smallmatrix} \right) + \left( \begin{smallmatrix} \Psi_{\tilde V} & 0 \end{smallmatrix} \right) } - P_{P_{ Q^{1/2} \left(\begin{smallmatrix} \Pi^0_{ZX} & \Pi^0_{ZW} & 0 \\ \Pi^0_{DX} & \Pi^0_{DW} & \Id_\mD \end{smallmatrix} \right) + \left( \begin{smallmatrix} \Psi_{\tilde V} & 0 \end{smallmatrix} \right) }Q^{1/2} \Pi^0_W})  \Psi_{\varepsilon^\star}
    \end{align*}
    We have $\Psi_{\varepsilon^\star} \tod \CN(0, \Omega_{\varepsilon^\star}^\star \cdot \Id_k)$, 
    $$\rank(P_{Q^{1/2} \left(\begin{smallmatrix} \Pi^0_{ZX} & \Pi^0_{ZW} & 0 \\ \Pi^0_{DX} & \Pi^0_{DW} & \Id_\mD \end{smallmatrix} \right) + \left( \begin{smallmatrix} \Psi_{\tilde V} & 0 \end{smallmatrix} \right) } - P_{P_{ Q^{1/2} \left(\begin{smallmatrix} \Pi^0_{ZX} & \Pi^0_{ZW} & 0 \\ \Pi^0_{DX} & \Pi^0_{DW} & \Id_\mD \end{smallmatrix} \right) + \left( \begin{smallmatrix} \Psi_{\tilde V} & 0 \end{smallmatrix} \right) }Q^{1/2} \Pi^0_W}) = \mX + \mD,$$
    and that, by Step 2, $\Psi_{\tilde V}$ and $\Psi_{\varepsilon^\star}$ are asymptotically independent.
    Thus, the above is asymptotically $\Omega^\star_{\varepsilon^\star} \cdot \chi^2(\mX + \mD)$-distributed.

    \paragraph*{Step 4:}
    Note that 
    $$
    \frac{1}{n-k - \mD} (y - X \beta_0 - W \gamma^\star - D \delta_0)^T M_{[Z, D]} (y - X \beta_0 - W \gamma^\star - D \delta_0) = \frac{1}{n-k - \mD} {\varepsilon^\star}^T M_{[Z, D]} \varepsilon^\star \overset{\BP}{\to} \Omega^\star_{\varepsilon^\star}.
    $$
    Thus
    \begin{align*}
    \LM(\beta_0) &= (n-k - \mD) \min_{\gamma \in \BR^\mW} \frac{(y - X \beta_0 - W \gamma )^T P_{P_{[Z, D]} \tilde S(\beta_0, \gamma, \delta_0)} (y - X \beta_0 - W \gamma)}{(y - X \beta_0 - W \gamma)^T M_{[Z, D]} (y - X \beta_0 - W \gamma)} \\
    & \leq
    (n-k - \mD) \frac{(y - X \beta_0 - W \gamma^\star)^T P_{P_{[Z, D]} \tilde S(\beta_0, \gamma^\star, \delta_0)} (y - X \beta_0 - W \gamma^\star)}{(y - X \beta_0 - W \gamma^\star)^T M_{[Z, D]} (y - X \beta_0 - W \gamma^\star)} \tod \chi^2(\mX + \mD)
    \end{align*}
    by Slutsky's Lemma.
\end{proofEnd}

%% file: theorems/prop_subvector_klm_test_statistic_chi_squared_many_instruments_exogenous.tex
\begin{theoremEnd}[malte_intro,category=lm_exogenous]{theorem}%
    \label{prop:subvector_klm_test_statistic_chi_squared_many_instruments_exogenous}
    Consider \Cref{model:3} and assume that \Cref{ass:2} holds as $k = k(n) \to \infty$ as $n \to\infty$.
    Here, one needs to reformulate \cref{ass:2} (b) with convergence of the difference between $\Psi$ and a Gaussian random variable to zero for all probabilities of Borel-measurable sets (see also \citeauthor{bentkus2003dependence}, \citeyear{bentkus2003dependence} and \citeauthor{chernozhukov2017central}, \citeyear{chernozhukov2017central}).
    Under the null $\beta = \beta_0$ and $\delta = \delta_0$, under both strong and weak instrument asymptotics, the subvector Lagrange multiplier test statistic is bounded from above by a random variable that is asymptotically $\chi^2(\mX + \mD)$ distributed.
\end{theoremEnd}
\begin{proofEnd}
    By \Cref{lem:lm_without_c} and \Cref{lem:ass_2_reduction} it suffices to show the result for $\mC = 0$.

    We proceed in 4 steps:
    \begin{itemize}
        \item[1.] We construct $\gamma^\star$ such that
        \begin{align}
            \nonumber
            \LM(\beta_0, \delta_0, \gamma^\star) :=&
            \frac{(y - X \beta_0 - W \gamma^\star - D \delta_0)^T P_{P_{[Z, D]} \tilde S(\beta_0, \gamma^\star, \delta_0)}(y - X \beta_0 - W \gamma^\star - D \delta_0)}{\frac{1}{n-k - \mD} (y - X \beta_0 - W \gamma^\star - D \delta_0)^T M_{[Z, D]} (y - X \beta_0 - W \gamma^\star - D \delta_0)}\\
            \label{eq:subvector_klm_test_statistic_chi_squared_many_instruments:0}        
            &= \frac{\varepsilon^T M_{P_{[Z, D]} \tilde W_0} P_{P_{[Z, D]} \tilde S(\beta_0, \gamma^\star, \delta_0)}  M_{P_{[Z, D]} \tilde W_0}  \varepsilon}{\frac{1}{n-k - \mD} (\varepsilon^T M_{[Z, D]} \varepsilon + \alpha^T \tilde W_0^T M_{[Z, D]} \tilde W_0 \alpha)} \\
            \label{eq:subvector_klm_test_statistic_chi_squared_many_instruments:1}
            &\toP \frac{\Psi_\varepsilon^T M_{\Psi_{\tilde W_0}} P_{\Psi_{\tilde S(\beta_0, \gamma^\star, \delta_0)}} M_{\Psi_{\tilde W_0}} \Psi_\varepsilon}{\sigma_\varepsilon^2 + \alpha^T (\Omega_{V_W} - \Omega_{V_W, \varepsilon} \Omega_{\varepsilon, V_W} / \sigma^2_\varepsilon) \alpha} \text{ as } n \to \infty,
        \end{align}
            where
        \begin{align*}
            \tilde W_0 &:= W - \varepsilon \frac{\varepsilon^T M_{[Z, D]} W}{\varepsilon^T M_{[Z, D]} \varepsilon}, \\
            \Psi_{\tilde W_0} &:= \left(\begin{pmatrix} Z & D \end{pmatrix}^T \begin{pmatrix} Z & D \end{pmatrix} \right)^{-1/2}\begin{pmatrix} Z & D \end{pmatrix}^T \tilde W_0, \\
            \Psi_{\tilde S(\beta_0, \gamma^\star, \delta_0 )} &:= \left(\begin{pmatrix} Z & D \end{pmatrix}^T \begin{pmatrix} Z & D \end{pmatrix} \right)^{-1/2} \begin{pmatrix} Z & D \end{pmatrix}^T \tilde S(\beta_0, \gamma^\star, \delta_0), \text{ and } \\
            \alpha &:= (\tilde W_0^T P_{[Z, D]} \tilde W_0)^{-1} \tilde W_0 P_{[Z, D]} \varepsilon.
        \end{align*}
        \item[2.] We calculate that
        \begin{align*}
        \Psi_{\tilde S(\beta_0, \gamma^\star, \delta_0)} &=: \begin{pmatrix} \Psi_{\tilde X} & \Psi_{\tilde W} & \Psi_{D} \end{pmatrix} \\
        &= \left(\begin{pmatrix} Z & D \end{pmatrix}^T \begin{pmatrix} Z & D \end{pmatrix} \right)^{1/2} \begin{pmatrix} \Pi_{ZX} & \Pi_{ZW} & 0 \\ \Pi_{DX} & \Pi_{DW} & \Id_{\mD} \end{pmatrix} + \begin{pmatrix} \Psi_{V_X} & \Psi_{V_W} & 0 \end{pmatrix} \\
        &\hspace{1cm}
        M_{\Psi_{\tilde W_0}} \Psi_\varepsilon \frac{
            (\varepsilon + \tilde W_0 \alpha )^T M_{[Z, D]} \begin{pmatrix} X & W & 0 \end{pmatrix}}{ \|M_{[Z, D]} (\varepsilon + \tilde W_0 \alpha ) \|^2 
        }
        \end{align*}
        whereas
        $$
        \Psi_{\tilde W_0} = \left(\begin{pmatrix} Z & D \end{pmatrix}^T \begin{pmatrix} Z & D \end{pmatrix} \right)^{1/2} \begin{pmatrix} \Pi_{ZW}  \\ \Pi_{DW} \end{pmatrix} + \Psi_{V_W} - \Psi_\varepsilon \frac{
            \varepsilon^T M_{[Z, D]} W }{ \|M_{[Z, D]} \varepsilon \|^2 
        }
        $$
        \item[3.] We show that, for strong instrument asymptotics, $\alpha \toP 0$ as $n \to \infty$ and $\Psi_{\tilde W_0}$ and $\Psi_{\tilde S(\beta_0, \gamma^\star, \delta_0)}$ converge to constants.
        This yields $\LM(\beta_0, \delta_0, \gamma^\star) \tod \chi^2(\mX + \mD)$.
        \item[4.] We show that, for weak instruments asymptotics, letting first $n \to \infty$ and then $k \to \infty$, we have $\alpha \toP 0$.
        Thus, $\Psi_{\tilde W} \toP \Psi_{\tilde W_0}$ and $\Psi_{\tilde X} \toP \Psi_{\tilde X_0}$, where
        $$
        \Psi_{\tilde X_0} := 
        \left( \left( \begin{matrix} Z & D \end{matrix} \right)^T  \left( \begin{matrix} Z & D \end{matrix} \right) \right)^{-1/2}  \left( \begin{matrix} Z & D \end{matrix} \right)^T \left(X - \varepsilon \frac{ \varepsilon^T M_{[Z, D]} X}{\varepsilon^T M_{[Z, D]} \varepsilon} \right).
        $$
        Also, $\Psi_{\tilde X_0}$ are $\Psi_{\tilde W_0}$ asymptotically independent of $\Psi_\varepsilon$.
        As $\frac{1}{\sqrt{n}} \Psi_D \toP Q^{1/2} \begin{pmatrix} 0 \\ \Id_{\mD} \end{pmatrix}$, this yields $\LM(\beta_0, \delta_0, \gamma^\star) \tod \chi^2(\mX + \mD)$ as $n \to \infty$ and $k \to \infty$.

    \end{itemize}

    \paragraph*{Step 1:}
    Define
    $\eta := \frac{\varepsilon^T M_{[Z, D]} W}{\varepsilon^T M_{[Z, D]} \varepsilon} \alpha$.
    This is not equal to $-1$ almost surely.
    Define $\gamma^\star := \gamma_0 + \frac{1}{1 + \eta} \alpha$.
    Note that $\tilde W_0$ ``uses'' the ground truth $\gamma_0$ and we will later see that the corresponding $\Psi_{\tilde W_0}$ and $\Psi_{\varepsilon}$ are asymptotically independent.
    For weak instrument asymptotics and finite $k$, this asymptotic independence does not hold for $\Psi_{\tilde S(\beta_0, \gamma^\star)}$, which ``uses'' $\gamma^\star$.

    We calculate
    \begin{align*}
    (1 + \eta) \cdot (y - X \beta_0 - W \gamma^\star - D \delta_0) &= (1 + \eta) \cdot (\varepsilon + W (\gamma_0 - \gamma^\star) ) \\
    &= \varepsilon + \varepsilon \cdot \frac{\varepsilon^T M_{[Z, D]} W}{\varepsilon^T M_{[Z, D]} \varepsilon} \alpha - W \alpha \\
    &= \varepsilon - \tilde W_0 \alpha\numberthis\label{eq:subvector_klm_test_statistic_chi_squared_many_instruments:2}
    \end{align*}
    and thus
    \begin{equation}
        \label{eq:subvector_klm_test_statistic_chi_squared_many_instruments:3}
        P_{[Z, D]} (y - X \beta_0 - W \gamma^\star - D \delta_0) = \frac{1}{1 + \eta} P_{[Z, D]} M_{P_{[Z, D]} \tilde W_0} \varepsilon.
    \end{equation}
    We plug \eqref{eq:subvector_klm_test_statistic_chi_squared_many_instruments:2} and \eqref{eq:subvector_klm_test_statistic_chi_squared_many_instruments:3} into the definition of the subvector Lagrange multiplier test statistic.
    Together with $\varepsilon^T M_{[Z, D]} \tilde W_0 = 0$, this yields \eqref{eq:subvector_klm_test_statistic_chi_squared_many_instruments:0}.

    Note that for any $A \in \BR^{n \times r}$, we have 
    \begin{align*}
    P_{[Z, D]} A &=
    \begin{pmatrix} Z & D \end{pmatrix} \left(\begin{pmatrix} Z & D \end{pmatrix}^T \begin{pmatrix} Z & D \end{pmatrix}\right)^{-1/2} \left(\begin{pmatrix} Z & D \end{pmatrix}^T \begin{pmatrix} Z & D \end{pmatrix}\right)^{-1/2} \begin{pmatrix} Z & D \end{pmatrix}^T A \\
    &= \begin{pmatrix} Z & D \end{pmatrix} \left(\begin{pmatrix} Z & D \end{pmatrix}^T \begin{pmatrix} Z & D \end{pmatrix}\right)^{-1/2} \Psi_A
    \end{align*}
    and thus 
    \begin{align*}
        P_{P_{[Z, D]} A} = \begin{pmatrix} Z & D \end{pmatrix} \left(\begin{pmatrix} Z & D \end{pmatrix}^T \begin{pmatrix} Z & D \end{pmatrix}\right)^{-1/2} P_{\Psi_A} \left(\begin{pmatrix} Z & D \end{pmatrix}^T \begin{pmatrix} Z & D \end{pmatrix}\right)^{-1/2} \begin{pmatrix} Z & D \end{pmatrix}^T
    \end{align*}
    and
    \begin{align*}
        P_{[Z, D]} \cdot M_{P_{[Z, D]} A} = \begin{pmatrix} Z & D \end{pmatrix} (\begin{pmatrix} Z & D \end{pmatrix}^T \begin{pmatrix} Z & D \end{pmatrix})^{-1/2} M_{\Psi_A} (\begin{pmatrix} Z & D \end{pmatrix}^T \begin{pmatrix} Z & D \end{pmatrix})^{-1/2} \begin{pmatrix} Z & D \end{pmatrix}^T.
    \end{align*}
    Applying this for $A = \varepsilon, A = \tilde W_0$, and $A = \tilde S(\beta_0, \gamma^\star, \delta_0)$ yields equality for the numerator of Equation \eqref{eq:subvector_klm_test_statistic_chi_squared_many_instruments:1}.

    For the denominator, note that by \Cref{ass:1} (a), 
    \begin{align*}
        &\frac{1}{n-k - \mC} \varepsilon^T M_Z \varepsilon \toP \sigma^2_\varepsilon, \\
        &\frac{1}{n-k - \mC} W^T M_{[Z, D]} W = \frac{1}{n-k - \mC} V_W^T M_{[Z, D]} V_W \toP \Omega_{V_W}, \text{ and }\\
        &\frac{1}{n-k - \mC} \varepsilon^T M_{[Z, D]} W = \frac{1}{n-k - \mC} \varepsilon^T M_{[Z, D]} V_W \toP \Omega_{\varepsilon, V_W}.
    \end{align*}
    Also, as $\tilde W_0^T M_{[Z, D]} \varepsilon = 0$ and thus
    $$
    \tilde W_0^T M_{[Z, D]} \tilde W_0 =
    W^T M_{[Z, D]} W + \frac{W^T M_{[Z, D]} \varepsilon}{\varepsilon^T M_{[Z, D]} \varepsilon} \varepsilon^T M_{[Z, D]} W
    \toP \Omega_{V_W} - \Omega_{V_W, \varepsilon} \Omega_{\varepsilon, V_W} / \sigma^2_\varepsilon.$$
    \paragraph*{Step 2:}
    Calculate
    \begin{align*}
        P_{[Z, D]} \tilde S(\beta_0, \gamma^\star, \delta_0) &= P_{[Z, D]} \begin{pmatrix} X & W & D \end{pmatrix} \\
        &\hspace{1.5cm} - P_{[Z, D]} (\varepsilon + W (\gamma_0 - \gamma^\star)) \cdot \frac{(\varepsilon + W (\gamma_0 - \gamma^\star) )^T M_{[Z, D]} \begin{pmatrix} X & W & 0 \end{pmatrix}}{ \| M_{[Z, D]} (\varepsilon + W (\gamma_0 - \gamma^\star) ) \|^2} \\
        &\overset{\eqref{eq:subvector_klm_test_statistic_chi_squared_many_instruments:3}}{=}P_{[Z, D]} \begin{pmatrix} X & W & D \end{pmatrix} - P_{[Z, D]} M_{P_{[Z, D]} \tilde W_0} \varepsilon \cdot \frac{
            (\varepsilon + \tilde W_0 \alpha )^T M_{[Z, D]} \begin{pmatrix} X & W & 0 \end{pmatrix}}{ \|M_{[Z, D]} (\varepsilon + \tilde W_0 \alpha ) \|^2 }
    \end{align*}
    (the $1 + \eta$ cancel out) and thus
    \begin{align*}
            \Psi_{\tilde S(\beta_0, \gamma^\star, \delta_0)} &= \left(\begin{pmatrix} Z & D \end{pmatrix}^T \begin{pmatrix} Z & D \end{pmatrix} \right)^{1/2} \begin{pmatrix} \Pi_{ZX} & \Pi_{ZW} & 0 \\ \Pi_{DX} & \Pi_{DW} & \Id_{\mD} \end{pmatrix} + \begin{pmatrix} \Psi_{V_X} & \Psi_{V_W} & 0 \end{pmatrix} \\
            &\hspace{3cm} - M_{\Psi_{\tilde W_0}} \Psi_\varepsilon \frac{
            (\varepsilon + \tilde W_0 \alpha )^T M_{[Z, D]} \begin{pmatrix} X & W & 0 \end{pmatrix}}{ \|M_{[Z, D]} (\varepsilon + \tilde W_0 \alpha ) \|^2 }
    \end{align*}
    The equation
    \begin{align*}
        \Psi_{\tilde W_0} &= \left(\begin{pmatrix} Z & D \end{pmatrix}^T \begin{pmatrix} Z & D \end{pmatrix} \right)^{-1/2} \begin{pmatrix} Z & D \end{pmatrix}^T \left(W - \varepsilon \cdot \frac{ \varepsilon^T M_{[Z, D]} W } {\varepsilon^T M_{[Z, D]} \varepsilon} \right) \\
        &= \left(\begin{pmatrix} Z & D \end{pmatrix}^T \begin{pmatrix} Z & D \end{pmatrix} \right)^{1/2} \begin{pmatrix} \Pi_{ZW}  \\ \Pi_{DW} \end{pmatrix} + \Psi_{V_W} - \Psi_\varepsilon \frac{
            \varepsilon^T M_{[Z, D]} W }{ \|M_{[Z, D]} \varepsilon \|^2 
        }
    \end{align*}
    follows from the definition of $\tilde W_0$.

    Now, as long as $\alpha \neq 0$, we have
    \begin{align*}
    \plim_{n \to\infty} M_{\Psi_{\tilde W_0}} \Psi_{\tilde S(\beta_0, \gamma^\star, \delta_0)} &\neq 
    \plim_{n \to \infty} M_{\Psi_{\tilde W_0}} \begin{pmatrix} \Psi_{\tilde X_0} & \Psi_{\tilde W_0} & Q^{1/2} \begin{pmatrix} 0 \\ \Id_\mD \end{pmatrix} \end{pmatrix} \\
    &= \plim_{n \to \infty} \begin{pmatrix} M_{\Psi_{\tilde W_0}} \Psi_{\tilde X_0} & 0 & M_{\Psi_{\tilde W_0}} Q^{1/2} \begin{pmatrix} 0 \\ \Id_\mD \end{pmatrix} \end{pmatrix}.
    \end{align*}
    We will show that for strong instrument asymptotics $\alpha \toP 0$ as $n \to \infty$ and that for weak instrument asymptotics $\alpha \toP 0$ as $n \to \infty$ and $k \to \infty$.

    \paragraph*{Step 3:}
    Under strong instrument asymptotics,
    \begin{align*}
    \frac{1}{\sqrt{n}} \Psi_{\tilde W_0} &= \frac{1}{\sqrt{n}} \left( \begin{pmatrix} Z & D \end{pmatrix}^T \begin{pmatrix} Z & D \end{pmatrix} \right)^{1/2} \begin{pmatrix} \Pi_{ZW} \\ \Pi_{DW} \end{pmatrix}
    + \frac{1}{\sqrt{n}} \Psi_{V_W} 
    - \frac{1}{\sqrt{n}} \Psi_{\varepsilon} \frac{\varepsilon^T M_{[Z, D]}W}{\varepsilon^T M_{[Z, D]} \varepsilon} \\
    &\toP Q^{1/2} \begin{pmatrix} \Pi_{ZW} \\ \Pi_{DW} \end{pmatrix} \text{ as }n \to \infty \numberthis\label{eq:subvector_klm_test_statistic_chi_squared_many_instruments:4}
    \end{align*}
    Thus,
    $$
    \sqrt{n} \cdot \alpha = \left(\frac{1}{\sqrt{n}}\Psi_{\tilde W_0}^T \frac{1}{\sqrt{n}}\Psi_{\tilde W_0}\right)^{-1} \frac{1}{\sqrt{n}}\Psi_{\tilde W_0}^T \Psi_\varepsilon \toP \left(\begin{pmatrix} \Pi_{ZW} \\ \Pi_{DW} \end{pmatrix}^T Q \begin{pmatrix} \Pi_{ZW} \\ \Pi_{DW} \end{pmatrix}\right)^{-1} \begin{pmatrix} \Pi_{ZW} \\ \Pi_{DW} \end{pmatrix}^T Q^{1/2} \Psi_\varepsilon
    $$
    and $\alpha \toP 0$ and
    \begin{align*}
    \plim_{n \to\infty} M_{\Psi_{\tilde W_0}} \Psi_{\tilde S(\beta_0, \gamma^\star, \delta_0)}    &= \plim_{n \to \infty} \begin{pmatrix} M_{\Psi_{\tilde W_0}} \Psi_{\tilde X_0} & 0 & M_{\Psi_{\tilde W_0}} Q^{1/2} \begin{pmatrix} 0 \\ \Id_\mD \end{pmatrix} \end{pmatrix}.
    \end{align*}
    Similarly to Equation \eqref{eq:subvector_klm_test_statistic_chi_squared_many_instruments:4} we have 
    \begin{align*}
    \frac{1}{\sqrt{n}} \Psi_{\tilde{X}_0}  
    &= \frac{1}{\sqrt{n}} \left(\begin{pmatrix} Z & D \end{pmatrix}^T \begin{pmatrix} Z & D \end{pmatrix} \right)^{1/2}
    \begin{pmatrix} \Pi_{ZX} \\ 0\end{pmatrix}
    + \frac{1}{\sqrt{n}} \Psi_{V_X} - \frac{1}{\sqrt{n}} \Psi_{\varepsilon} \frac{\varepsilon^T M_{[Z, D]} X}{\varepsilon^T M_{[Z, D]} \varepsilon} \toP Q^{1/2} \begin{pmatrix} \Pi_{ZX} \\ \Pi_{DX} \end{pmatrix}.
    \end{align*}

    \noindent Applying \Cref{lem:chain_projections} yields
    \begin{align*}
    &\Psi_{\varepsilon}^T M_{\Psi_{\tilde W_0}} P_{\Psi_{\tilde S(\beta_0, \gamma^\star)}} M_{\Psi_{\tilde W_0}} \Psi_\varepsilon \\
    &\hspace{2cm} 
    \toP \Psi_\varepsilon^T M_{Q^{1/2} {\left( \begin{smallmatrix} \Pi_{ZW} \\ \Pi_{DW} \end{smallmatrix} \right) } }
    P_{Q^{1/2} \left( \begin{smallmatrix} \Pi_{ZX} & \Pi_{ZW} & 0 \\ \Pi_{DX} &  \Pi_{DW} & \Id_\mD \end{smallmatrix} \right) }
        M_{Q^{1/2} \left( \begin{smallmatrix} \Pi_{ZW} \\ \Pi_{DW} \end{smallmatrix} \right) } \Psi_\varepsilon \\
        &\hspace{2cm} \overset{\Cref{lem:chain_projections}}{=} 
    \Psi_\varepsilon^T 
        P_{ M_{Q^{1/2} {\left( \begin{smallmatrix} \Pi_{ZW} \\ \Pi_{DW} \end{smallmatrix} \right) } } Q^{1/2} \left( \begin{smallmatrix} \Pi_{ZX} & 0 & 0 \\ \Pi_{DX} & 0 &  \Id_\mD \end{smallmatrix} \right) } 
    \Psi_\varepsilon 
    \end{align*}
    As $\rank\left(M_{Q^{1/2} {\left( \begin{smallmatrix} \Pi_{ZW} \\ \Pi_{DW} \end{smallmatrix} \right) } } Q^{1/2} \left( \begin{smallmatrix} \Pi_{ZX} & 0 & 0 \\ \Pi_{DX} & 0 &  \Id_\mD \end{smallmatrix} \right) \right)= \mX + \mD$ and $\Psi_\varepsilon$ is asymptotically $\CN(0, \sigma_\varepsilon \cdot \Id_k)$ by \Cref{ass:1} (b), this converges in distribution to a random variable that is distributed as $ \sigma^2_\varepsilon \chi^2(\mX + \mD)$.
    Also 
    $$\frac{1}{n-k - \mD} (y - X \beta_0 - W \gamma^\star - D \delta_0)^T M_{[Z, D]} (y - X \beta_0 - W \gamma^\star - \delta_0) \toP \sigma^2_\varepsilon,$$ and thus
    \begin{align*}
    LM(\beta_0, \gamma^\star, \delta_0) \tod \chi^2(\mX + \mD)
    \end{align*}
    This is an upper bound to $\LM(\beta_0, \delta_0)$ due to the minimization over $\gamma$ and thus, under strong instrument asymptotics, $\LM(\beta_0, \delta_0)$ converges to a random variable that is bounded from above by a $\chi^2(\mX + \mD)$ distributed random variable.
    Note that this holds for any fixed $k$ and does not require $k \to \infty$.

    \paragraph*{Step 4:}
    Now assume weak instrument asymptotics, where $\sqrt{n} \begin{pmatrix} \Pi_{ZX} & \Pi_{ZW} \\ \Pi_{DX} & \Pi_{DW} \end{pmatrix} =: \begin{pmatrix} \Pi^0_{ZX} & \Pi^0_{ZW} \\ \Pi^0_{DX} & \Pi^0_{DW} \end{pmatrix}$ is fixed.

    Calculate
    \begin{multline*}
        \begin{pmatrix} \Psi_{\tilde X_0} & \Psi_{\tilde W_0} & \Psi_\varepsilon \end{pmatrix} = \frac{1}{\sqrt{n}} \left(\begin{pmatrix} Z & D \end{pmatrix}^T \begin{pmatrix} Z & D \end{pmatrix} \right)^{1/2}  \begin{pmatrix} \Pi^0_{ZX} & \Pi^0_{ZW} & 0 \\ \Pi^0_{DX} & \Pi^0_{DW} & 0 \end{pmatrix} \\
        + \begin{pmatrix} \Psi_{V_X} & \Psi_{V_W} & \Psi_\varepsilon \end{pmatrix} - \Psi_\varepsilon \begin{pmatrix} \frac{\varepsilon^T M_{[Z, D]} X}{\varepsilon^T M_{[Z, D]} \varepsilon} & \frac{\varepsilon^T M_{[Z, D]} W}{\varepsilon^T M_{[Z, D]} \varepsilon} & 0 \end{pmatrix} \\
        \toP Q^{1/2} \begin{pmatrix} \Pi^0_{ZX} & \Pi^0_{ZW} & 0 \\ \Pi^0_{DX} & \Pi^0_{DW} & 0 \end{pmatrix} + \begin{pmatrix} \Psi_{V_X} & \Psi_{V_W} & \Psi_\varepsilon \end{pmatrix} - \Psi_\varepsilon \begin{pmatrix} \Omega_{\varepsilon, V_X} / \sigma^2_\varepsilon & \Omega_{\varepsilon, V_W} / \sigma^2_\varepsilon & 0 \end{pmatrix}.
    \end{multline*}
    Thus, $\Psi_{\tilde X_0}$, $\Psi_{\tilde W_0}$, and $\Psi_\varepsilon$ are asymptotically jointly Gaussian as $n \to \infty$ by \Cref{ass:2} (b). Also,
    \begin{align*}
        \Cov(\Psi_\varepsilon, \Psi_{\tilde X_0} ) &= \left( \Omega_{\varepsilon, V_X} - \sigma^2_\varepsilon \cdot \frac{\varepsilon^T M_{[Z, D]} X}{\varepsilon^T M_{[Z, D]} \varepsilon} \right) \otimes \Id_{k + \mD} \toP \left( \Omega_{\varepsilon, V_X} - \sigma^2_\varepsilon \cdot \frac{\Omega_{\varepsilon, V_X}}{\sigma_\varepsilon^2} \right) \otimes \Id_{k + \mD} = 0
    \end{align*}
    and
    \begin{align*}
        \Cov(\Psi_\varepsilon, \Psi_{\tilde W_0} ) &= \left( \Omega_{\varepsilon, V_W} - \sigma^2_\varepsilon \cdot \frac{\varepsilon^T M_{[Z, D]} W}{\varepsilon^T M_{[Z, D]} \varepsilon} \right) \otimes \Id_{k + \mD} \toP \left( \Omega_{\varepsilon, V_W} - \sigma^2_\varepsilon \cdot \frac{\Omega_{\varepsilon, V_W}}{\sigma_\varepsilon^2} \right) \otimes \Id_{k + \mD} = 0.
    \end{align*}
    By asymptotic Gaussianity, they are also asymptotically independent.

    Thus, $\alpha = (\tilde W_0^T P_{[Z, D]} \tilde W_0)^{-1}\tilde W_0^T P_{[Z, D]} \varepsilon = (\Psi_{\tilde W_0}^T \Psi_{\tilde W_0})^{-1} \Psi_{\tilde W_0}^T \Psi_\varepsilon$ is the linear regression coefficient regressing $\Psi_\varepsilon$ on $\Psi_{\tilde W_0}$, which are asymptotically jointly Gaussian and independent.
    Thus, $\alpha$ is consistent for $0$ and $\alpha = O(k^{-1/2}) \toP 0$ as $k \to \infty$.
    Therefore, letting $k\to\infty$ after $n \to \infty$ yields
    $\Psi_W \toP \Psi_{W_0} \text{ and } \Psi_X \toP \Psi_{X_0}.$
    The projection $P_{\Psi_{\tilde S(\beta_0, \gamma^\star, \delta_0)}}$ only depends on the column space of $\Psi_{S(\beta_0, \gamma^\star, \delta_0)}$, which is independent of scaling.
    Thus
    $$
    P_{\Psi_{\tilde S(\beta_0, \gamma^\star, \delta_0)}} = P_{\begin{pmatrix} \Psi_{\tilde X} & \Psi_{\tilde W} & \Psi_D \end{pmatrix}} = P_{\begin{pmatrix} \Psi_{\tilde X} & \Psi_{\tilde W} & \frac{1}{\sqrt{n}} \Psi_D \end{pmatrix}} \toP P_{\begin{pmatrix} \Psi_{\tilde X_0} & \Psi_{\tilde W_0} & Q^{1/2} \left(\begin{smallmatrix} 0 \\ \Id_\mD \end{smallmatrix} \right) \end{pmatrix}}
    $$
    and
    $$
    M_{\Psi_{\tilde W_0}} P_{\Psi_{\tilde S(\beta_0, \gamma^\star, \delta_0)}} \toP P_{\begin{pmatrix} M_{\Psi_{\tilde W_0}} \Psi_{\tilde X_0} & 0 & M_{\Psi_{\tilde W_0}} Q^{1/2} \left(\begin{smallmatrix} 0 \\ \Id_\mD \end{smallmatrix} \right) \end{pmatrix}}
    $$
    as $n \to \infty$ and $k \to \infty$.
    Thus, the numerator of \eqref{eq:subvector_klm_test_statistic_chi_squared_many_instruments:1} 
    \begin{align*}
            \Psi_\varepsilon^T M_{\Psi_{\tilde W_0}} P_{\Psi_{\tilde S(\beta_0, \gamma^\star, \delta_0)}} M_{\Psi_{\tilde W_0}} \Psi_\varepsilon \toP \Psi_\varepsilon^T P_{\begin{pmatrix} M_{\Psi_{\tilde W_0}} \Psi_{\tilde X_0} & 0 & M_{\Psi_{\tilde W_0}} Q^{1/2} \left(\begin{smallmatrix} 0 \\ \Id_\mD \end{smallmatrix} \right) \end{pmatrix}} \Psi_\varepsilon
    \end{align*}
    as $n \to \infty$ and $k \to \infty$.
    Conditionally on $\Psi_{\tilde X_0}$ and $\Psi_{\tilde W_0}$, by their asymptotic independence with $\Psi_\varepsilon$, asymptotically $\Psi_\varepsilon \tod \CN(0, \sigma_\varepsilon^2 \cdot \Id_{k + \mD})$ and
    $$
    \Psi_\varepsilon^T P_{\begin{pmatrix} M_{\Psi_{\tilde W_0}} \Psi_{\tilde X_0} & 0 & M_{\Psi_{\tilde W_0}} Q^{1/2} \left(\begin{smallmatrix} 0 \\ \Id_\mD \end{smallmatrix} \right) \end{pmatrix}} \Psi_\varepsilon \tod \sigma_\varepsilon^2 \cdot \chi^2(\mX + \mD),
    $$
    as
    $$
    \rank\left( \begin{pmatrix} M_{\Psi_{\tilde W_0}} \Psi_{\tilde X_0} & 0 & M_{\Psi_{\tilde W_0}} Q^{1/2} \left(\begin{smallmatrix} 0 \\ \Id_\mD \end{smallmatrix} \right) \end{pmatrix} \right) = \mX + \mD
    $$
    As this holds conditionally on $\Psi_{\tilde X_0}$ and $\Psi_{\tilde W_0}$ for all values of $\Psi_{\tilde X_0}$ and $\Psi_{\tilde W_0}$,this also holds unconditionally.

    On the other hand, the denominator of \eqref{eq:subvector_klm_test_statistic_chi_squared_many_instruments:1}
    $$
    \sigma_\varepsilon^2 + \alpha^T (\Omega_{V_W} - \Omega_{V_W, \varepsilon} \Omega_{\varepsilon, V_W} / \sigma^2_\varepsilon) \alpha \toP \sigma_\varepsilon^2
    $$
    as $\alpha \toP 0$ as $k \to \infty$.
    Thus, $\LM(\beta_0, \gamma^\star, \delta_0) \tod \chi^2(\mX + \mD)$ as $n \to \infty$ and $k \to \infty$.
    Due to the minimization over $\gamma$ in $\LM(\beta_0, \delta_0)$, this is an upper bound to $\LM(\beta_0, \delta_0)$ and thus $\LM(\beta_0, \delta_0)$ is, asymptotically, as $n \to \infty$ and $k \to \infty$, upper bounded by a $\chi^2(\mX + \mD)$ distributed random variable.

\end{proofEnd}

%% file: sec_auxiliary_results.tex
\section{Proofs}
\label{sec:proofs}
\subsection{Auxiliary results}
\input{theorems/lem_kappa_pos_definite.tex}
\input{theorems/lem_projected_quadrics.tex}
\input{theorems/lem_chain_projections.tex}
\input{theorems/prop_lm_derivative.tex}
\subsection{Proofs of \Cref{sec:main_results:subvector_lagrange_multiplier} and Appendix \ref{sec:exogenous_variables}}
\label{app:proofs:subvector_lagrange_multiplier}
We prove the more general results of Appendix \ref{sec:exogenous_variables}.
\Cref{prop:subvector_klm_test_statistic_chi_squared_many_instruments_exogenous} and \Cref{prop:subvector_klm_test_statistic_chi_squared_exogenous} directly imply \Cref{prop:subvector_klm_test_statistic_chi_squared_many_instruments} and \Cref{prop:subvector_klm_test_statistic_chi_squared} or \Cref{sec:main_results:subvector_lagrange_multiplier} by setting $\mC = \mD = 0$.

\defsubvectorklmteststatisticexogenous*
\printProofs[lm_exogenous]
\subsection{Proofs of \Cref{sec:main_results:closed_form_subvector_confidence_sets}}
\label{app:proofs:closed_form_subvector_confidence_sets}
\subvectorandersonrubinteststatistic*
\input{theorems/def_wald_test_statistic.tex}
\input{theorems/def_lr_test_statistic.tex}
\input{theorems/prop_cis_closed_forms.tex}

\tctwo*
\input{theorems/prop_cis_bounded.tex}
\input{theorems/cor_ar_bounded_iff_rank_test_rejects.tex}
\printProofs[propinversearequaltowald]

%% file: theorems/lem_kappa_pos_definite.tex
\begin{theoremEnd}[malte_all_end]{lemma}
    \label{lem:kappa_pos_definite}
    Assume that $M_Z X$ is of full column rank.
    Let $\lambda_1 := \lambdamin{ (X^T M_Z X)^{-1} X^T X} = \lambdamin{ (X^T M_Z X)^{-1} X^T P_Z X} + 1$.
    \begin{itemize}
    \item Then $M(\kappa) := X^T (\kappa P_Z + (1 - \kappa) \Id) X$ is positive definite if and only if $\kappa < \lambda_1$.
    \item The matrix $M(\kappa)$ is singular positive semi-definite if and only if $\kappa = \lambda_1$.
    \item The matrix $M(\kappa)$ has at least one negative eigenvalue if and only if $\kappa > \lambda_1$.
    \end{itemize}
\end{theoremEnd}
\begin{proofEnd}
    Let $v$ be an eigenvector corresponding to an eigenvalue $\lambda$ of $(X^T M_Z X)^{-1} X^T X$.
    Then
    $$
        (X^T M_Z X)^{-1} X^T X v = \lambda v \Leftrightarrow \lambda X^T M_Z X v = X^T X v \Leftrightarrow X^T (\Id - \lambda M_Z) X v = 0.
    $$
    That is, the eigenvalues of $(X^T M_Z X)^{-1} X^T X v$ correspond to $\kappa$ such that $X^T (\Id - \kappa M_Z) X v = M(\kappa) v =  0$.
    Also, as $X^T M_Z X \preceq X^T X$, we have $\lambda \geq 1$.

    Note that, as $M_Z X$ is of full column rank, $X^T M_Z X$ is positive definite.
    Write $M(\kappa) = X^T (P_Z + (1 - \lambda) M_Z) X + (\lambda - \kappa) X^T M_Z X$ such that for any eigenvector, eigenvalue pair $(v, \lambda)$ of $(X^T M_Z X)^{-1} X^T X$, we have $v^T M(\kappa) v = (\lambda - \kappa) v^T X^T M_Z Xv$.
    \begin{itemize}
    \item If $\kappa < \lambda_1 \leq \lambda$, then $v^T M(\kappa) v > 0$ by positive definiteness of $X^T M_Z X$ and thus $M(\kappa)$ is positive definite.
    \item If $\kappa \leq \lambda_1 \leq \lambda$, then $v^T M(\kappa) v_1 \geq 0$ and thus $M(\kappa)$ is positive semi-definite.
    Let $v_1$ be an eigenvector corresponding to $\lambda_1$. Then $v_1^T M(\lambda_1) v_1 = 0$ and thus $M(\lambda_1)$ is not positive definite.
    \item Let $v_1$ be the eigenvector corresponding to $\lambda_1$.
    If $\kappa > \lambda_1$, then $v^T M(\kappa) v = (\lambda_1 - \kappa) v^T X^T M_Z X v < 0$ and thus $M(\kappa)$ has at least one negative eigenvalue.
    \end{itemize}
    The other direction follows as we exhaust all possible cases.
\end{proofEnd}

%% file: theorems/lem_projected_quadrics.tex
\begin{theoremEnd}[malte_all_end]{lemma}
    \label{lem:projected_quadric}
    Let $n\geq m$, $A = \begin{pmatrix} A_{11} & A_{12} \\ A_{21} & A_{22} \end{pmatrix} \in \BR^{n \times n}$ be symmetric with $A_{11} \in \BR^{m \times m}$, $A_{12} = A_{21}^T \in \BR^{m \times (n-m)}$, and $A_{22} \in \BR^{(n-m) \times (n-m)}$.
    Let $z_0 \in \BR^n$, $c \in \BR$, and let $B \in \BR^{m \times n}$ have ones on the diagonal and zeros elsewhere.
    Then
    \begin{align*}
        &\ \  \{ x \in \BR^{m} \mid \min_{y \in \BR^{n - m}} \left(\begin{pmatrix} x \\ y \end{pmatrix} - z_0 \right)^T A \left(\begin{pmatrix} x \\ y \end{pmatrix} - z_0 \right) \leq c \} \\
        &= \begin{cases}
            \ \ \{ x \in \BR^m \mid x^T (A_{11} - A_{12} A_{22}^{\dagger} A_{21}) x \leq c \} + B z_0 & \ \ \text{if } \lambdamin{A_{22}} \geq 0 \\
            \ \ \BR^{m} & \ \ \text{if } \lambdamin{A_{22}} < 0
            \end{cases} \numberthis \label{eq:lem:projected_quadric:1} \\
        &= \begin{cases}
            \ \ \{ x \in \BR^m \mid (x - B z_0)^T (B A^{\dagger} B^T)^{\dagger} (x - B z_0) \leq c \} &  \ \ \text{if } \lambdamin{A_{22}} \geq 0 \\
            \ \ \BR^{m} & \ \  \text{if } \lambdamin{A_{22}} < 0,
            \end{cases} \numberthis \label{eq:lem:projected_quadric:2}
    \end{align*}
    where $\dagger$ denotes the Moore-Penrose pseudo-inverse.
\end{theoremEnd}
\begin{proofEnd}
        First, note that
        \begin{align*}
            &\{ x \in \BR^{m} \mid \min_{y \in \BR^{n - m}} \left(\begin{pmatrix} x \\ y \end{pmatrix} - z_0 \right)^T A \left(\begin{pmatrix} x \\ y \end{pmatrix} - z_0 \right) \leq c \}
            \\
            &\hspace{2cm} = \{ \tilde x \in \BR^{m} \mid \min_{y \in \BR^{n - m}} \begin{pmatrix} \tilde x \\ y \end{pmatrix}^T A \begin{pmatrix} \tilde x \\ y \end{pmatrix} \leq c \} + B z_0
             \numberthis \label{eq:lem:projected_quadric:2:0}
        \end{align*}
        via a change of variables $\tilde x = x - B z_0$.
        \paragraph*{Case 1: $\lambdamin{A_{22}} < 0.$}
        Let $v \in \BR^{n - m}$ be an eigenvector of $A_{22}$ with eigenvalue $\lambda < 0$.
        Then,
        $$
            \inf_{y \in \BR^{n-m}} \begin{pmatrix} x \\ y \end{pmatrix}^T A \begin{pmatrix} x \\ y \end{pmatrix} \leq \inf_{t \in \BR} \begin{pmatrix} x \\ tv \end{pmatrix}^T A \begin{pmatrix} x \\ tv \end{pmatrix} = \inf_{t \in \BR} x^T A_{11} x + 2 x^T A_{12} v t + \lambda \| v \|^2 t^2 = -\infty
        $$
        Thus, the set is $\BR^m$.

        \paragraph*{Case 2: $\lambdamin{A_{22}} \geq 0.$}
        Then, for any $x\in\BR^n$, the function $y \mapsto \begin{pmatrix} x \\ y \end{pmatrix}^T A \begin{pmatrix} x \\ y \end{pmatrix}$ is convex and any minimizer $y^\star(x)$ satisfies 
        \begin{align}
            \label{eq:lem:projected_quadric:2:1}
            0 &= \frac{\dd}{\dd y} \left. \left[ \begin{pmatrix} x \\ y \end{pmatrix}^T A \begin{pmatrix} x \\ y \end{pmatrix}\right] \right|_{y = y^\star(x)} = A_{21} x + A_{12}^T x + 2 A_{22} y^\star(x) \\
            \label{eq:lem:projected_quadric:2:2}
            & \Rightarrow y^\star(x) = - A_{22}^\dagger A_{12}^T x.
        \end{align}
        such that
        \begin{align*}
            \min_y \begin{pmatrix} x \\ y \end{pmatrix}^T A \begin{pmatrix} x \\ y \end{pmatrix} &= 
            x^T A_{11} x +  x^T A_{12} y^\star(x) + y^\star(x)^T A_{21} x +y^\star(x)^T A_{22} y^\star(x) \\
            &=x^T A_{11} x +  y^\star(x)^T (A_{12}^T x + A_{21} x + A_{22}  y^\star(x) )\\
            &\overset{\eqref{eq:lem:projected_quadric:2:1}}{=} x^T A_{11} x -  y^\star(x)^T A_{22} y^{\star}(x) \\
            &\overset{\eqref{eq:lem:projected_quadric:2:2}}{=} x^T A_{11} x - x^T A_{12} A_{22}^\dagger A_{21} x
        \end{align*}
        Together with \eqref{eq:lem:projected_quadric:2:0}, thus proves \eqref{eq:lem:projected_quadric:1}.
        By the formula for the pseudo inverse of a block matrix in \citep{rohde1965generalized} (here, $A_{11} - A_{12} A_{22}^\dagger A_{21}$ is the Schur complement $A / A_{22}$), we have $A_{11} - A_{12} A_{22}^\dagger A_{21} = (B A^\dagger B)^\dagger$.
        Finally, \eqref{eq:lem:projected_quadric:2} follows from another change of variables $\tilde x = x - B z_0$.
\end{proofEnd}

%% file: theorems/lem_chain_projections.tex
\begin{theoremEnd}[malte_all_end]{lemma}
    \label{lem:chain_projections}
    For $p \geq q$ and $p \geq r$ let $A \in \BR^{p \times q}$ and $B \in \BR^{p \times r}$. Then $P_{\left( \begin{smallmatrix} A & B \end{smallmatrix} \right) } M_A = P_{M_A B}$.
\end{theoremEnd}
\begin{proofEnd}
    First, note that $P_{\left( \begin{smallmatrix} A & B \end{smallmatrix} \right)} = P_{\left( \begin{smallmatrix} A & M_A B \end{smallmatrix} \right)}$ as the column space spanned by $\left( \begin{smallmatrix} A & B \end{smallmatrix} \right)$ is the same as the column space spanned by $\left( \begin{smallmatrix} A & M_A B \end{smallmatrix} \right)$.
    Then, as $A$ and $M_A B$ are orthogonal, we have $P_{\left( \begin{smallmatrix} A & M_A B \end{smallmatrix} \right)} = P_A + P_{M_A B}$.
    Finally, $( P_A + P_{M_A B} ) M_A = P_A M_A + M_A B (B^T M_A B)^{-1} B^T M_A M_A = 0 + M_A B (B^T M_A B)^{-1} B^T M_A = P_{M_A B}$.
\end{proofEnd}

%% file: theorems/prop_lm_derivative.tex
\begin{theoremEnd}[malte_all_end]{proposition}%
    \label{prop:lm_derivative}
    The derivative of the Lagrange multiplier test statistic
    $$
    \LM(\beta) := (n - k) \frac{(y - X \beta)^T P_{P_Z \tilde X(\beta)} (y - X \beta)}{(y - X \beta)^T M_Z (y - X \beta )}
    $$
    with respect to $\beta$ is
    \begin{multline*}        
        \frac{\dd}{\dd \beta} \LM(\beta) = -2 (n - k) \Big( \frac{(y - X \beta)^T \tilde X(\beta)} {(y - X \beta)^T M_Z (y - X \beta )} -
        \frac{(y - X \beta)^T (P_Z - P_{P_Z \tilde X(\beta)}) (y - X \beta)}{(y - X \beta)^T M_Z (y - X \beta )} \\
        \cdot \frac{(y - X \beta)^T P_Z \tilde X(\beta) (\tilde X(\beta)^T P_Z \tilde X(\beta))^{-1} \tilde X(\beta)^T M_Z \tilde X(\beta)} {(y - X \beta)^T M_Z (y - X \beta )} \Big)
    \end{multline*}
\end{theoremEnd}
\begin{proofEnd}

Let
\begin{align*}
&\Sigma(\beta) := \frac{(y - X \beta)^T M_Z X}{(y - X \beta)^T M_Z (y - X \beta)} \quad \text{ such that } \quad \tilde X(\beta) := X - (y - X \beta) \Sigma.
\end{align*}
We calculate the derivative of $\tilde X(\beta)$ with respect to $\beta_j$:
\begin{align*}
\frac{\dd}{\dd \beta_j} &\tilde X(\beta) = X_j \Sigma(\beta) + (y - X \beta) \frac{X_j^T M_Z X}{(y - X \beta)^T M_Z (y - X \beta)} 
\\
&-2 (y - X \beta) \frac{(y - X \beta)^T M_Z X_j}{(y - X \beta)^T M_Z (y - X \beta)} \Sigma(\beta) 
= \tilde X(\beta)_j \Sigma(\beta) + \frac{(y - X \beta) \cdot X_j^T M_Z \tilde X(\beta)}{(y - X \beta)^T M_Z (y - X \beta)}
\end{align*}

Let $P(a) := A(a) (A(a)^T A(a))^{-1} A(a)^T$. Write $A = A(a)$ and $A' = \frac{\dd}{\dd a} A(a)$.
Then,
\begin{align*}
\frac{\dd}{\dd a} P(a) &= A' (A^T A)^{-1} A^T + A (A^T A)^{-1} A'^T - A (A^T A)^{-1} (A'^T A + A^T A') (A^T A)^{-1} A^T \\
&= M_A A' (A^T A)^{-1} A^T +  A (A^T A)^{-1} A'^T M_A
\end{align*}
For $a = \beta_j$ and $A = P_Z \tilde X(\beta)$, we have
\begin{align*}
    M_A &A' (A^T A)^{-1} A^T \\
    &= M_{P_Z \tilde X(\beta)} P_Z \left(  \tilde X(\beta)_j \Sigma(\beta) + \frac{(y - X \beta) \cdot X_j^T M_Z \tilde X(\beta)}{(y - X \beta)^T M_Z (y - X \beta)} \right) (\tilde X(\beta)^T P_Z \tilde X(\beta))^{-1} \tilde X(\beta)^T P_Z \\
    &= \frac{M_{P_Z \tilde X(\beta)} P_Z (y - X \beta) X_j^T M_Z \tilde X(\beta) (\tilde X(\beta)^T P_Z \tilde X(\beta))^{-1} \tilde X(\beta)^T P_Z }{(y - X \beta)^T M_Z (y - X \beta)}
\end{align*}
For symmetric $Q$, one has $\frac{\dd}{\dd \beta_j} (y - X \beta)^T Q(\beta) (y - X \beta) = -2 (y - X \beta)^T Q(\beta) X_j + (y - X \beta)^T \frac{\dd}{\dd \beta_j} Q(\beta) (y - X \beta)$.
Thus,
\begin{align*}
    \frac{\dd}{\dd \beta_j} &(y - X \beta)^T P_{P_Z \tilde X(\beta)} (y - X \beta) \\
    &= -2 (y - X \beta)^T P_{P_Z \tilde X(\beta)} X_j + (y - X \beta)^T \frac{\dd}{\dd \beta_j} P_{P_Z \tilde X(\beta)} (y - X \beta) \\
    &=  -2 (y - X \beta)^T P_{P_Z \tilde X(\beta)} X_j \\
    &+ 2 \frac{(y - X \beta)^T M_{P_Z \tilde X(\beta)} P_Z (y - X \beta)}{(y - X \beta)^T M_Z (y - X \beta)} X_j^T M_Z \tilde X(\beta) (\tilde X(\beta)^T P_Z \tilde X(\beta))^{-1} \tilde X(\beta)^T P_Z (y - X \beta) \\
\end{align*}
Write $\hat\sigma^2(\beta) := (y - X \beta)^T M_Z (y - X \beta)$ and $\varepsilon_\beta := y - X \beta$.
Note that $\tilde X(\beta) = X - \varepsilon_\beta \frac{\varepsilon_\beta^T M_Z X}{\hat\sigma^2(\beta)}$ implies
$$
\varepsilon_\beta^T P_{P_Z \tilde X(\beta)} \varepsilon_\beta \cdot \varepsilon_\beta^T M_Z X_j =\varepsilon_\beta P_{P_Z \tilde X(\beta)} (X_j - \tilde X(\beta)_j) \hat\sigma^2(\beta).
$$

We calculate
\begin{align*}
    (\hat\sigma^2&(\beta))^2 \frac{\dd}{\dd \beta_j} \frac{
        (y - X \beta)^T P_{P_Z \tilde X(\beta)} (y - X \beta)
    }{
        (y - X \beta)^T M_Z (y - X \beta)
    } = -2 \varepsilon_\beta^T P_{P_Z \tilde X(\beta)} X_j \hat\sigma^2(\beta)
         \\
    &+ 2  \varepsilon_\beta^T M_{P_Z \tilde X(\beta)} P_Z \varepsilon_\beta X_j^T M_Z \tilde X(\beta) (\tilde X(\beta)^T P_Z \tilde X(\beta))^{-1} \tilde X(\beta)^T P_Z \varepsilon_\beta \\
    &+ 2 \varepsilon_\beta^T P_{P_Z \tilde X(\beta)} X_j \hat\sigma^2(\beta) - 2\varepsilon_\beta^T P_{P_Z \tilde X(\beta)} \tilde X(\beta) \hat\sigma^2(\beta) \\
    &= 2  \varepsilon_\beta^T M_{P_Z \tilde X(\beta)} P_Z \varepsilon_\beta X_j^T M_Z \tilde X(\beta) (\tilde X(\beta)^T P_Z \tilde X(\beta))^{-1} \tilde X(\beta)^T P_Z \varepsilon_\beta - 2 \varepsilon_\beta^T P_Z \tilde X(\beta)_j \hat\sigma^2(\beta) %
\end{align*}
Finally, $\tilde X^T M_Z \varepsilon_\beta = 0$ implies that $X_j^T M_Z \tilde X(\beta) = \tilde X(\beta)_j^T M_Z \tilde X(\beta)$.
\end{proofEnd}

%% file: theorems/def_wald_test_statistic.tex
\begin{definition}\label{def:wald_test_statistic}
    Let $B \in \BR^{\mX \times m}$ have ones on the diagonal and zeros elsewhere.
    Let $\kappa > 0$ and let 
    $$
    \hat \beta_\kclass(\kappa) := (S^T (\kappa P_Z + (1 - \kappa) \Id) S)^{\dagger} S^T (\kappa P_Z + (1 - \kappa) \Id) y,
    $$
    be the k-class estimator of $(\beta_0^T, \gamma_0^T)^T$ using outcomes $y$, covariates $S := \begin{pmatrix} X & W \end{pmatrix}$, and instruments $Z$, where $\dagger$ denotes the Moore-Penrose pseudoinverse.
    Let $\hat\sigma^2_{\Wald}(\kappa) := \frac{1}{n - m} \| y - S \hat\beta_\kclass(\kappa) \|^2$.
    The (subvector) Wald test statistic is
    \begin{align*}
        \Wald_{\hat \beta_\kclass(\kappa)}(\beta) := \frac{1}{\hat\sigma^2_{\Wald}(\kappa)} 
        (\beta - B \hat \beta_\kclass(\kappa))^T \left( B \left(S^T (\Id_{m} - \kappa M_Z) S\right)^{-1} B^T \right)^{-1} ( \beta - B \hat \beta_\kclass(\kappa))
    \end{align*}
\end{definition}

%% file: theorems/def_lr_test_statistic.tex
\begin{definition}
    \label{def:lr_test_statistic}
    The \emph{likelihood-ratio (LR)} test statistic is defined as
    $$
    \LR(\beta) := (k - \mW) \left( \AR(\beta) - \min_{b \in \BR^\mX} \AR(b) \right) = (k - \mW) \left( \AR(\beta) - \AR(\hat\beta_\liml) \right).
    $$
\end{definition}

%% file: theorems/prop_cis_closed_forms.tex
\begin{theoremEnd}[malte]{proposition}%
    \label{prop:cis_closed_forms}
    The inverse Wald, likelihood-ratio, and Anderson-Rubin test confidence sets describe quadrics in $\BR^{\mX}$.
    Let $S := \begin{pmatrix} X & W \end{pmatrix}$ and $B \in \BR^{\mX \times m}$ have ones on the diagonal and zeros elsewhere, such that $B S = X$.
    Let
    $$\kappa_{\AR}(\alpha) = 1 + F^{-1}_{\chi^2(k - \mW)}(1 - \alpha) / (n - k) \ \text{ and } \ \kappa_{\LR}(\alpha) = \hat\kappa_\liml + F^{-1}_{\chi^2(\mX)}(1 - \alpha) / (n - k).$$
    Let
    \begin{align*}
    A(\kappa) &:= \left( B \left(S^T (\Id_n - \kappa M_Z) S \right)^{-1} B^T \right)^{-1} =
    X^T (\Id_n - \kappa M_Z) X \\ & \hspace{1.5cm} - X^T (\Id_n - \kappa M_Z) W \left(W^T (\Id_n - \kappa M_Z) W \right)^{-1} W^T (\Id_n - \kappa M_Z) X.
    \end{align*}
    Let $\hat\sigma^2_\mathrm{Wald}(\kappa) := \frac{1}{n - m} \| y - S \hat\beta_\kclass(\kappa) \|^2$ and $\hat\sigma^2(\kappa) := \frac{1}{n-k} \| M_Z (y - S \hat\beta_\kclass(\kappa)) \|^2$.
    If $\mW > 0$, let $\kappa_\mathrm{max} := \lambdamin{(W^T M_Z W)^{-1} W^T P_Z W} + 1$. Else $\kappa_{\min} := \infty$.
    \begin{itemize}
    \item If $\kappa \leq \kappa_\mathrm{max}$, then
    \begin{multline*}
        \CI_{\mathrm{Wald}_{\hat\beta_\kclass(\kappa)}}(1 - \alpha) = \Big\{ \beta \in \BR^{\mX}  \mid \left(\beta - B \hat\beta_\kclass(\kappa) 
        \right)^T A(\kappa) \left(\beta - B\hat\beta_\kclass(\kappa) \right) \\
        \leq \hat\sigma^2_\mathrm{Wald}(\kappa) \cdot F^{-1}_{\chi^2(\mX)}(1 - \alpha)\Big\}.
    \end{multline*}
    Else, $\CI_{\mathrm{Wald}_{\hat\beta_\kclass(\kappa)}}(1 - \alpha) = \BR^{\mX}$.
    \item  If $\kappa_{\AR}(\alpha) \leq \kappa_\mathrm{max}$, then
    \begin{multline*}
        \CI_{\AR}(1 - \alpha) = \Big\{ \beta \in \BR^{\mX} \mid \left(\beta - B \hat\beta_\kclass(\kappa(\alpha))\right)^T A(\kappa_{\AR}(\alpha)) \left(\beta - B\hat\beta_\kclass(\kappa_{\AR}(\alpha)) \right) \\
        \leq \hat\sigma^2(\kappa_{\AR}(\alpha)) \cdot (F^{-1}_{\chi^2(k - \mW)}(1 - \alpha) - k \AR(\hat\beta_\kclass(\kappa_{\AR}(\alpha))))\Big\}.
    \end{multline*}
    Else, $\CI_{\AR}(1 - \alpha) = \BR^{\mX}$.
    \item If $\kappa_{\LR}(\alpha) \leq \kappa_\mathrm{max}$, then
    \begin{multline*}
        \CI_{\LR}(1 - \alpha) = \Big\{ \beta \in \BR^{\mX} \mid \left(\beta - B \hat\beta_\kclass(\kappa_{\LR}(\alpha)) \right)^T A(\kappa_{\LR}(\alpha))  \left(\beta - B \hat\beta_\kclass(\kappa_{\LR}(\alpha)) \right) \\
        \leq \hat\sigma^2(\kappa_{\LR}(\alpha)) \cdot (F^{-1}_{\chi^2(\mX)}(1 - \alpha) + k \AR(\hat\beta_\liml) - k \AR(\hat\beta_\kclass(\kappa_{\LR}(\alpha)))) \Big\}.
    \end{multline*}
    Else, $\CI_{\LR}(1 - \alpha) = \BR^{\mX}$.
\end{itemize}
\end{theoremEnd}
\begin{proofEnd}
    Let $\mW = 0$.
    By definition
    \begin{multline*}
        \mathrm{Wald}_{\hat \beta_\kclass(\kappa)}(\beta) \leq F^{-1}_{\chi^2(\mX)}(1 - \alpha) \Leftrightarrow (\beta - \hat\beta_\kclass(\kappa))^T \left( X^T (\kappa P_Z + (1 - \kappa) \Id_n) X \right) (\beta - \hat\beta_\kclass(\kappa)) \\   
        \leq \hat\sigma^2_\mathrm{Wald}(\kappa) \cdot F^{-1}_{\chi^2(\mX)}(1 - \alpha). \numberthis \label{eq:ci_closed_form:wald_1}    
    \end{multline*}
    Let $\tilde \AR(\beta) := (n - k) \frac{(y - X \beta)^T P_Z (y - X \beta)}{(y - X \beta)^T M_Z (y - X \beta)}$.
    We prove that
    \begin{align*}
        \tilde\AR(\beta) \leq (n - k)(\kappa - 1) \Leftrightarrow (\beta - \hat\beta_\kclass(\kappa))^T \left( X^T (\kappa P_Z + (1 - \kappa) \Id) X \right) (\beta -  \hat\beta_\kclass(\kappa)) \\
        \leq \hat\sigma^2(\kappa) \cdot \left((n - k) (\kappa - 1)- \tilde\AR(\hat\beta_\kclass(\kappa)) \right). \numberthis \label{eq:ci_closed_form:1}
    \end{align*}
    Calculate
    \begin{align*}
        \tilde\AR&(\beta) = (n - k) \frac{(y - X \beta)^T P_Z (y - X \beta)}{(y - X \beta)^T M_Z (y - X \beta)} \leq (n - k) (\kappa - 1) \\
        &\Leftrightarrow (y - X \beta)^T P_Z (y - X \beta) \leq (\kappa - 1) (y - X \beta)^T M_Z (y - X \beta)\\
        &\Leftrightarrow (y - X \beta)^T (P_Z + (1 - \kappa) M_Z) (y - X \beta) \leq 0.
    \end{align*}
    Here,
    \begin{align*}
        &(y - X \beta)^T (P_Z + (1 - \kappa) M_Z) (y - X \beta)  \\
        &= \beta^T \left( X^T ( P_Z + (1 - \kappa) M_Z) X \right) \beta - 2 \beta^T \underbrace{X^T (P_Z + (1 - \kappa) M_Z) y}_{\mathclap{=(X^T (P_Z + (1 - \kappa) M_Z) X) \hat\beta_\kclass(\kappa)}} + y^T (P_Z + (1 - \kappa) M_Z) y \\
        &= (\beta -  \hat\beta_\kclass(\kappa))^T \left( X^T (P_Z + (1 - \kappa) M_Z) X \right) (\beta -  \hat\beta_\kclass(\kappa))
        + y^T (P_Z + (1 - \kappa) M_Z) y \\
        & \hspace{3cm} - \underbrace{\hat\beta_\kclass(\kappa)^T X^T (P_Z + (1 - \kappa) M_Z) X \hat\beta_\kclass(\kappa)}_{\mathclap{= \hat\beta_\kclass(\kappa)^T X^T (P_Z + (1 - \kappa) M_Z) y}}\\
        &= (\beta -  \hat\beta_\kclass(\kappa))^T \left( X^T (P_Z + (1 - \kappa) M_Z) X \right) (\beta - \hat\beta_\kclass(\kappa)) \\
        &\hspace{3cm} +  (y -  X \hat\beta_\kclass(\kappa))^T  (P_Z + (1 - \kappa) M_Z) (y - X \hat\beta_\kclass(\kappa)),
    \end{align*}
    as $(y - X \hat\beta_\kclass(\kappa))^T (P_Z + (1 - \kappa) M_Z) X \hat\beta_\kclass(\kappa) = 0$. Now, rewrite
    \begin{align*}
        (y& - X \hat\beta_\kclass(\kappa))^T (P_Z + (1 - \kappa) M_Z)  (y - X \hat\beta_\kclass(\kappa)) \\
        &= (y - X \hat\beta_\kclass(\kappa)) P_Z (y - X \hat\beta_\kclass(\kappa)) + (1 - \kappa) (y - X \hat\beta_\kclass(\kappa))^T M_Z (y - X \hat\beta_\kclass(\kappa)) \\
        &= \hat\sigma^2(\kappa) \cdot (\tilde\AR(\hat\beta_\kclass(\kappa)) + (n - k) (1 - \kappa)).
    \end{align*}
    This proves \eqref{eq:ci_closed_form:1}.
    If $\mW = 0$, then $\beta \in \CI_{\AR}(1 - \alpha)$ if and only if $\tilde\AR(\beta) \leq F^{-1}_{\chi^2(k)}(1 - \alpha)$.
    The expression of $\CI_{\AR}(1 - \alpha)$ follows from \eqref{eq:ci_closed_form:1} with 
    \begin{align*}(n - k) (\kappa_{\AR}(\alpha) - 1) = F^{-1}_{\chi^2(k)}(1 - \alpha) \Leftrightarrow \kappa_{\AR}(\alpha) = 1 + F^{-1}_{\chi^2(k)}(1 - \alpha) / (n-k).
    \end{align*}
    Also, $\beta \in \CI_{\LR}(1 - \alpha)$ if and only if $\tilde\AR(\beta) - \min_b \tilde\AR(b) \leq F^{-1}_{\chi^2(\mX)}(1 - \alpha) / (n - k)$.
    The expression of $\CI_{\LR}(1 - \alpha)$ thus follows from \eqref{eq:ci_closed_form:1} with $\kappa_{\LR}(\alpha) = 1 + (\min_b \tilde\AR(b) + F^{-1}_{\chi^2(\mX)}(1 - \alpha)) / (n-k) = \hat\kappa_\liml + F^{-1}_{\chi^2(\mX)}(1 - \alpha) / (n-k)$.

    Apply \Cref{lem:kappa_pos_definite} with $X \leftarrow W$. This implies that $W^T (\kappa P_Z + (1 - \kappa) \Id) W$ is positive semi-definite for $\kappa \leq \kappa_\mathrm{max}$ and negative definite if $\kappa > \kappa_\mathrm{max}$.

    Calculate
    \begin{multline*}
        \{ \beta \in \BR^\mX \mid \min_\gamma \tilde{\AR} \left( \begin{pmatrix} \beta \\ \gamma \end{pmatrix} \right) \leq (n - k) \cdot (\kappa - 1) \}  \\
        \overset{\eqref{eq:ci_closed_form:1}}{=} \Big\{\beta \in \BR^\mX \mid \min_\gamma \left(\begin{pmatrix} \beta \\ \gamma \end{pmatrix} - \hat\beta_\kclass(\kappa) \right)^T (S^T (\kappa P_Z + (1 - \kappa) \Id) S) \left(\begin{pmatrix} \beta \\ \gamma \end{pmatrix} - \hat\beta_\kclass(\kappa) \right) \\
        \leq \hat\sigma^2(\kappa) \cdot \left((n - k) (\kappa - 1) - \tilde\AR(\hat\beta_\kclass(\kappa)) \right) \Big\} \\
        \overset{\Cref{lem:projected_quadric}}{=}
        \begin{cases}
            \ \ \!\begin{aligned}
                \Big\{ \beta \in \BR^{\mX} \mid (\beta - &B \hat\beta_\kclass(\kappa))^T ( B (S^T (P_Z + (1 - \kappa) M_Z) S)^\dagger B^T)^{\dagger} \\
                &(\beta - B \hat\beta_\kclass(\kappa)) 
                 \leq \hat\sigma^2(\kappa) \cdot ((n - k) (\kappa - 1) - \tilde\AR(\hat\beta_\kclass(\kappa))) \Big\}
            \end{aligned} & \ \ \text{if } \kappa \leq \kappa_\mathrm{max} \\
            \ \ \BR^{\mX} & \ \ \text{if } \kappa > \kappa_\mathrm{max}
        \end{cases},
    \end{multline*}
    with $A = S^T (\kappa P_Z + (1 - \kappa) \Id) S$ and $A_{22} = W^T ( \kappa P_Z + (1 - \kappa) \Id) W$, where $\dagger$ denotes the Moore-Penrose pseudo-inverse.
    Again, the results for $\CI_{\AR}(1 - \alpha)$ and $\CI_{\LR}(1 - \alpha)$ follow from the cases via $\kappa = \kappa_{\AR}(\alpha)$ and $\kappa = \kappa_{\LR}(\alpha)$.
\end{proofEnd}

%% file: theorems/prop_cis_bounded.tex
\begin{theoremEnd}[malte]{proposition}%
    \label{prop:cis_bounded}
    Assume the \Cref{tc:2} holds.
    Let
    \begin{align*}
        J_\liml &= (k - \mW) \cdot \min_b \AR(b) \\
        &= (n - k) \ \lambdamin{ \left(\begin{pmatrix} X & W & y \end{pmatrix}^T M_Z \begin{pmatrix} X & W & y \end{pmatrix} \right)^{-1}\begin{pmatrix} X & W & y \end{pmatrix}^T P_Z \begin{pmatrix} X & W & y \end{pmatrix}}
    \end{align*}
    be the LIML variant of the J-statistic \citep{guggenberger2012asymptotic} and let
    $$
        \lambda = (n - k) \ \lambdamin{ \left(\begin{pmatrix} X & W \end{pmatrix}^T M_Z \begin{pmatrix} X & W \end{pmatrix} \right)^{-1}\begin{pmatrix} X & W \end{pmatrix}^T P_Z \begin{pmatrix} X & W \end{pmatrix}}.
    $$
    be \citeauthor{anderson1951estimating}'s \citeyearpar{anderson1951estimating} likelihood-ratio test statistic of reduced rank.

    \begin{itemize}
        \item The (subvector) inverse Anderson-Rubin test confidence set is nonempty and bounded if and only if
    \begin{align*}
        J_\liml \leq F^{-1}_{\chi^2(k - \mW)}(1 - \alpha) < \lambda \ \Leftrightarrow \
        F_{\chi^2(k - \mW)} \left( J_\liml \right) \leq 1 - \alpha < F_{\chi^2(k - \mW)} \left(\lambda \right)
    \end{align*}
        \item The (subvector) inverse likelihood-ratio test confidence set is always nonempty. It is bounded if and only if
    \begin{align*}
        F^{-1}_{\chi^2(\mX)}(1 - \alpha) < \lambda - J_\liml \ \Leftrightarrow \  1 - \alpha \in \left[0, F_{\chi^2(\mX)} \left( \lambda - J_\liml \right) \right)
    \end{align*}

    \end{itemize}
\end{theoremEnd}
\begin{proofEnd}
    We prove the conditions for boundedness and non-emptiness separately.
    \paragraph*{Step 1: (non)emptiness of the confidence sets}
    The (subvector) inverse likelihood-ratio test confidence set always contains the LIML and hence is never empty.
    The (subvector) inverse Anderson-Rubin test confidence set is nonempty if and only if
    $$
    J_\liml = (k - \mW) \cdot \AR(\hat\beta_\liml) \leq F^{-1}_{\chi^2(k - \mW)}(1 - \alpha) \Leftrightarrow \alpha \leq 1 - F_{\chi^2(k - \mW)} \left( J_\liml  \right).
    $$

    \paragraph*{Step 2: Boundedness of the confidence sets}
    We show the boundedness equivalence simultaneously for the subvector inverse Anderson-Rubin and likelihood-ratio test confidence sets
    Setting $\kappa = \kappa_{\AR}(\alpha) = 1 + F^{-1}_{\chi^2(k - \mW)}(1 - \alpha)$ and $\kappa = \kappa_{\LR}(\alpha) = \hat\kappa_\liml + F^{-1}_{\chi^2(\mX)}(1 - \alpha) = 1 + (n - k )(J_\liml + F^{-1}_{\chi^2(\mX)}(1 - \alpha))$, we need to show that the confidence set is bounded if and only if $(n - k)(\kappa - 1) < \lambda$.
    By \Cref{prop:cis_closed_forms}, the confidence sets are bounded if and only if $\kappa \leq \kappa_\mathrm{max}$ and $A(\kappa)$ is positive definite.

    \paragraph*{Step 2a: $(n - k)(\kappa - 1) < \lambda \Leftrightarrow S^T (P_Z + (1 - \kappa) M_Z) S$ is positive definite}
    This follows from \Cref{lem:kappa_pos_definite} with $X \leftarrow S = \begin{pmatrix} X & W \end{pmatrix}$.

    \paragraph*{Step 2b: $(n - k)(\kappa - 1) < \lambda$ implies the confidence set is bounded.}

    By step 2a, $(n - k)(\kappa - 1) < \lambda \Rightarrow S^T (P_Z + (1 - \kappa) M_Z) S$ is positive definite.
    By \citep[Proposition 2.1]{gallier2010schur}, the matrix $ S^T (P_Z + (1 - \kappa) M_Z) S$ is positive definite if and only if the submatrix $W^T (P_Z + (1 - \kappa) M_Z) W$ and the Schur complement $A(\kappa)$ are positive definite.
    Applying \Cref{lem:kappa_pos_definite} with $X \leftarrow W$ yields $\kappa < \kappa_\mathrm{max}$.

    \paragraph*{Step 2c: If $\kappa = \kappa_\mathrm{max}$ and \Cref{tc:2} \ref{tc:2:a} holds, the confidence set is unbounded.}
    If $\kappa = \kappa_\mathrm{max}$, then $\lambdamin{S^T (P_Z + (1 - \kappa) M_Z) S} \leq \lambdamin{W^T (P_Z + (1 - \kappa) M_Z) W} = 0.$
    By \Cref{tc:2} (\ref{tc:2:a}), this inequality is strict and thus $\lambdamin{S^T (P_Z + (1 - \kappa) M_Z) S} < 0$.
    Let $v = \begin{pmatrix} v_X \\ v_W \end{pmatrix}$ be an eigenvector corresponding to $\lambdamin{S^T (P_Z + (1 - \kappa) M_Z) S}$.
    Then $v_X \neq 0$, as else $0 > \lambdamin{S^T (P_Z + (1 - \kappa) M_Z) S} = v^T S^T (P_Z + (1 - \kappa) M_Z) S v / \|v\|^2 = v_W^T W^T (P_Z + (1 - \kappa) M_Z) W v_W / \|v_W \|^2 \leq \lambdamin{W^T (P_Z + (1 - \kappa) M_Z) W} = 0$.
    Write
    \begin{align*}
        \CI_{\AR}(1 - \alpha) &= \{ \beta\in\BR^\mX \mid \min_\gamma \begin{pmatrix} \beta \\ \gamma \end{pmatrix}^T S^T (P_Z + (1 - \kappa) M_Z) S \begin{pmatrix} \beta \\ \gamma \end{pmatrix} \leq \mathrm{const} \} + B \hat\beta_\kclass(\kappa) \\
        &\supseteq \{ t \cdot v_X \mid t^2 v^T S^T (P_Z + (1 - \kappa) M_Z) S v \leq \mathrm{const} \} + B \hat\beta_\kclass(\kappa).
    \end{align*}
    See also the proof of \Cref{prop:cis_closed_forms}.
    The latter is equal to $\BR v_X$ if $\mathrm{const} \geq 0$ and $\left( (-\infty, t_\mathrm{min}] \cup [t_\mathrm{min}, \infty) \right) \cdot v_X + B \hat\beta_\kclass(\kappa)$ for $t_\mathrm{min} = \frac{1}{\|v\|} \sqrt{\frac{\mathrm{const}}{\lambdamin{S^T (P_Z + (1 - \kappa) M_Z) S}}}$ otherwise.
    In both cases, the confidence set is unbounded.

    \paragraph*{Step 2d: A bounded confidence set implies $(n - k)(\kappa - 1) < \lambda$}
    If the confidence set is bounded, then $\kappa \leq \kappa_\mathrm{max}$ and $A(\kappa)$ is positive definite.
    By step 2c, the $\kappa < \kappa_\mathrm{max}$.
    Then by \Cref{lem:kappa_pos_definite}, the submatrix $W^T (P_Z + (1 - \kappa) M_Z) W$ is positive definite and by the if direction of \citep[Proposition 2.1]{gallier2010schur},  $S^T (P_Z + (1 - \kappa) M_Z) S$ is positive definite.
    By step 2a, $(n - k)(\kappa - 1) < \lambda$.

\end{proofEnd}

%% file: theorems/cor_ar_bounded_iff_rank_test_rejects.tex
\begin{theoremEnd}[malte]{corollary}%
    \label{cor:ar_bounded_iff_rank_test_rejects}
    Consider structural equations $y_i = S_i^T \beta_0 + \varepsilon_i$ and $S_i = Z_i^T \Pi + V_i$.
    We are interested in inference for each component of the causal parameter $\beta$.
    Thus, for each covariate $i = 1, \ldots, m$, we separate $X \leftarrow S^{(i)}$ and $W \leftarrow S^{(-i)}$ and construct confidence intervals $\CI_{\AR}^{(i)}(1 - \alpha)$ and $\CI_{\LR}^{(i)}(1 - \alpha)$ for the $i$-th component of $\beta_0$.
    Then
    \begin{enumerate}[(a)]
        \item The subvector inverse Anderson-Rubin test confidence sets are jointly (un)bounded. That is, if $\CI_{\AR}^{(i)}(1 - \alpha)$ is (un)bounded for any $i$, then $\CI_{\LR}^{(i)}(1 - \alpha)$ is (un)bounded for all $i$.
        \item The subvector inverse Anderson-Rubin test confidence sets at level $\alpha$ are bounded (and thus confidence intervals) if and only if \citeauthor{anderson1951estimating}'s \citeyearpar{anderson1951estimating} likelihood-ratio test rejects the null hypothesis that $\Pi$ is of reduced rank at level $\alpha$. \label{cor:ar_bounded_iff_rank_test_rejects:b}
        \item The subvector inverse Anderson-Rubin test confidence sets are jointly (non)empty. That is, if $\CI_{\AR}^{(i)}(1 - \alpha)$ is (non)empty for any $i$, then $\CI_{\AR}^{(i)}(1 - \alpha)$ is (non)empty for all $i$.
        \item The subvector inverse Anderson-Rubin test confidence sets at level $\alpha$ are jointly empty if and only if the LIML variant of the J-statistic $J_\liml > F^{-1}_{\chi^2(k - m + 1)}(1 - \alpha)$.
        \item The subvector inverse likelihood-ratio test confidence sets are jointly (un)bounded. That is, if $\CI_{\LR}^{(i)}(1 - \alpha)$ is (un)bounded for any $i$, then $\CI_{\LR}^{(i)}(1 - \alpha)$ is (un)bounded for all $i$.
    \end{enumerate}
\end{theoremEnd}
\begin{proofEnd}
    For each $i=1,\ldots,m$, the conditions for the subvector inverse Anderson-Rubin or likelihood-ratio test confidence sets to be bounded and nonempty follow from \Cref{prop:cis_bounded}.
    That the confidence sets are jointly (un)bounded or (non)empty follows from the fact that these conditions do not depend on how $S = \begin{pmatrix} S^{(i)} & S^{(-i)} \end{pmatrix}$ is partitioned.
    \Cref{cor:ar_bounded_iff_rank_test_rejects:b} follows from \citet[][Corollary 39]{londschien2025overview} with $r=1$.
\end{proofEnd}

%% file: tables/table_guggenberger12_empirical_sizes_n50.tex
\begin{table}[h]
    \center
    \begin{tabular}{l | rrrrr | rrrrr}
         & \multicolumn{5}{r}{$\alpha = 0.01$} & \multicolumn{5}{|r}{$\alpha = 0.05$} \\
        \multicolumn{1}{ r |}{$k=$ } & 5 & 10 & 15 & 20 & 30 & 5 & 10 & 15 & 20 & 30 \\
        \hline
        AR & 0.5\% & 0.7\% & 0.7\% & 0.9\% & 1.9\% & 3.0\% & 3.2\% & 3.3\% & 3.4\% & 5.0\% \\
        AR (GKM) & 1.2\% & 1.5\% & 1.6\% & 1.7\% & 2.6\% & 5.5\% & 6.4\% & 5.9\% & 6.4\% & 7.0\% \\
        CLR & 0.5\% & 0.8\% & 0.8\% & 1.2\% & 2.0\% & 3.1\% & 3.8\% & 4.0\% & 4.9\% & 6.7\% \\
        LM (ours) & 0.3\% & 0.7\% & 0.7\% & 1.0\% & 1.8\% & 2.4\% & 3.0\% & 3.8\% & 5.0\% & 7.5\% \\
        LM (LIML) & 1.4\% & 4.1\% & 5.8\% & 8.5\% & 11\% & 6.6\% & 12\% & 16\% & 21\% & 25\% \\
    \end{tabular}
        
    \caption{
        \label{tab:guggenberger12_empirical_sizes_n50}
        \tableguggenbergerempiricalsizes{$n=50$}
    }
\end{table}

%% file: tables/table_guggenberger12_empirical_sizes_n100.tex
\begin{table}[h]
    \center
    \begin{tabular}{l | rrrrr | rrrrr}
         & \multicolumn{5}{r}{$\alpha = 0.01$} & \multicolumn{5}{|r}{$\alpha = 0.05$} \\
        \multicolumn{1}{ r |}{$k=$ } & 5 & 10 & 15 & 20 & 30 & 5 & 10 & 15 & 20 & 30 \\
        \hline
        AR & 0.4\% & 0.3\% & 0.2\% & 0.3\% & 0.4\% & 2.7\% & 2.4\% & 2.1\% & 2.1\% & 2.2\% \\
        AR (GKM) & 1.2\% & 1.1\% & 1.0\% & 1.2\% & 1.1\% & 5.5\% & 5.6\% & 4.9\% & 5.3\% & 5.1\% \\
        CLR & 0.4\% & 0.3\% & 0.2\% & 0.2\% & 0.5\% & 2.9\% & 2.5\% & 2.3\% & 2.6\% & 2.4\% \\
        LM (ours) & 0.3\% & 0.2\% & 0.4\% & 0.4\% & 0.5\% & 2.1\% & 2.2\% & 2.5\% & 2.8\% & 3.2\% \\
        LM (LIML) & 1.4\% & 3.0\% & 4.5\% & 6.7\% & 9.8\% & 6.2\% & 11\% & 14\% & 18\% & 23\% \\
    \end{tabular}
        
    \caption{
        \label{tab:guggenberger12_empirical_sizes_n100}
        \tableguggenbergerempiricalsizes{$n=100$}
    }
\end{table}

%% file: tables/table_card_black.tex
\begin{table}[h]
    \begin{center}
    \begin{tabular}{lrrrrrr}
        specification & \multicolumn{2}{c}{\ (i) $k=4$} & \multicolumn{2}{c}{\ \ (ii) $k=8$} & \multicolumn{2}{c}{\ \ (iii) $k=14$} \\
        \midrule
        estimate (TSLS) &  \multicolumn{2}{r}{0.013 (0.041)} & \multicolumn{2}{r}{0.000 (0.041)} & \multicolumn{2}{r}{-0.010 (0.021)}  \\
        estimate (LIML) &  \multicolumn{2}{r}{0.013 (0.041)} & \multicolumn{2}{r}{0.001 (0.043)} & \multicolumn{2}{r}{-0.012 (0.021)}  \\
        
        \rule{0pt}{4ex}test & \multicolumn{1}{l}{\ \ \ \ statistic} & \multicolumn{1}{l}{\ $p$-value} & \multicolumn{1}{l}{\ \ \ \ statistic} & \multicolumn{1}{l}{\ $p$-value} & \multicolumn{1}{l}{\ \ \ \ statistic} & \multicolumn{1}{l}{\ $p$-value}  \\
        \midrule
        Wald (TSLS) & 0.10 & 0.756 & 0.00 & 0.998 & 0.24 & 0.624 \\
         \multicolumn{1}{r}{\ \ 95\% conf. set} & \multicolumn{2}{r}{[-0.07, 0.09]} & \multicolumn{2}{r}{[-0.08, 0.08]  } & \multicolumn{2}{r}{[-0.05, 0.03]} \\
        \rule{0pt}{3ex}Wald (LIML) & 0.10 & 0.756 & 0.00 & 0.984 & 0.33 & 0.567 \\
          \multicolumn{1}{r}{\ \ 95\% conf. set} & \multicolumn{2}{r}{[-0.07, 0.09]} & \multicolumn{2}{r}{[-0.08, 0.09]} & \multicolumn{2}{r}{[-0.05, 0.03]} \\
        \rule{0pt}{3ex}AR & 0.10 & 0.757 & 1.09 & 0.364 & 0.69 & 0.754 \\
          \multicolumn{1}{r}{\ \ 95\% conf. set} & \multicolumn{2}{r}{ [-0.07, 0.10] } & \multicolumn{2}{r}{ [-0.11, 0.11] } & \multicolumn{2}{r}{[-0.10, 0.06] }  \\
        \rule{0pt}{3ex}CLR  & 0.10 & 0.757 & 0.00 & 0.985 & 0.33 & 0.607 \\
        \multicolumn{1}{r}{\ \ 95\% conf. set} & \multicolumn{2}{r}{ [-0.07, 0.10] } & \multicolumn{2}{r}{ [-0.09, 0.10] } & \multicolumn{2}{r}{[-0.06, 0.03]} \\
        \rule{0pt}{3ex}LM (ours)  & 0.10 & 0.757 & 0.00 & 0.984 & 0.32 & 0.57 \\
        \multicolumn{1}{r}{\ \ 95\% conf. set} & \multicolumn{2}{r}{ [-0.07, 0.10] } & \multicolumn{2}{r}{[-0.09, 0.10]} & \multicolumn{2}{r}{$ [-0.06, 0.03]\ \cup \ $} \\
         & \multicolumn{2}{r}{\ } & \multicolumn{2}{r}{\ } & \multicolumn{2}{r}{$ [0.16, 0.28]$}  \\
        \midrule
        rank  & 12.14 & 4.93e-4 & 28.81 & 2.53e-05 & 45.91 & 3.35e-06 \\
        $J_\liml$ & & & 5.45 & 0.244 & 7.21 & 0.706 \\
      \end{tabular}
      \caption{
        \label{tab:applications:card_3}
        Estimates, test statistics, $p$-values, and confidence sets for the causal effect of education on wages for blacks. See \Cref{sec:applications:card} for details.
      }
    \end{center}

\end{table}

%% file: supplement.tex



\input{sec_optimization.tex}
\input{sec_numerical_details.tex}
\FloatBarrier
\newpage
\section{Additional figures}
\FloatBarrier

\label{sec:additional_figures}
\input{figures/figure_tc3_counterexample.tex}
\input{figures/figure_optimization.tex}
\input{figures/figure_guggenberger12_qqplots.tex}
\input{figures/figure_guggenberger12_power_rho=0.8.tex}
\input{figures/figure_guggenberger12_power_rho=0.99.tex}
\input{figures/figure_kleibergen19_identity_k5_20.tex}
\input{figures/figure_kleibergen19_guggenberger12_k5_20.tex}
\input{figures/figure_kleibergen19_guggenberger12_k100.tex}
\input{figures/figure_card.tex}
\input{figures/figure_tanaka.tex}

%% file: sec_optimization.tex
\section{Optimization}
\label{sec:optimization}
The optimization problem in \Cref{def:subvector_klm_test_statistic} is non-convex and thus difficult to solve.
In our \texttt{ivmodels} software package, which was used for the simulations in \Cref{sec:numerical_analysis}, we use the quasi-Newton method of Broyden–Fletcher–Goldfarb–Shannon (BFGS).
This is the default option of \texttt{scipy.optimize.minimize} in Python.

The BFGS algorithm requires an initial guess $x_0$.
We start the optimization once at the LIML estimate, using the outcomes $y - X \beta_0$, endogenous covariates $W$, and instruments $Z$, and once at $x_0 = 0$.
We then take the minimum of the two solutions as the final estimate.
Empirically, we observe that even if we only use the solution of the optimization started at the LIML estimate, the resulting algorithm has the correct size for the simulation settings from \Cref{sec:numerical_analysis}.

In \Cref{fig:optimization} in Appendix \ref{sec:additional_figures}, we plot the Lagrange multiplier statistic for different values of $\gamma$ at $\beta_0$ for draws from \citeauthor{guggenberger2012asymptotic}'s \citeyearpar{guggenberger2012asymptotic} data-generating process.

%% file: sec_numerical_details.tex
\section{Simulation details}
\label{sec:simulation_details}
To sample $\Pi_X, \Pi_W$ such that $ \sqrt{n} \| \Pi_X \| = 100, \sqrt{n} \| \Pi_W \| = 1$ and $n \langle \Pi_X, \Pi_W \rangle = 95$, 
we first sample the entries of $\Pi_1, \Pi_2 \overset{\iid}{\sim} \mathcal{N}(0, 1)$, subtract their respective means, and then assign $\Pi_X = \Pi_1 / \| \Pi_1 \|$ and $\Pi_W = 0.95 * \Pi_X + \sqrt{1 - 0.95^2} * M_{\Pi_1} \Pi_2 / \| M_{\Pi_1} \Pi_2 \|$.
We finally scale $\Pi_X \leftarrow 100 \ \Pi_X / \sqrt{n}$ and $\Pi_W \leftarrow \Pi_W / \sqrt{n}$.

For the simulation in \Cref{sec:numerical_analysis:kleibergen19}, we sample $Z_i \overset{\text{i.i.d.}}{\sim} \CN(0, \Id_k)$ such that $Q = \Id_k$.
For some $(\lambda_1, \lambda_2, \tau)$, we calculate 
$$
\Omega_{VV\cdot \varepsilon}^{1/2} R \Lambda^T \Lambda R^T \Omega_{VV\cdot \varepsilon}^{1/2} \overset{!}{=} n \Pi^T Q \Pi = n \begin{pmatrix} \| \Pi_X \|^2 & \langle \Pi_X, \Pi_W \rangle \\ \langle \Pi_X, \Pi_W \rangle  & \| \Pi_W \|^2 \end{pmatrix}.
$$
We then sample $\Pi$ according to the above paragraph.
We use a grid of 50 values for $\lambda_1, \lambda_2 \in [0, 100]$ and 25 values for $\tau \in [0, \pi)$.
Note that this is equivalent to 50 values for $\tau \in [0, 2\pi)$, as $R(\tau) = - R(\tau + \pi)$.

%% file: figures/figure_tc3_counterexample.tex
\begin{figure}[H]
    \centering
    \includegraphics[width=1\textwidth]{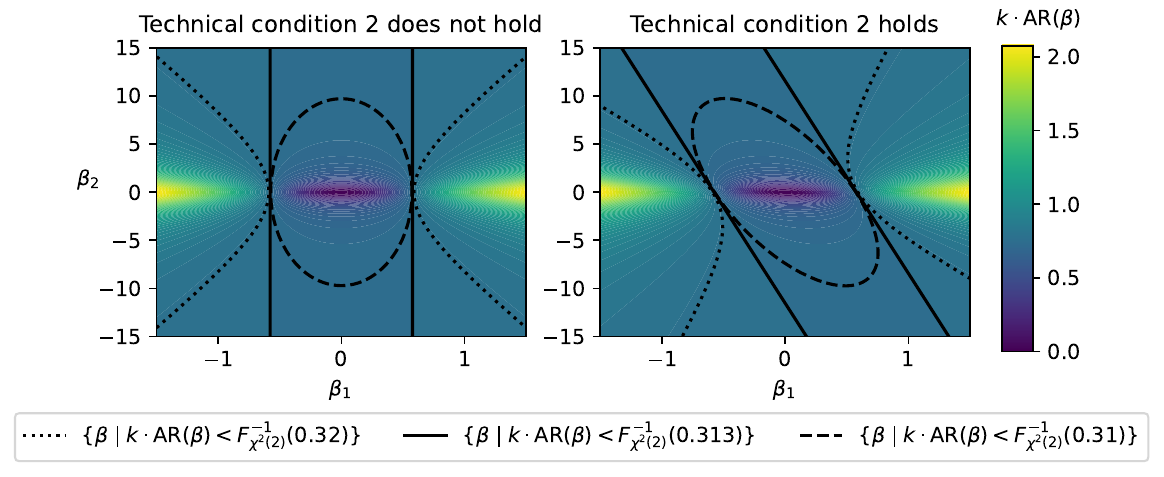}
    \caption{
        \label{fig:figure_tc3_counterexample}
        Let \\
        \parbox{\linewidth}{\centering
        \begin{equation*}
        X = \begin{pmatrix} 0.5 & 0 \\ 0 & 1 \\ 0 & 0 \\ 1 & 0 \\ 0 & 1 \\ 0 & 0 \end{pmatrix}, \ y = \begin{pmatrix} 0 \\ 0 \\ 0 \\0 \\ 0 \\ 1\end{pmatrix}, \text{ and } Z = \begin{pmatrix} 1 & 0 & 0 \\ 0 & 1 & 0 \\ 0 & 0 & 1\\ 0 & 0 & 0 \\ 0 & 0 & 0 \\ 0 & 0 & 0 \end{pmatrix} \Rightarrow \ 
        \begin{aligned}[c]
        &\begin{pmatrix} X & y \end{pmatrix}^T M_Z \begin{pmatrix} X & y \end{pmatrix} = \Id_3 \\
        &\begin{pmatrix} X & y \end{pmatrix}^T P_Z \begin{pmatrix} X & y \end{pmatrix} = \begin{pmatrix} 0.25 & 0 & 0 \\ 0 & 1 & 0 \\ 0 & 0 & 0 \end{pmatrix}.
        \end{aligned}
        \end{equation*}}
        Here, \Cref{tc:2} does not hold.
        The confidence sets $\CI_{\AR}(1 - \alpha)$ for parameters $\beta_1$ and $\beta_2$ individually result from projecting the quadrics shown in the left panel onto the respective axes.
        If $1 - \alpha = 0.32$, the confidence sets for both coefficients are unbounded.
        If $1 - \alpha = 0.31$, the confidence sets for both coefficients are bounded.
        If $1 - \alpha = 1 - F_{\chi^2(2)}^{-1}(0.25) \approx 0.313$, the confidence set for $\beta_1$ is bounded, while the confidence set for $\beta_2$ is unbounded.
        This occurs as the principal axis of the quadric is aligned with the $\beta_2$ axis and \Cref{tc:2} does not apply.
        In the right panel, we rotated $X \leftarrow X \begin{pmatrix} 1 & 0.05 \\ 0.05 & 1 \end{pmatrix}$.
        \Cref{tc:2} holds, and for $1 - \alpha = 1 - F_{\chi^2(2)}^{-1}(0.25) \approx 0.683$, the confidence sets for both coefficients are unbounded.
        If the data was generated according to \Cref{model:1} with noise distribution absolutely continuous with respect to the Lebesgue measure, \Cref{tc:2} applies with probability one.
    }
\end{figure}

%% file: figures/figure_optimization.tex
\begin{figure}[htpb]
    \centering
    \includegraphics[width=0.95\textwidth]{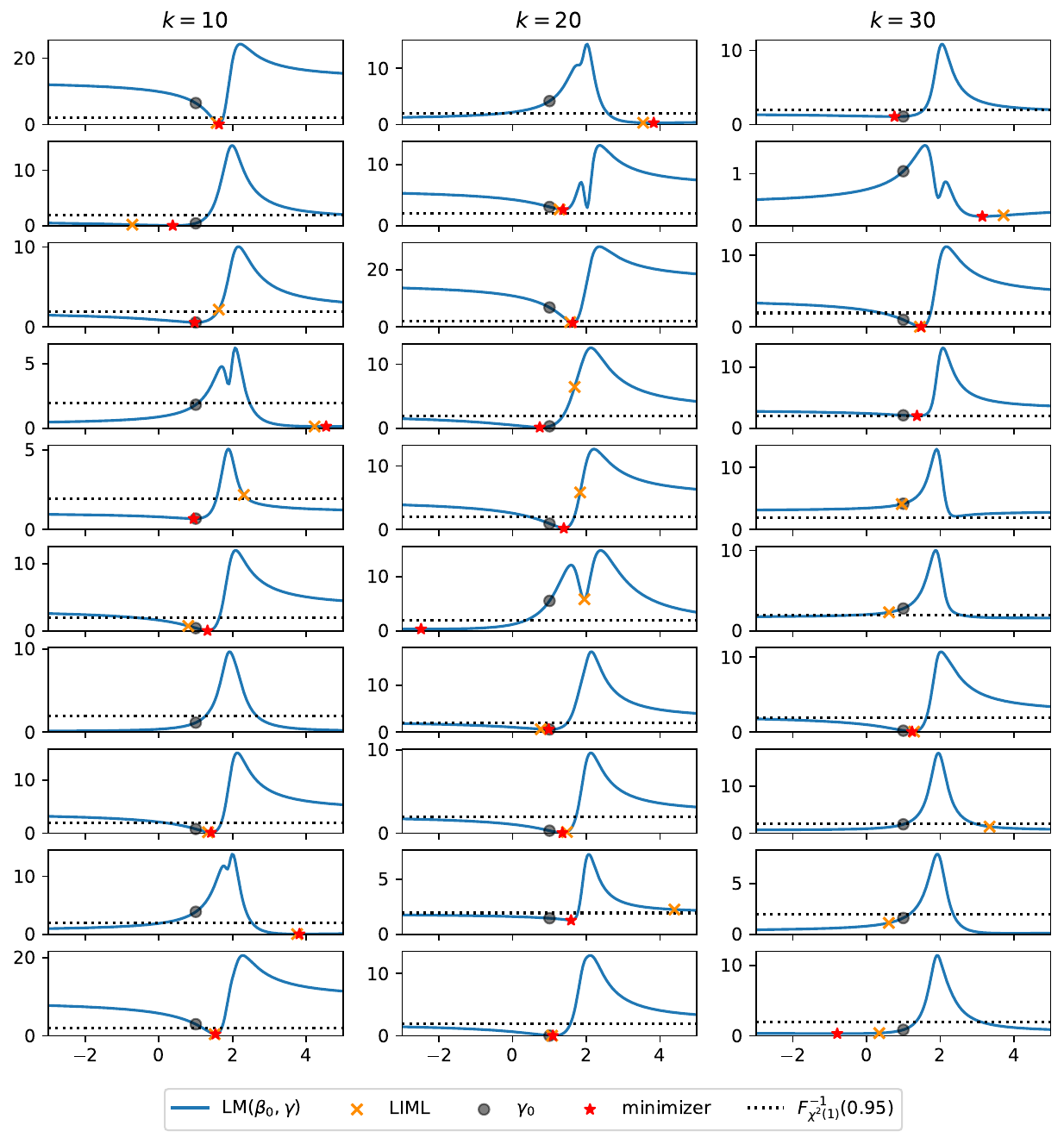}
    \caption{
        \label{fig:optimization}
        The Lagrange multiplier statistic for different $\gamma$ at $\beta_0$ for draws from \citeauthor{guggenberger2012asymptotic}'s \citeyearpar{guggenberger2012asymptotic} data-generating process.
        Rows correspond to different seeds and columns to different values of $k$.
        The LIML estimate using outcomes $y - X \beta_0$, endogenous covariates $W$, and instruments $Z$ is marked with an orange cross.
        The minimizer found with the BFGS algorithm is marked with a red star.
        The horizontal dashed line marks the 95\% quantile of the $\chi^2(1)$ distribution.
        The statistic at the minimizer exceeds the 95\% quantile only for the data in the fifth row and third column.
    }
    \vspace{-5cm}
\end{figure}

%% file: figures/figure_guggenberger12_qqplots.tex
\begin{figure}[!htb]
    \centering
    \includegraphics[width=0.95\textwidth]{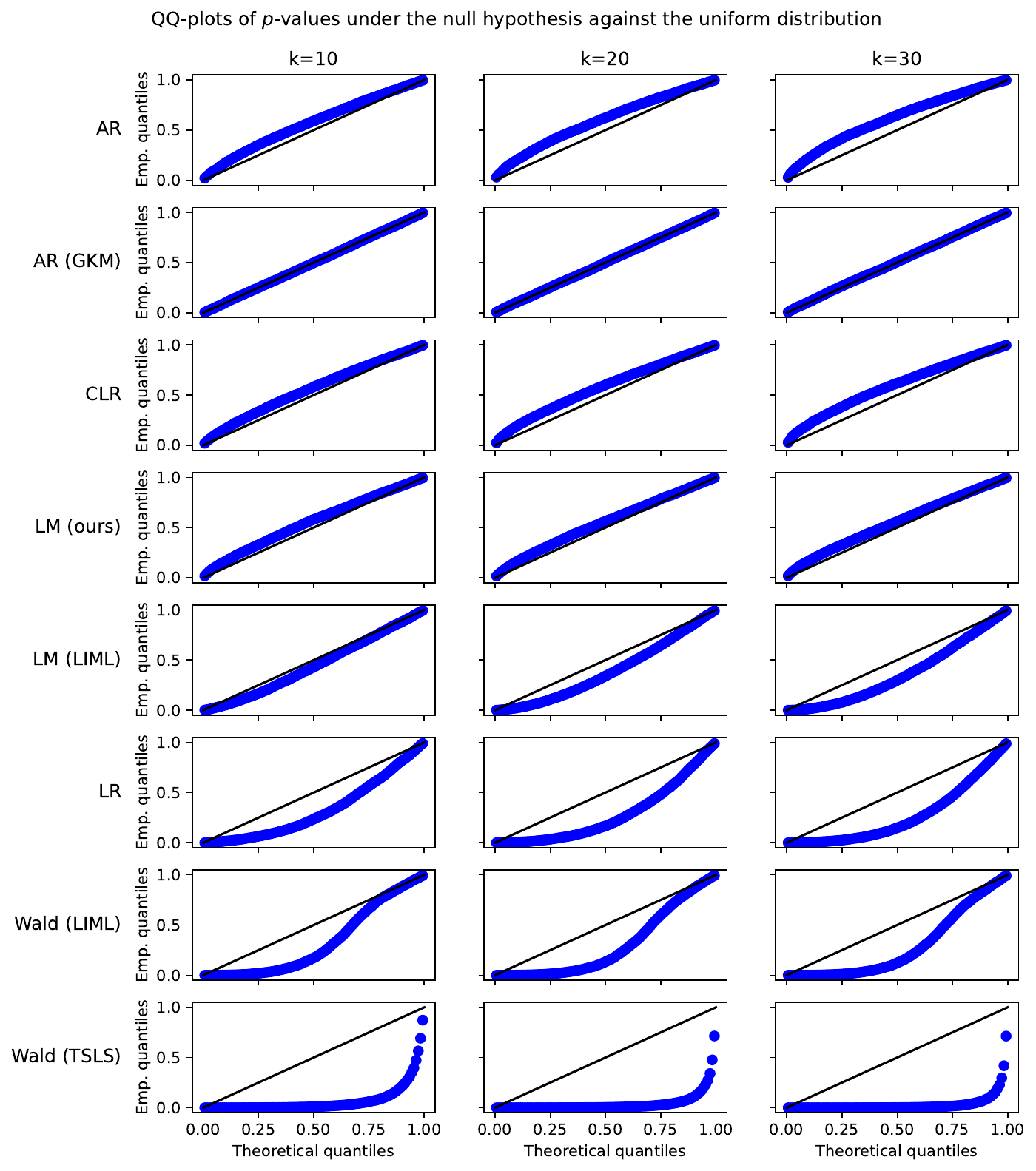}
    \caption{
        \label{fig:guggenberger12_qqplots}
        QQ plots comparing the empirical distribution of $p$-values of various subvector tests to the uniform distribution.
        Based on 10'000 simulations from the data-generating process proposed by \citet{guggenberger2012asymptotic} and described in \Cref{sec:numerical_analysis:guggenberger12_size}.
        Points above the diagonal indicate that the test is conservative.
        Points below the diagonal indicate that the test is size distorted.
        AR (GKM) uses the critical values of \citet{guggenberger2019more}.
    }
\end{figure}

%% file: figures/figure_guggenberger12_power_rho=0.8.tex
\begin{figure}[ht]
    \centering
    \includegraphics[width=0.8\textwidth]{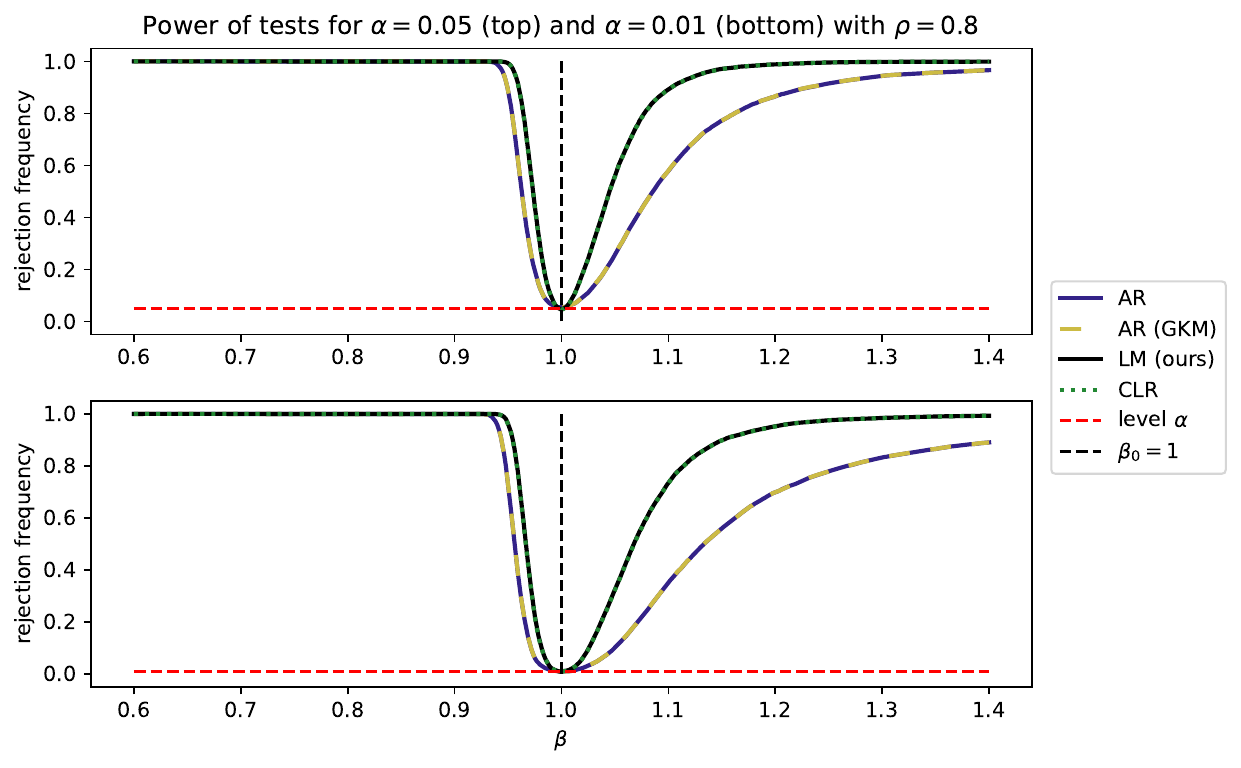}
    \caption{
        \label{fig:guggenberger12_power_rho=0.8}
        Power curves of various weak-instrument-robust subvector tests, based on 10'000 simulations from the data-generating process proposed by \citet{guggenberger2012asymptotic} with $n=1000, k=10$, but with $\Pi_W$ scaled such that $\sqrt{n} \| \Pi_W \| = 10$ and $\rho = \langle \Pi_W,  \Pi_X \rangle / ( \| \Pi_W \| \| \Pi_X \| ) = 0.8$.
        AR (GKM) uses the critical values of \citet{guggenberger2019more}.
    }
\end{figure}

%% file: figures/figure_guggenberger12_power_rho=0.99.tex
\begin{figure}[ht]
    \centering
    \includegraphics[width=0.99\textwidth]{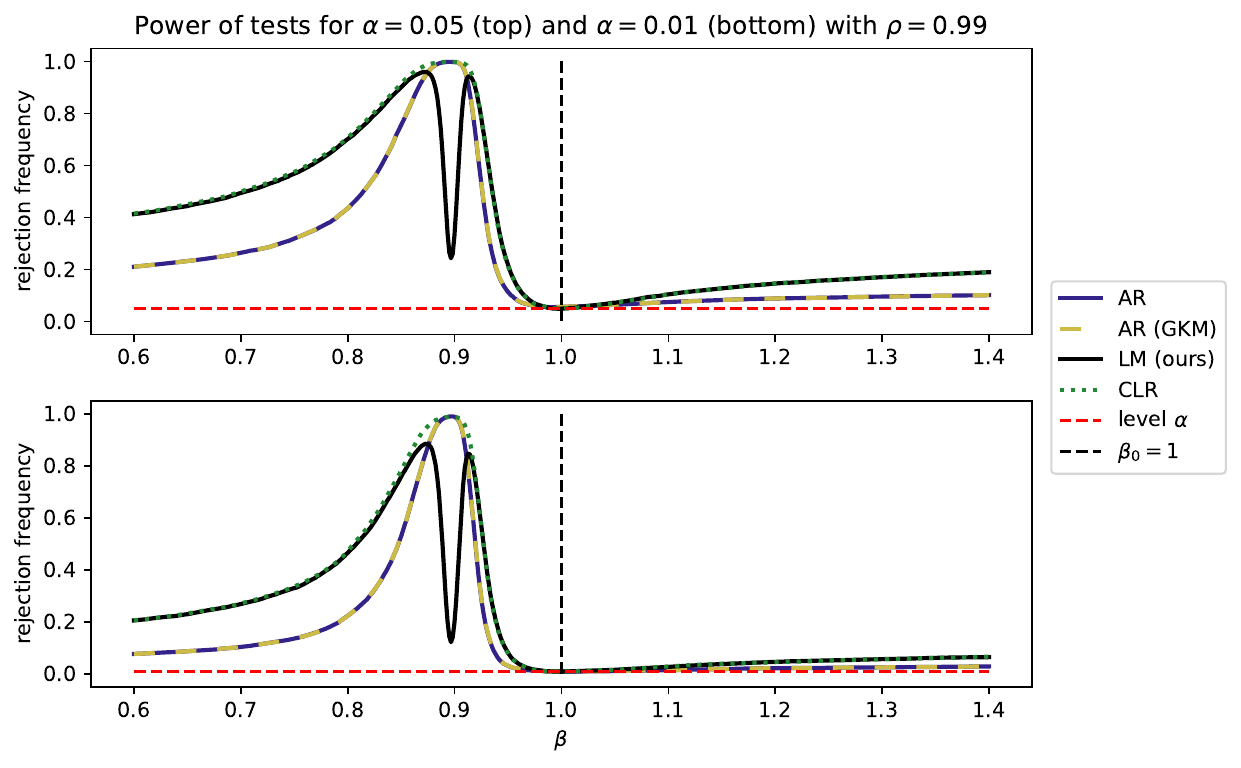}
    \caption{
        \label{fig:guggenberger12_power_rho=0.99}
        Power curves of various weak-instrument-robust subvector tests, based on 10'000 simulations from the data-generating process proposed by \citet{guggenberger2012asymptotic} with $n=1000, k=10$, but with $\Pi_W$ scaled such that $\sqrt{n} \| \Pi_W \| = 10$ and $\rho = \langle \Pi_W,  \Pi_X \rangle / ( \| \Pi_W \| \| \Pi_X \| ) = 0.99$.
        AR (GKM) uses the critical values of \citet{guggenberger2019more}.
    }
\end{figure}

%% file: figures/figure_kleibergen19_identity_k5_20.tex
\begin{figure}[p]
    \centering
    \includegraphics[width=0.95\textwidth]{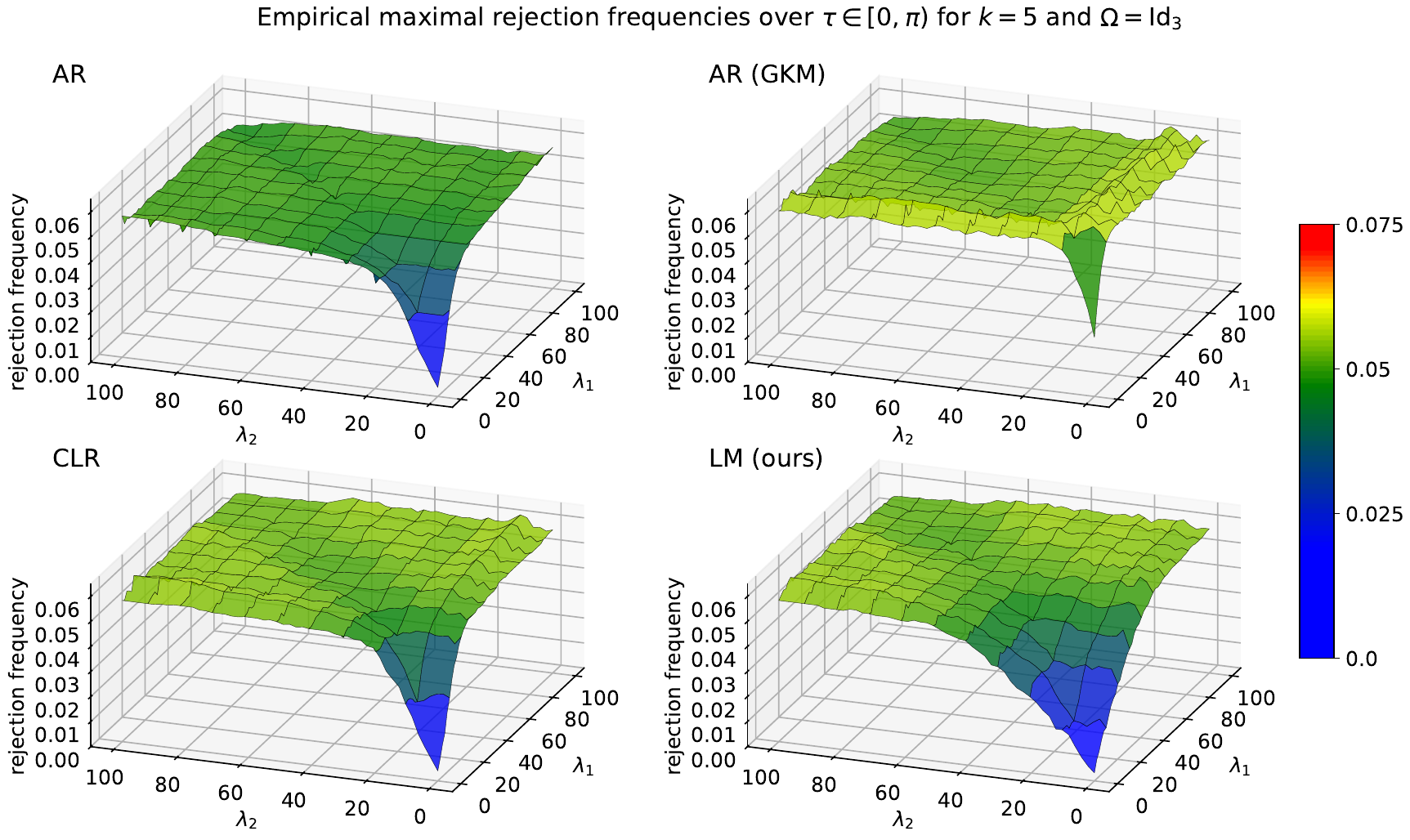}\\
    \vspace{0.5cm}
    \includegraphics[width=0.95\textwidth]{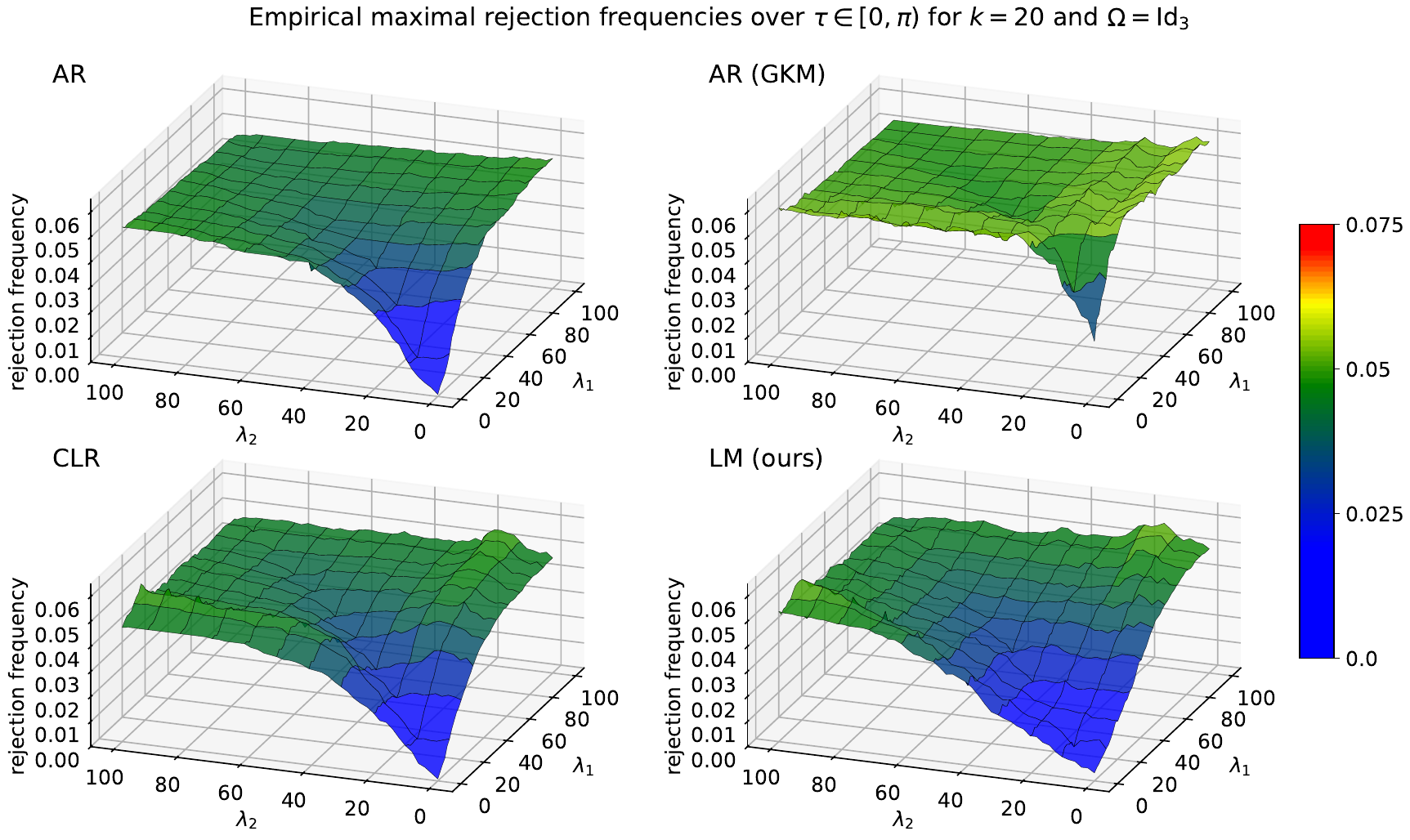}\\
    \vspace{0.2cm}
    \caption{
        \label{fig:kleibergen19_identity_k5_20}
        \figurekleibergencaption{$k=5$ (top) and $k=20$ (bottom)}{$\Omega = \Id_3$}
    }
\end{figure}

%% file: figures/figure_kleibergen19_guggenberger12_k5_20.tex
\begin{figure}[p]
    \centering
    \includegraphics[width=0.95\textwidth]{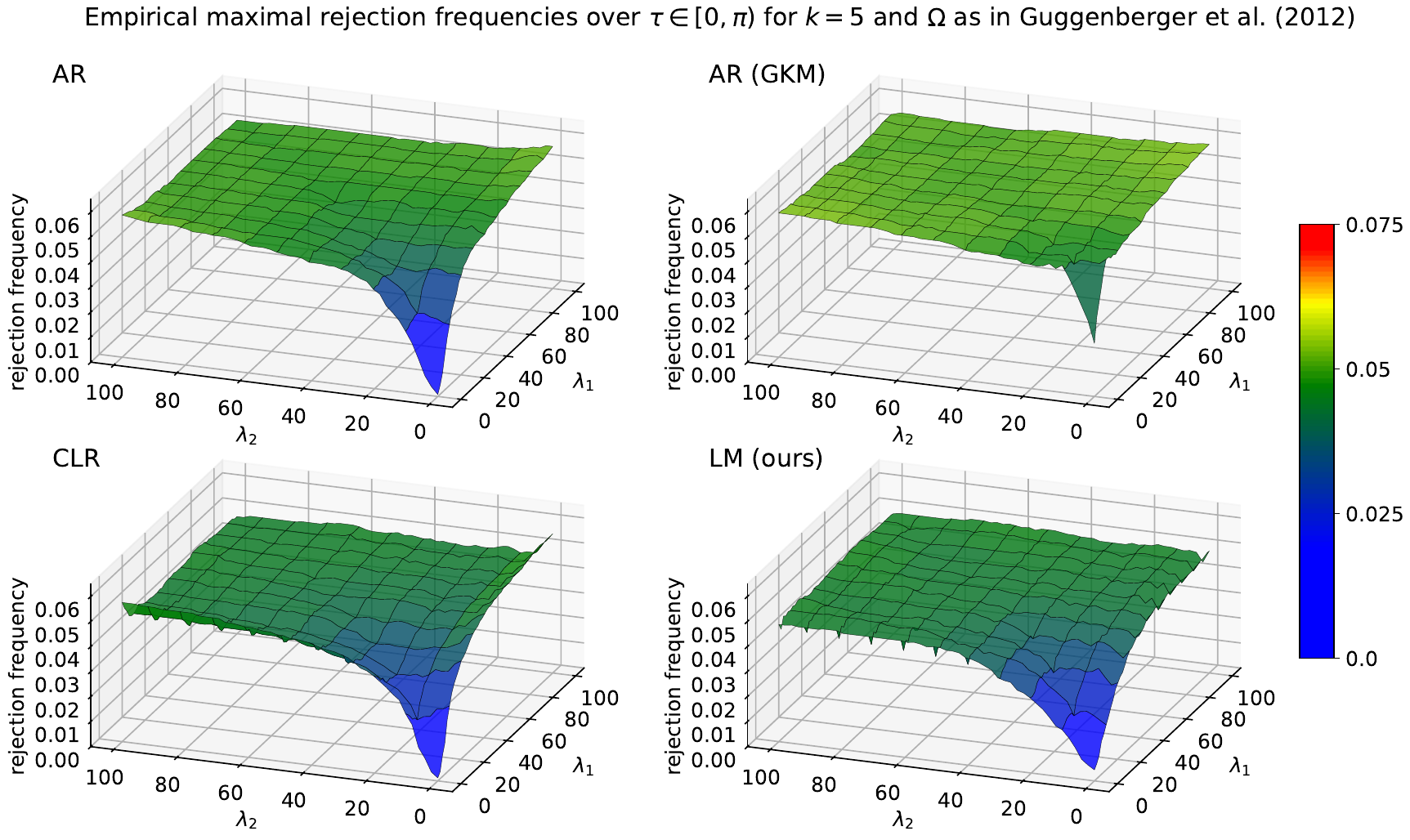}\\
    \vspace{0.5cm}
    \includegraphics[width=0.95\textwidth]{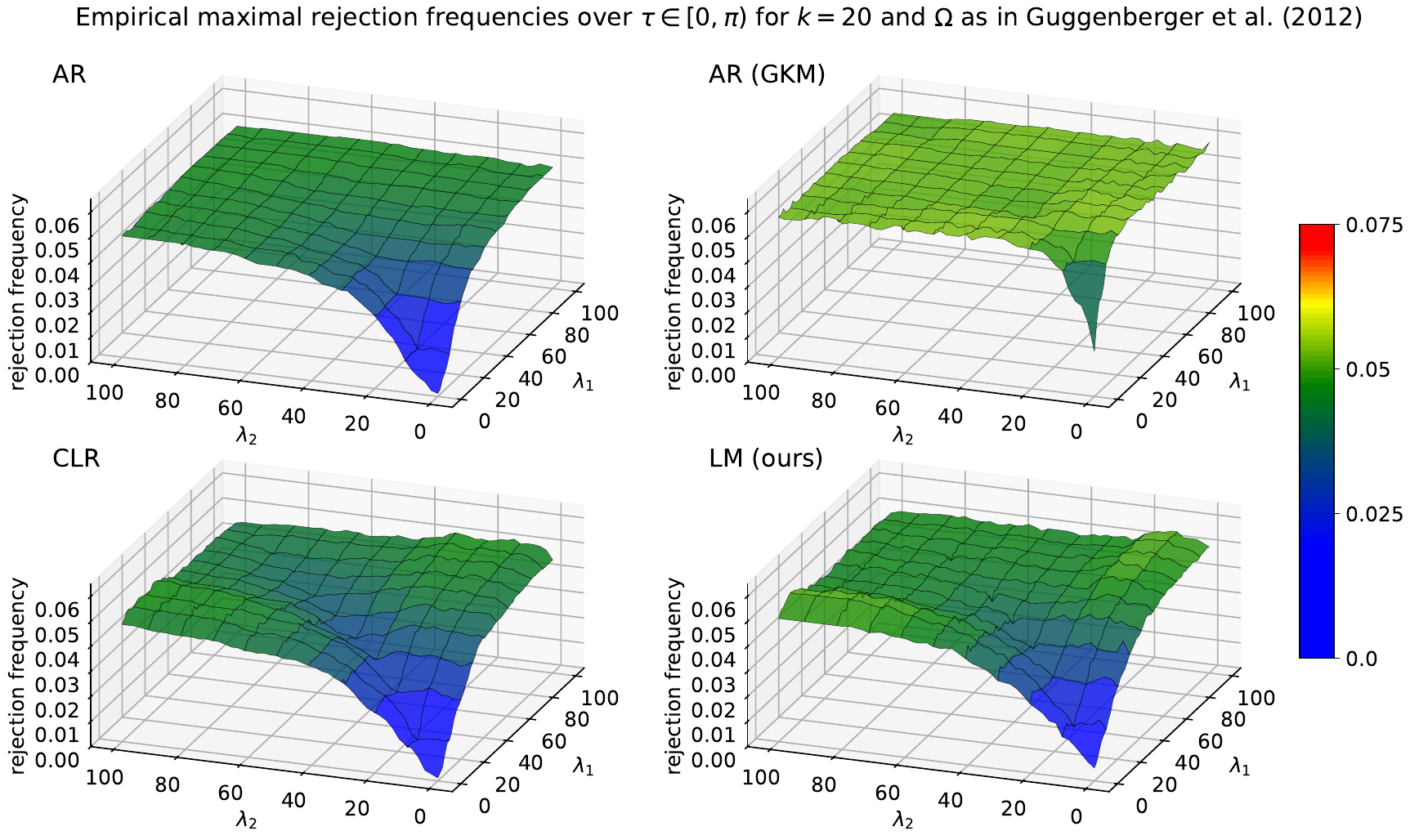}\\
    \vspace{0.2cm}
    \caption{
        \label{fig:kleibergen19_guggenberger12_k5_20}
        \figurekleibergencaption{$k=5$ (top) and $k=20$ (bottom)}{$\Omega$ as in \citep{guggenberger2012asymptotic}}
    }
\end{figure}

%% file: figures/figure_kleibergen19_guggenberger12_k100.tex
\begin{figure}[p]
    \centering
    \includegraphics[width=0.95\textwidth]{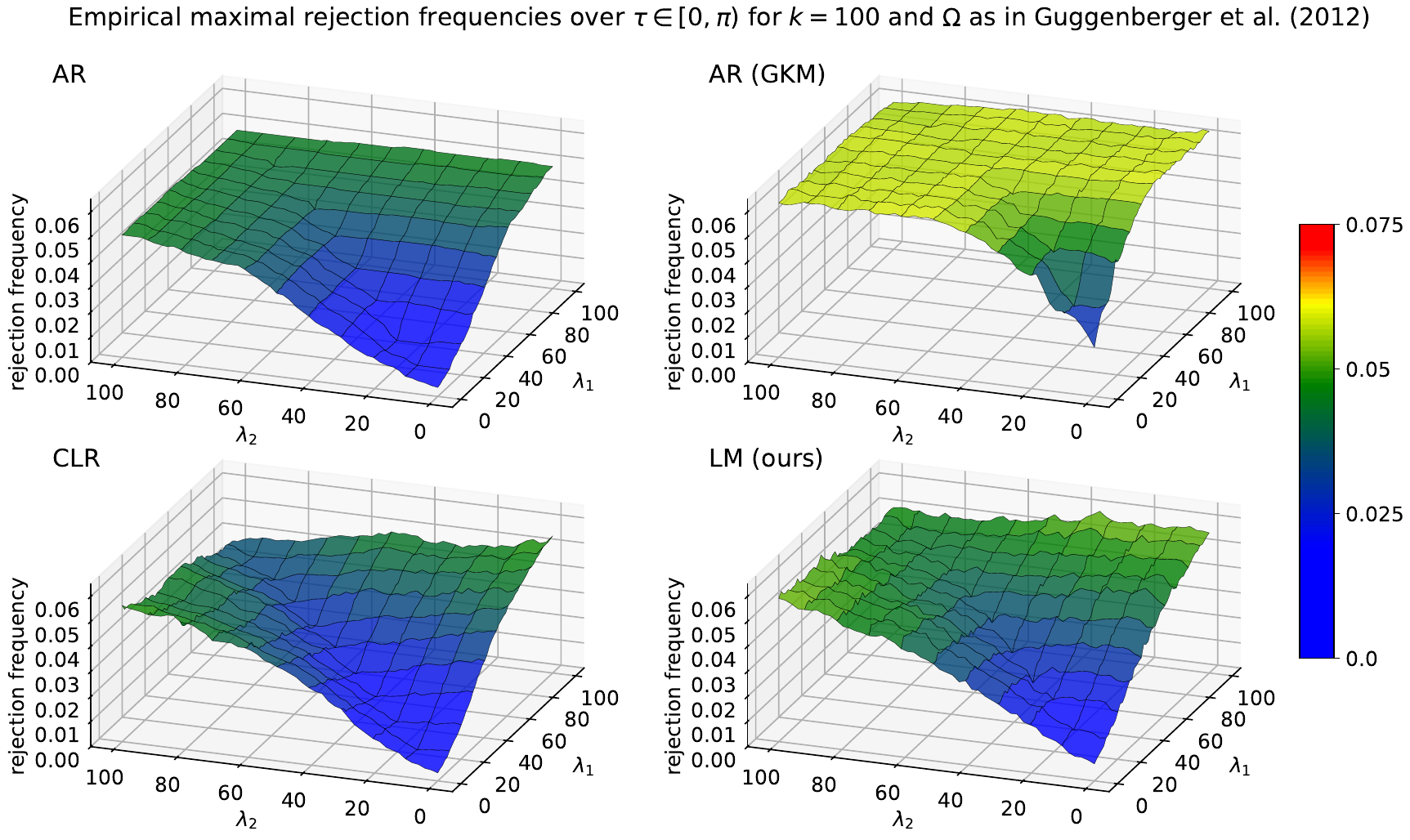}
    \caption{
        \label{fig:kleibergen19_guggenberger12_k100}
        \figurekleibergencaption{$k=100$}{$\Omega$ as in \citep{guggenberger2012asymptotic}}
    }
\end{figure}

%% file: figures/figure_card.tex
\begin{figure}[!htb]
    \centering
    \includegraphics[width=0.95\textwidth]{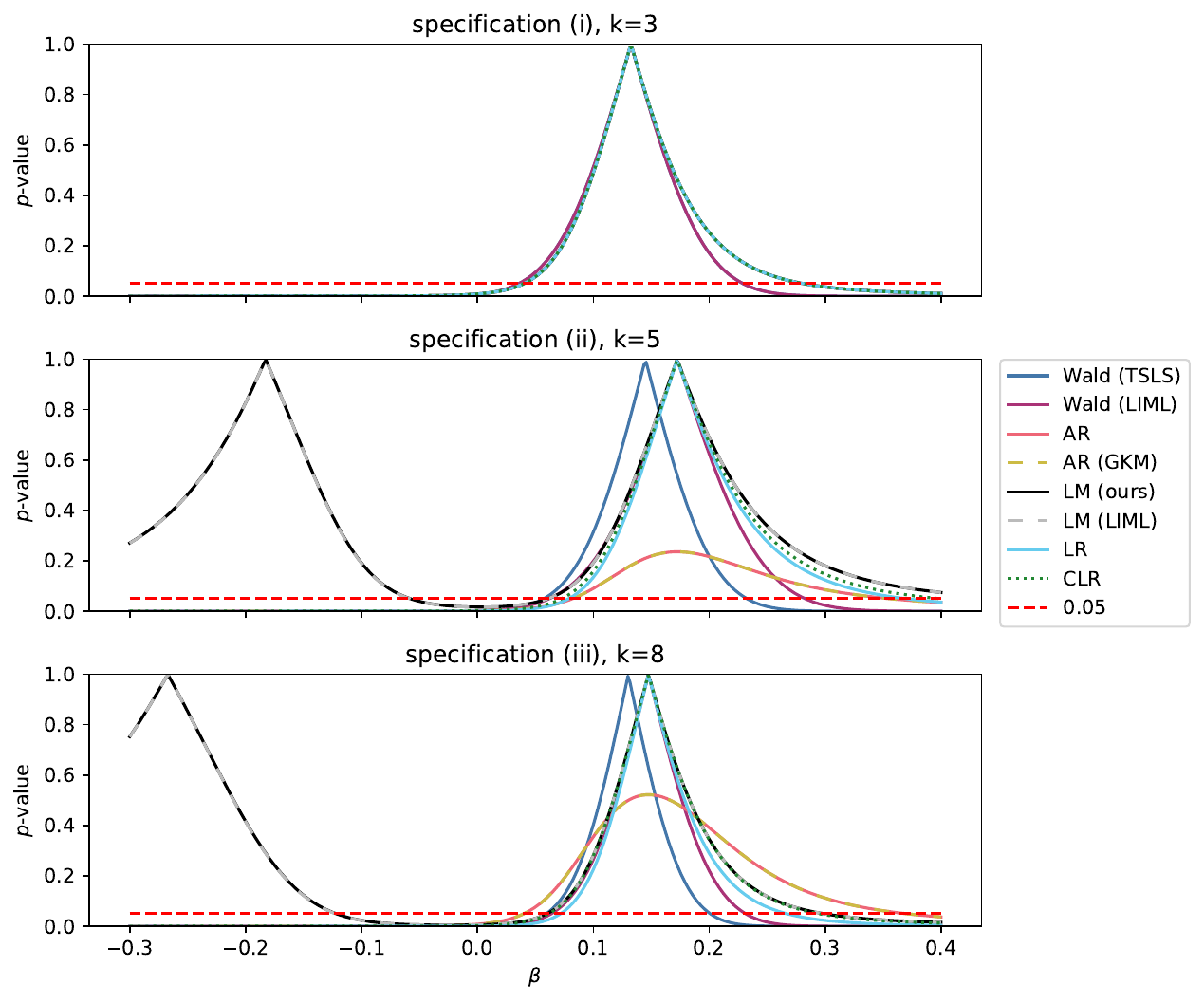}
    \caption{
        \label{fig:card_0}
        The $p$-values of the tests for varying values of $\beta$ for different specifications on the Card application.
        Corresponds to \Cref{tab:applications:card}.
    }
\end{figure}

%% file: figures/figure_tanaka.tex
\begin{figure}[!htb]
    \centering
    \includegraphics[width=0.95\textwidth]{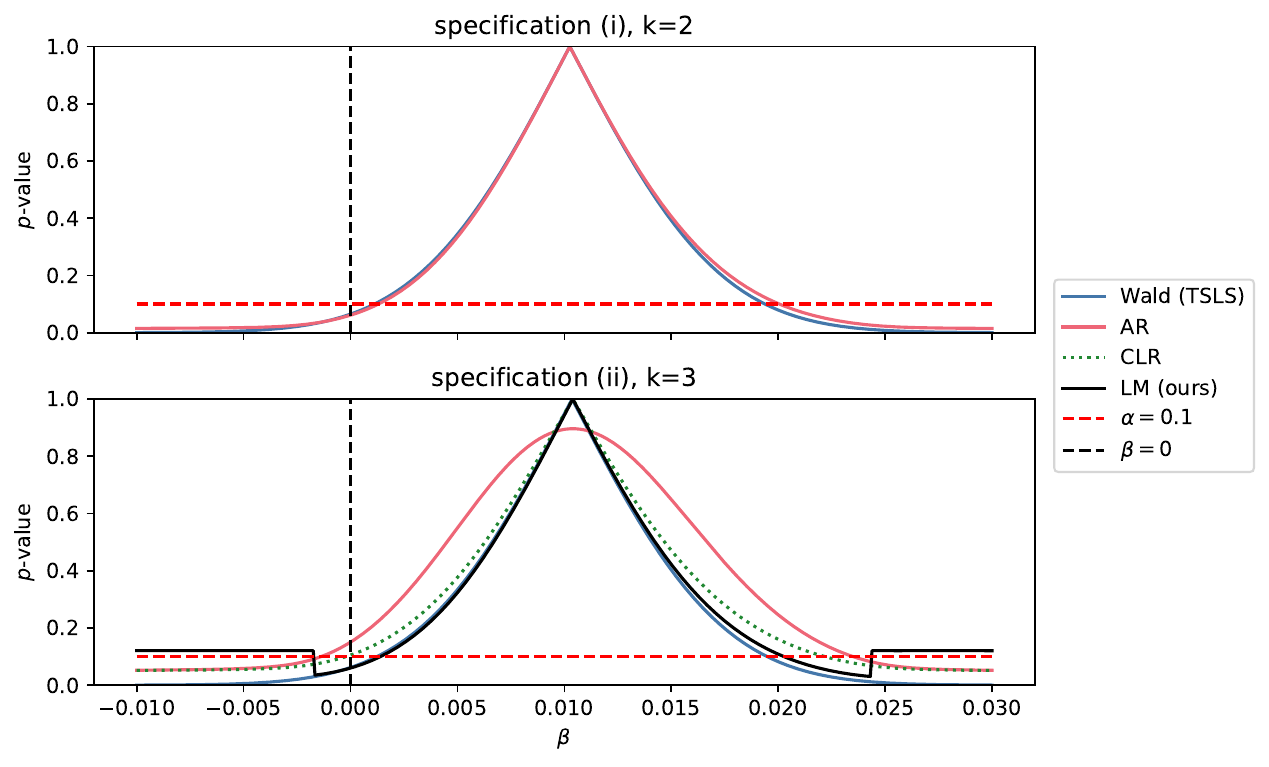}
    \caption{
        \label{fig:tanaka}
        The $p$-values of the causal effect of village mean income on risk preference from different tests for varying values of $\beta$. See also \Cref{tab:applications:tanaka}.
    }
\end{figure}

%% file: main.bbl
\begin{thebibliography}{}

\bibitem[\protect\citeauthoryear{Anderson}{Anderson}{1951}]{anderson1951estimating}
Anderson, T.~W. (1951).
\newblock Estimating linear restrictions on regression coefficients for
  multivariate normal distributions.
\newblock {\em The Annals of Mathematical Statistics\/}~{\em 22\/}(3),
  327--351.

\bibitem[\protect\citeauthoryear{Anderson and Rubin}{Anderson and
  Rubin}{1949}]{anderson1949estimation}
Anderson, T.~W. and H.~Rubin (1949).
\newblock Estimation of the parameters of a single equation in a complete
  system of stochastic equations.
\newblock {\em The Annals of Mathematical Statistics\/}~{\em 20\/}(1), 46--63.

\bibitem[\protect\citeauthoryear{Andrews, Moreira, and Stock}{Andrews
  et~al.}{2006}]{andrews2006optimal}
Andrews, D.~W., M.~J. Moreira, and J.~H. Stock (2006).
\newblock Optimal two-sided invariant similar tests for instrumental variables
  regression.
\newblock {\em Econometrica\/}~{\em 74\/}(3), 715--752.

\bibitem[\protect\citeauthoryear{Bentkus}{Bentkus}{2003}]{bentkus2003dependence}
Bentkus, V. (2003).
\newblock On the dependence of the {B}erry--{E}sseen bound on dimension.
\newblock {\em Journal of Statistical Planning and Inference\/}~{\em 113\/}(2),
  385--402.

\bibitem[\protect\citeauthoryear{Card}{Card}{1995}]{card1993using}
Card, D. (1995).
\newblock Using geographic variation in college proximity to estimate the
  return to schooling.
\newblock {\em Aspects of labour market behaviour: essays in honour of {J}ohn
  {V}anderkamp\/}, 201--222.

\bibitem[\protect\citeauthoryear{Chernozhukov, Chetverikov, and
  Kato}{Chernozhukov et~al.}{2017}]{chernozhukov2017central}
Chernozhukov, V., D.~Chetverikov, and K.~Kato (2017).
\newblock Central limit theorems and bootstrap in high dimensions.
\newblock {\em The Annals of Probability\/}~{\em 45\/}(4), 2309--2352.

\bibitem[\protect\citeauthoryear{Dufour and Taamouti}{Dufour and
  Taamouti}{2005}]{dufour2005projection}
Dufour, J.-M. and M.~Taamouti (2005).
\newblock Projection-based statistical inference in linear structural models
  with possibly weak instruments.
\newblock {\em Econometrica\/}~{\em 73\/}(4), 1351--1365.

\bibitem[\protect\citeauthoryear{Gallier}{Gallier}{2019}]{gallier2010schur}
Gallier, J. (2019).
\newblock The {Schur} complement and symmetric positive semidefinite (and
  definite) matrices.
\newblock {\em Penn Engineering\/}.

\bibitem[\protect\citeauthoryear{Guggenberger, Kleibergen, and
  Mavroeidis}{Guggenberger et~al.}{2019}]{guggenberger2019more}
Guggenberger, P., F.~Kleibergen, and S.~Mavroeidis (2019).
\newblock A more powerful subvector {Anderson} {Rubin} test in linear
  instrumental variables regression.
\newblock {\em Quantitative Economics\/}~{\em 10\/}(2), 487--526.

\bibitem[\protect\citeauthoryear{Guggenberger, Kleibergen, and
  Mavroeidis}{Guggenberger et~al.}{2024}]{guggenberger2024powerful}
Guggenberger, P., F.~Kleibergen, and S.~Mavroeidis (2024).
\newblock A powerful subvector {A}nderson--{R}ubin test in linear instrumental
  variables regression with conditional heteroskedasticity.
\newblock {\em Econometric Theory\/}~{\em 40\/}(5), 957--1002.

\bibitem[\protect\citeauthoryear{Guggenberger, Kleibergen, Mavroeidis, and
  Chen}{Guggenberger et~al.}{2012}]{guggenberger2012asymptotic}
Guggenberger, P., F.~Kleibergen, S.~Mavroeidis, and L.~Chen (2012).
\newblock On the asymptotic sizes of subset {Anderson--Rubin} and {Lagrange}
  multiplier tests in linear instrumental variables regression.
\newblock {\em Econometrica\/}~{\em 80\/}(6), 2649--2666.

\bibitem[\protect\citeauthoryear{Kleibergen}{Kleibergen}{2002}]{kleibergen2002pivotal}
Kleibergen, F. (2002).
\newblock Pivotal statistics for testing structural parameters in instrumental
  variables regression.
\newblock {\em Econometrica\/}~{\em 70\/}(5), 1781--1803.

\bibitem[\protect\citeauthoryear{Kleibergen}{Kleibergen}{2007}]{kleibergen2007generalizing}
Kleibergen, F. (2007).
\newblock Generalizing weak instrument robust {IV} statistics towards multiple
  parameters, unrestricted covariance matrices and identification statistics.
\newblock {\em Journal of Econometrics\/}~{\em 139\/}(1), 181--216.

\bibitem[\protect\citeauthoryear{Kleibergen}{Kleibergen}{2021}]{kleibergen2021efficient}
Kleibergen, F. (2021).
\newblock Efficient size correct subset inference in homoskedastic linear
  instrumental variables regression.
\newblock {\em Journal of Econometrics\/}~{\em 221\/}(1), 78--96.

\bibitem[\protect\citeauthoryear{Lee, McCrary, Moreira, and Porter}{Lee
  et~al.}{2022}]{lee2022valid}
Lee, D.~S., J.~McCrary, M.~J. Moreira, and J.~Porter (2022).
\newblock Valid t-ratio inference for {IV}.
\newblock {\em American Economic Review\/}~{\em 112\/}(10), 3260--3290.

\bibitem[\protect\citeauthoryear{Londschien}{Londschien}{2025}]{londschien2025overview}
Londschien, M. (2025).
\newblock A statistician's guide to weak-instrument-robust inference in
  instrumental variables regression with illustrations in {Python}.
\newblock {\em arXiv preprint arXiv:2508.12474\/}.

\bibitem[\protect\citeauthoryear{Moreira}{Moreira}{2003}]{moreira2003conditional}
Moreira, M.~J. (2003).
\newblock A conditional likelihood ratio test for structural models.
\newblock {\em Econometrica\/}~{\em 71\/}(4), 1027--1048.

\bibitem[\protect\citeauthoryear{Rohde}{Rohde}{1965}]{rohde1965generalized}
Rohde, C.~A. (1965).
\newblock Generalized inverses of partitioned matrices.
\newblock {\em Journal of the Society for Industrial and Applied
  Mathematics\/}~{\em 13\/}(4), 1033--1035.

\bibitem[\protect\citeauthoryear{Staiger and Stock}{Staiger and
  Stock}{1997}]{staiger1997instrumental}
Staiger, D.~O. and J.~H. Stock (1997).
\newblock Instrumental variables regression with weak instruments.
\newblock {\em Econometrica\/}~{\em 65\/}(3), 557--586.

\bibitem[\protect\citeauthoryear{Stock and Yogo}{Stock and
  Yogo}{2005}]{stock2002testing}
Stock, J.~H. and M.~Yogo (2005).
\newblock Testing for weak instruments in linear {IV} regression.
\newblock {\em Identification and Inference for Econometric Models: Essays in
  Honor of Thomas Rothenberg\/}, 80--108.

\bibitem[\protect\citeauthoryear{Tanaka, Camerer, and Nguyen}{Tanaka
  et~al.}{2010}]{tanaka2010risk}
Tanaka, T., C.~F. Camerer, and Q.~Nguyen (2010).
\newblock Risk and time preferences: Linking experimental and household survey
  data from vietnam.
\newblock {\em American Economic Review\/}~{\em 100\/}(1), 557--571.

\bibitem[\protect\citeauthoryear{Thams, S{\o}ndergaard, Weichwald, and
  Peters}{Thams et~al.}{2024}]{thams2022identifying}
Thams, N., R.~S{\o}ndergaard, S.~Weichwald, and J.~Peters (2024).
\newblock Identifying causal effects using instrumental time series: Nuisance
  {IV} and correcting for the past.
\newblock {\em Journal of Machine Learning Research\/}~{\em 25}, 1--51.

\end{thebibliography}
